%% file: paper.tex
\documentclass{article}
\usepackage{jlcode}
\RequirePackage{amsmath}
\RequirePackage{amsfonts}
\RequirePackage{amssymb}
\RequirePackage{mathtools}  
\RequirePackage{mathrsfs}   
\RequirePackage{siunitx}
\RequirePackage{bbm}        
\RequirePackage{nicefrac}   
\RequirePackage[hyphens]{url}
\RequirePackage{adjustbox}
\RequirePackage{booktabs}
\RequirePackage{multirow}
\RequirePackage{tabularx}
\usepackage{arxiv}
\usepackage{kthcolors}
\usepackage{algorithm}
\usepackage{algpseudocode}

\algnewcommand\algorithmicinput{\hspace*{1.0em}\textbf{input:}}
\algnewcommand\Input{\item[\algorithmicinput]}
\algnewcommand\algorithmicoutput{\hspace*{1.0em}\textbf{output:}}
\algnewcommand\Output{\item[\algorithmicoutput]}

\DeclarePairedDelimiter{\parentheses}{(}{)}
\DeclarePairedDelimiter{\brackets}{[}{]}
\DeclarePairedDelimiter{\braces}{\lbrace}{\rbrace}
\DeclarePairedDelimiter{\verts}{\lvert}{\rvert}

\DeclarePairedDelimiterX{\Set}[1]{\lbrace}{\rbrace}{
  
  #1
}

\providecommand{\expect}[2][{}k]{\ensuremath{\operatorname{\mathbb{E}}_{#1}\brackets*{#2}}}
\providecommand{\var}[2][{}k]{\ensuremath{\operatorname{Var}}_{#1}\brackets*{#2}}
\providecommand{\abs}[1]{\ensuremath{\verts*{#1}}}

\providecommand{\argmin}[2][x]{\ensuremath{\arg\,\min_{#1}\braces*{#2}}}

\providecommand{\sbj}{\ensuremath\text{subject to}}
\providecommand{\st}{\ensuremath\,\text{s.t.}\,}

\DeclareMathOperator*{\minimize}{minimize}
\DeclareMathOperator*{\maximize}{maximize}

\usepackage{pgfplots}
\usepackage{tikz}

\pgfplotsset{compat=1.14}

\usetikzlibrary{arrows.meta,intersections,fit,math,shapes,calc,positioning}

\begin{document}

\title{Large-scale hydropower models in StochasticPrograms.jl}

\author{Martin Biel \\
  Division of Decision and Control Systems\\
  School of EECS, KTH Royal Institute of Technology\\
  SE-100 44 Stockholm, Sweden\\
  \texttt{mbiel@kth.se}}

\maketitle

\begin{abstract}
  \sloppy
  We present three large-scale hydropower planning models implemented in our open-source software framework \jlinl{StochasticPrograms.jl} developed using the Julia programming langugage. The framework provides an expressive syntax for formulating stochastic programming models and has distributed capabilities that can handle large-scale instances. The three models describe different case studies of the hydroelectric power plants in the Swedish river Skellefteälven. The models are two-stage stochastic programs with sampled scenarios that describe uncertain electricity prices and local water inflows. The first model is a day-ahead planning problem that concerns how to determine optimal order strategies in a day-ahead energy market. We pose this problem from the perspective of a hydropower producer, who participates in the Nordic day-ahead market and operates in the Swedish river Skellefteälven. We implement the day-ahead model using our computational tools and then solve large-scale instances of the problem in a distributed environment. A statistically significant value of running stochastic planning is obtained using a sample-based algorithm. Next, we consider a variation of the day-ahead problem that includes preventive maintenance scheduling. We show how intricate coordination between the submitteed market orders and the maintenance schedule results in a larger value of the stochastic solution than the day-ahead problem. The final model is a capacity expansion problem with a long planning horizon. The same methodology is applied as when solving the first two hydropower problems. However, the planning horizon is considerably longer, from one year up to 20 years compared to a 24 hour horizon. We note that the relative significance of the value of the stochastic solution is much greater when comparing to the extra profits incurred from the capacity expansion instead of the total profit.
\end{abstract}

\section{Introduction}
\label{sec:introduction}

Engineering problems often involve making decisions under uncertainty. In particular, hydropower operations are associated with various random elements that make it difficult to operate optimally. For example, the future electricity price is not known when scheduling the next-day production, which could incur a loss of profit if not considered. While power production from renewable sources, such as solar and wind, is constantly expanding, the production is irregular in nature. Therefore, a large increase of renewable power production can lead to large imbalances in the power grid. Because hydropower production can be delayed by storing water in reservoirs it is possible to coordinate the power production to balance the power grid. This is only possible if the random renewable production is taken into account when planning the hydropower production. Seasonal planning of reservoir contents is also associated with random demand in load from end users as well as random water inflow. Hydropower operations in the Nordic regions are also affected by snow melt. Therefore, if the reservoir levels are high during the spring flood, there is an high risk of spillage. This results in lost opportunities for power production. In all examples, uncertainty can be taken into account during planning by formulating and solving a mathematical model.

In this paper, we present three large-scale hydropower models implemented in the open-source software framework \jlinl{StochasticPrograms.jl}. Stochastic programming is a well-established approach in power systems modeling~\cite{Fleten2007,Groewe-Kuska2005,petra_real-time_2014,LoB:16}. We showcase how our framework can be used to formulate complex energy models with uncertain elements. Moreover, we utilize the parallel capabilities and specialized algorithms in the framework to handle large-scale instances of the energy models. We will present three case studies on large-scale hydropower modeling to illustrate the potential of our framework.

\section{Preliminaries}
\label{sec:preliminaries}

We first provide some necessary background information to keep this work self-contained. We list some preliminary stochastic programming results that will be used throughout the paper and briefly introduce our software framework. We also describe the deterministic and uncertain data parameters that are required to create the hydropower models.

\subsection{Stochastic programming}
\label{sec:stoch-progr}

A two-stage linear stochastic program is given by
\begin{equation} \label{eq:linearsp}
  \begin{aligned}
    \minimize_{\mathclap{x \in \mathbb{R}^n}} & \quad c^T x + \expect[\xi]{Q(x,\xi(\omega))} \\
    \sbj & \quad Ax = b \\
    & \quad x \geq 0,
  \end{aligned}
\end{equation}
where
\begin{equation} \label{eq:subprob}
  \begin{aligned}
    Q(x,\xi(\omega)) = \min_{\mathclap{y \in \mathbb{R}^m}} & \quad q_{\omega}^T y \\
    \st & \quad T_{\omega}x + Wy = h_{\omega} \\
    & \quad y \geq 0.
  \end{aligned}
\end{equation}
The formulation seeks the first-stage decision $x$ that is optimal in expectation over a set of future uncertain scenarios. A specific scenario $\omega$ is drawn from the sample space $\Omega$. The random variable
\begin{equation*}
  \xi(\omega) = \begin{pmatrix}
    q_{\omega} \\
    h_{\omega} \\
    T^1_{\omega} \\
    \vdots \\
    T^q_{\omega}
  \end{pmatrix}
\end{equation*}
then parameterizes the second-stage subproblem~\eqref{eq:subprob} where corrective actions $y$ can be taken to mitigate the first-stage decision. The \emph{value of the recourse problem}, or the VRP, is the optimal value of~\eqref{eq:linearsp}.

If $\Omega$ is finite, say with $N$ scenarios of probability $\pi_s,\, s \in \Set{1,\dots,N}$ respectively, then we can represent~\eqref{eq:linearsp} compactly as
\begin{equation} \label{eq:finitesp}
  \begin{aligned}
    \minimize_{\mathclap{x \in \mathbb{R}^n, y_s \in \mathbb{R}^m}} & \quad c^T x + \sum_{s = 1}^{N} \pi_s q_s^T y_s & \\
    \sbj & \quad Ax = b & \\
    & \quad T_s x + W y_s = h_s, \quad &&s = 1,\dots,N \\
    & \quad x \geq 0, \, y_s \geq 0, \quad &&s = 1,\dots,N.
  \end{aligned}
\end{equation}
For small $N$, it is viable to solve this problem with standard solvers. For large $N$, decomposition approaches~\cite{van_slyke_l-shaped_1969, Rockafellar1991} are required. If $\Omega$ is not finite, the stochastic program~\eqref{eq:linearsp} is exactly computable only under certain assumptions. However, it is possible to formulate computationally tractable approximations of~\eqref{eq:linearsp} using the finite form~\eqref{eq:finitesp}. The most common approximation technique is the \textit{sample average approximation} (SAA)~\cite{saa}. Assume that we sample $N$ scenarios $\omega_s,\; s = 1,\dots,N$ independently from $\Omega$ with equal probability. These scenarios now constitute a finite sample space $\tilde{\Omega}$ with the measure
\begin{equation*}
  \tilde{\pi}(\omega_s) = \frac{1}{N},\quad s = 1,\dots,N,
\end{equation*}
and we can use them to create a sampled model in finite extensive form~\eqref{eq:finitesp}. An optimal solution to this sampled model approximates the optimal solution to~\eqref{eq:linearsp} in the sense that the empirical average second-stage cost $\frac{1}{N}\sum_{s = 1}^{N}q_s^T \hat{y}_s,$ where $\hat{y}_s = \argmin[y \in \mathbb{R}^m]{Q(x, \xi(\omega_s))}$, converges pointwise with probability $1$ to $\expect[\xi]{Q(x,\xi(\omega))}$ as $N$ goes to infinity~\cite{saaconvergence}. In practice, we do not reach this asymptotic limit, and instead, rely on statistical approaches. Such methods are based on the following central limit result~\cite{saadist}. As $N$ goes to infinity, it holds that $\sqrt{N}(V_N - VRP) \to_d \mathcal{N}(0,\var[\xi]{Q(\hat{x},\xi)})$, where
  \begin{equation*}
    V_N = \min_{x \in \mathcal{X}}{c^Tx + \sum_{s = 1}^{N} Q(x, \xi_s)}.
  \end{equation*}
This result provides a basis for calculating confidence intervals around the VRP of~\eqref{eq:linearsp}, as described next.

We follow the methodology developed in~\cite{saa}. First, given any $x \in \mathbb{R}^n$, consider
\begin{equation} \label{eq:saaeval}
  V_N(x) = c^T x + \frac{1}{N}\sum_{s = 1}^{N}Q(x, \xi_s).
\end{equation}
Evaluation of the subproblems is cheap compared to solving~\eqref{eq:finitesp}; so, we can use large values of $N$. Now, we solve $T$ sampled batches of the above problem, using $\braces{\mathbf{\xi}^i}_{i = 1}^{T}$ where all $\mathbf{\xi}^i = \braces{\xi_s^i}_{s = 1}^{N}$ are $N$ IID samples of $\xi$, and construct the unbiased estimator
\begin{equation*}
  U_T = \frac{1}{T}\sum_{i = 1}^{T} V_N^i = \frac{1}{T} \sum_{i = 1}^{T} \parentheses*{c^T x + \frac{1}{N}\sum_{s = 1}^{N}Q(x, \xi^i_s)},
\end{equation*}
which estimates $V(x)$. Using the $\alpha$-critical value of the $t$-distribution with $T-1$ degrees of freedom, a $(1-\alpha)$ confidence interval around $U_{T}$ is given by
\begin{equation*}
  \brackets*{U_{T} - \frac{t_{\alpha/2,T-1}\sigma_{T}^2}{\sqrt{T}}, U_{T} + \frac{t_{\alpha/2,T-1}\sigma_{T}^2}{\sqrt{T}}},
\end{equation*}
where the sample variance is given by
\begin{equation*}
  \sigma_{T}^2 = \frac{1}{T-1}\sum_{i = 1}^{T} \parentheses*{V^i_N(x) - U_T}^2.
\end{equation*}

A lower bound on the gap between $V(x)$ and the VRP of~\eqref{eq:linearsp} is computed by first solving $N$-sized sampled instances of the form
\begin{equation} \label{eq:smallsaa}
\begin{aligned}
 \hat{V}^i_N = \min_{x \in \mathbb{R}^n} & \quad c^T x + \frac{1}{N}\sum_{s = 1}^N Q(x, \xi^i_s) \\
 \st & \quad Ax = b \\
 & \quad x \geq 0
\end{aligned}
\end{equation}
for $M$ sampled IID batches $\braces{\mathbf{\xi}^i}_{i = 1}^{N}$. Now,
\begin{equation*}
  VRP = \min_{x \in \mathcal{X}_r} c^Tx + \expect[\xi]{Q(x,\xi(\omega))} = \min_{x \in \mathcal{X}_r} \expect[\xi_s^i]{c^Tx + \frac{1}{N}\sum_{s = 1}^N Q(x, \xi^i_s)} \geq \expect[\xi^i_s]{\hat{V}^i_N}.
\end{equation*}
Consequently,
\begin{equation*}
  VRP \geq \expect[\xi_s^i]{\hat{V}^i_N},
\end{equation*}
and an estimate of this lower bound can be computed by
\begin{equation*}
  L_{N, M} = \frac{1}{M}\sum_{i = 1}^{M} \hat{V}^i_N.
\end{equation*}
Because this involves solving stochastic programs, the values of $M$ or $N$ cannot be too large for the procedure to be computationally tractable on a single node. In a distributed environment, we can employ parallel solver strategies since the procedure trivially parallelizes over $M$. An approximate $(1-\alpha)$ confidence interval around the lower bound is then given by
\begin{equation*}
\brackets*{L_{N,M} - \frac{t_{\alpha/2,M-1}\sigma_{N,M}^2}{\sqrt{M}}, L_{N,M} + \frac{t_{\alpha/2,M-1}\sigma_{N,M}^2}{\sqrt{M}}},
\end{equation*}
where the sample variance is given by
\begin{equation*}
  \sigma_{N,M}^2 = \frac{1}{M-1}\sum_{i = 1}^{M} \parentheses*{\hat{V}^i_N - L_{N,M}}^2.
\end{equation*}
The two bounds are now combined to form a $(1-2\alpha)$ confidence interval around the gap between $V(x)$ and the VRP of~\eqref{eq:linearsp}:
\begin{equation*}
  \brackets*{0, U_{T} - L_{N,M} + \frac{t_{\alpha/2,T-1}\sigma_{T}^2}{\sqrt{T}} + \frac{t_{\alpha/2,M-1}\sigma_{N,M}^2}{\sqrt{M}}}.
\end{equation*}

If we acquire a candidate decision by solving a single sampled instance, then we can use the above procedures to calculate a confidence interval around the VRP of~\eqref{eq:linearsp}. Specifically, assume that $\hat{x}_N$ is the optimizer of some sampled instance of size $N$. A confidence interval around the VRP of~\eqref{eq:linearsp} is then given by
\begin{equation} \label{eq:saaconfidence}
  \brackets*{L_{VRP},U_{VRP}} = \brackets*{L_{N,M} - \frac{t_{\alpha/2,M-1}\sigma_{N,M}^2}{\sqrt{M}}, \hat{U}_{N,T} + \frac{t_{\alpha/2,T-1}\hat{\sigma}_{N,T}^2}{\sqrt{T}}},
\end{equation}
where
\begin{equation*}
  \hat{U}_{N,T} = \frac{1}{T} \sum_{i = 1}^{T} \parentheses*{c^T \hat{x}_N + \frac{1}{N}\sum_{s = 1}^{N}Q(\hat{x}_N, \xi^i_s)},
\end{equation*}
and
\begin{equation*}
  \hat{\sigma}_{N,T}^2 = \frac{1}{T-1}\sum_{i = 1}^{T} \parentheses*{\parentheses*{c^T \hat{x}_N + \frac{1}{N}\sum_{s = 1}^{N}Q(\hat{x}_N, \xi^i_s)} - \hat{U}_{N,T}}^2.
\end{equation*}
The remaining quantities are calculated as before. Because the sampled solution $\hat{x}_N$ converges to the optimizer $\hat{x}$ as $N \to \infty$, it follows that $\hat{U}_{N,T}$ decreases with $N$ while $L_{N,M}$ increases with $N$. Hence, the length of the resulting confidence interval~\eqref{eq:saaconfidence} will decrease with $N$. We can therefore employ an iterative procedure, where $N$ is increased until the length of the confidence interval decreases to some desired relative tolerance. We would then report the resulting confidence interval $\brackets{L_{VRP},U_{VRP}}$ around the VRP of~\eqref{eq:linearsp}.

After solving~\eqref{eq:linearsp}, either exactly for finite problems or approximately using SAA, we can compute a classical measures of stochastic performance: The \emph{value of the stochastic solution}. This quantities is well defined for finite models~\eqref{eq:finitesp}. When $\Omega$ is infinite, we again employ the statistical approach as outlined above. Given
\begin{equation} \label{eq:expectxi}
  \bar{\xi} = \expect[\xi]{\xi(\omega)}
\end{equation}
the \emph{expected value decision} $\bar{x}$ associated with~\eqref{eq:linearsp} is given by the solution to
\begin{equation} \label{eq:ev}
  \begin{aligned}
    \minimize_{x \in \mathbb{R}^n} & \quad c^T x + Q(x,\bar{\xi}) \\
    \sbj & \quad Ax = b \\
    & \quad x \geq 0.
  \end{aligned}
\end{equation}
This problem is known as the \emph{expected value problem}. The \emph{expected result of the expected value decision}, or the EEV, is then given by
\begin{equation} \label{eq:EEV}
  EEV = c^T \bar{x} + \expect[\xi]{Q(\bar{x},\xi(\omega))}.
\end{equation}
Now, the \emph{value of the stochastic solution}, or the VSS, is given by
\begin{equation} \label{eq:VSS}
  VSS = VRP - EEV.
\end{equation}
The VSS measures the expected loss of ignoring the uncertainty in the problem. It indicates if the second stage is sensitive to the stochastic data and if there is any value in considering a stochastic formulation.

When $\Omega$ is infinite, we first use SAA to determine a confidence interval $\brackets{L_{VRP},U_{VRP}}$ around the VRP. This calculation involves solving sampled instances of size $N$. We then sample $N$ scenarios from $\Omega$ to obtain $\braces{\xi_i}_{i = 1}^{N}$, calculate the expected outcome
\begin{equation*}
  \bar{\xi} = \frac{1}{N}\sum_{i = 1}^{N}\xi_i,
\end{equation*}
and determine an expected value decision $\bar{x}$ according to~\eqref{eq:ev}. An EEV approximation is then obtained by
\begin{equation*}
  EEV = c^T \bar{x} + \frac{1}{\bar{N}}\sum_{i = 1}^{\bar{N}}Q(\bar{x}, \xi_i).
\end{equation*}
where $\braces{\xi_i}_{i = 1}^{\bar{N}}$ is a batch of $\bar{N}$ IID sampled scenarios. A confidence interval around the EEV is then given by
\begin{equation*}
  \brackets*{L_{EEV},U_{EEV}} = \brackets*{EEV-\frac{z_{\alpha/2}\sigma^2_{EEV}}{\sqrt{\bar{N}}}, EEV+\frac{z_{\alpha/2}\sigma^2_{EEV}}{\sqrt{\bar{N}}}},
\end{equation*}
where
\begin{equation*}
  \sigma_{EEV}^2 = \frac{1}{\bar{N}-1}\sum_{i = 1}^{\bar{N}} \parentheses{Q(\bar{x}, \xi_i) - EEV}^2.
\end{equation*}
Now, if there is no overlap between the confidence interval around VRP and the confidence interval around EEV, there is a VSS that is statistically significant to the chosen significance level $\alpha$. A confidence interval around this VSS is then given by
\begin{equation*}
  \brackets*{L_{VSS},U_{VSS}} = \brackets*{L_{VRP} - U_{EEV}, U_{VRP} - L_{EEV}}.
\end{equation*}
We will utilize these techniques to calculate confidence intervals around the stochastic solution in all considered hydropower planning models.

\subsection{StochasticPrograms.jl}
\label{sec:spjl}

The open-source framework \jlinl{StochasticPrograms.jl}~\cite{spjl}, or SPjl for short, is implemented in the Julia programming language and allows the user to efficiently formulate and solve stochastic programs. Moreover, it is designed to scale seamlessly to distributed environments. The framework also includes a solver suite with efficient implementations of the structure-exploiting L-shaped, progressive-hedging, and quasi-gradient algorithms. Each algorithm has a parallel extension that can solve large stochastic programs distributed over multiple cores. SPjl provides a domain-specific language for stochastic programming, as exemplified in Listing~\ref{lst:example}, leveraged by the algebraic modeling language JuMP~\cite{jump}.

\begin{lstlisting}[language = julia, float, caption = {Example declaration of a stochastic program.}, label = {lst:example}]
@stochastic_model begin
    @stage 1 begin
        @decision(model, x[i in 1:10] >= 0)
        @objective(model, Min, sum(x))
    end
    @stage 2 begin
        @parameters W
        @uncertain q T h
        @recourse(model, y[j in 1:5] >= 0)
        @objective(model, Max, q⋅y)
        @constraint(model, T * x  + W * y .== h)
    end
end
\end{lstlisting}
\noindent
The code listing defines a general stochastic model object that can be used to generate specific stochastic programming instances. If a list of scenarios is provided, a finite stochastic program in the form~\eqref{eq:finitesp} is created. This model can then be solved efficiently using decomposition algorithms. The user can then conviniently query the resulting VRP and calculate the VSS. SPjl provides a variety of computational tools for analysing stochastic programs. See the online documentation\footnote{\url{https://martinbiel.github.io/StochasticPrograms.jl/dev/}} for more details. If we instead provide a sampler capable of generating scenarios, we can run the sample average approximation scheme outline above to compute confidence intervals around the stochastic solution. Sampled subproblems instantiated during the SAA algorithm can be distributed and solved efficiently using the framework tools available for finite models. We will showcase how SPjl can be utilized to pose and solve complex and large hydropower planning models.

\subsection{Physical data}
\label{sec:physical-data}

The physical data parameters used in all three model formulations are provided here. The deterministic parameters constitute physical hydro plant parameters and trade regulations. Physical parameters for the power stations in Skellefteälven is available in~\cite{Sag2018,Sandstrom2019} and are provided in Table~\ref{tab:skelleftealven}. We will use this data in all three hydropower models presented in the subsequent sections, and refer back to this table.
\begin{table}[htbp]
  \centering
  \begin{adjustbox}{max width = \textwidth}
    \begin{tabular}{c|cccc}
      \toprule
      \textbf{Plant} & Capacity $\bar{P}$ [MW] & Maximum discharge $\bar{Q}$ [m³/s] & Maximum volume $\bar{M}$ [HE] & Flow time (Q/S) [min] \\
      \midrule
      1. Rebnis & 64 & 80 & 205560 & 2880/2880 \\
      2. Sadva  & 31 & 70 & 168000 & 2880/2880 \\
      3. Bergnas & 8 & 160 & 425280 & 60/60 \\
      4. Slagnas & 7 & 160 & 768 & 240/240 \\
      5. Bastusel & 100 & 170 & 8208 & 60/150 \\
      6. Grytfors & 31 & 165 & 1248 & 15/15 \\
      7. Gallejaur & 214 & 310 & 3600 & 30/150 \\
      8. Vargfors & 131 & 320 & 4008 & 180/180 \\
      9. Rengard & 36 & 220 & 1400 & 180/180 \\
      10. Batfors & 42 & 280 & 1330 & 180/180 \\
      11. Finnfors & 54 & 300 & 300 & 180/180 \\
      12. Granfors & 40 & 240 & 280 & 180/180 \\
      13. Krangfors & 62 & 240 & 330 & 180/180 \\
      14. Selsfors & 61 & 300 & 500 & 180/180 \\
      15. Kvistforsen & 130 & 300 & 1120 & -/- \\
      \bottomrule
    \end{tabular}
  \end{adjustbox}
  \caption{Physical properties of the hydropower plants in Skellefteälven. Displays installed capacity, maximum discharge and reservoir content, as well as the time taken for water to flow to the next downstream plant when discharging or spilling water respectively.}
  \label{tab:skelleftealven}
\end{table}
The physical parameters include reservoir capacities, discharge limits, and water travel time between adjacent stations. Water volume is measured in hour equivalents (HE), which corresponds to a water flow of $1\,\mathrm{m}^{3}/\mathrm{s}$ during one hour. Trade regularizations, including for example trading fees and order limits, are available at NordPool~\cite{nordpool}. Next, we describe uncertain parameters that are also present in all three models.

\subsection{Uncertainty modeling}
\label{sec:uncertainty-modeling}

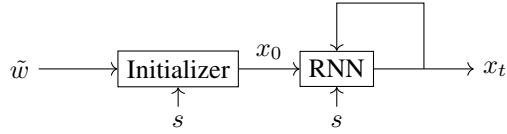
\begin{figure}
  \centering
  \input{rnnarch.tex}
  \caption{Noise-driven RNN forecaster architecture. $\tilde{w}$ is a Gaussian input and $s$ denotes seasonal parameters.}
  \label{fig:rnnarch}
\end{figure}

Two major random elements in hydropower planning are the unknown future electricity price and the unkown future inflow of water to reservoirs. With the aim of posing hydropower planning problems in SPjl, we will first consider how to model these random elements. To this end, we have proposed a noise-driven recurrent neural network (RNN) structure for forecasting electricity prices and local inflow to water reservoirs in the Swedish river Skellefteälven~\cite{dayahead}. We provide a brief re-cap of the essentials. The main aim of the proposed RNN structure is to enable forecasting of sequential data with seasonal variation, without having to rely on long input sequences. The general structure of the forecaster is shown in Figure~\ref{fig:rnnarch}.
The proposed network structure consists of two key components. First, an initializer network is used to compute the initial state of the forecasted sequence. The inputs to the initializer network are a set of seasonal indicators $s$ and a Gaussian noise signal $\tilde{w}$. The structure is used to pose both a price forecaster and an inflow forecaster. The price forecaster was trained on historical price from NordPool~\cite{nordpooldata}, while the inflow forecaster was trained on historical local inflow data in Skellefteälven from the Swedish Meteorological and Hydrological Institute (SMHI)~\cite{shype}. We showed in~\cite{dayahead} that the trained forecasters can generate realistic price and inflow scenarios, with seasonal variation, from a sampled noise signal. We will utilize these forecasters in the following sections, where we consider three case studies on stochastic planning related to the hydropower stations in the Swedish river Skellefteälven.

\section{Case study 1: Day-ahead planning}
\label{sec:dayahead}

In this section, we present a large-scale day-ahead problem. We have already considered this problem in~\cite{dayahead}. We provide a more detailed overview of the model here and re-cap the essential results of the earlier study. We will then expand upon the model in the second case study.

\subsection{The day-ahead market}
\label{sec:day-ahead-market}

Electricity trading in the Nordic energy market is mainly driven by day-ahead auctions. Market participants submit orders of price and electricity volumes for the upcoming day before the market price is known. Any imbalances in settled orders and available production can then be resolved on balancing markets. Hydropower producers can store water in the reservoirs for later use and are therefore able to submit strategic day-ahead orders.

The Nordic day-ahead market offers four order variants for trading electricity volumes, hourly orders, block orders, exclusive groups, and flexible orders. We give a brief introduction to hourly orders and regular block orders.

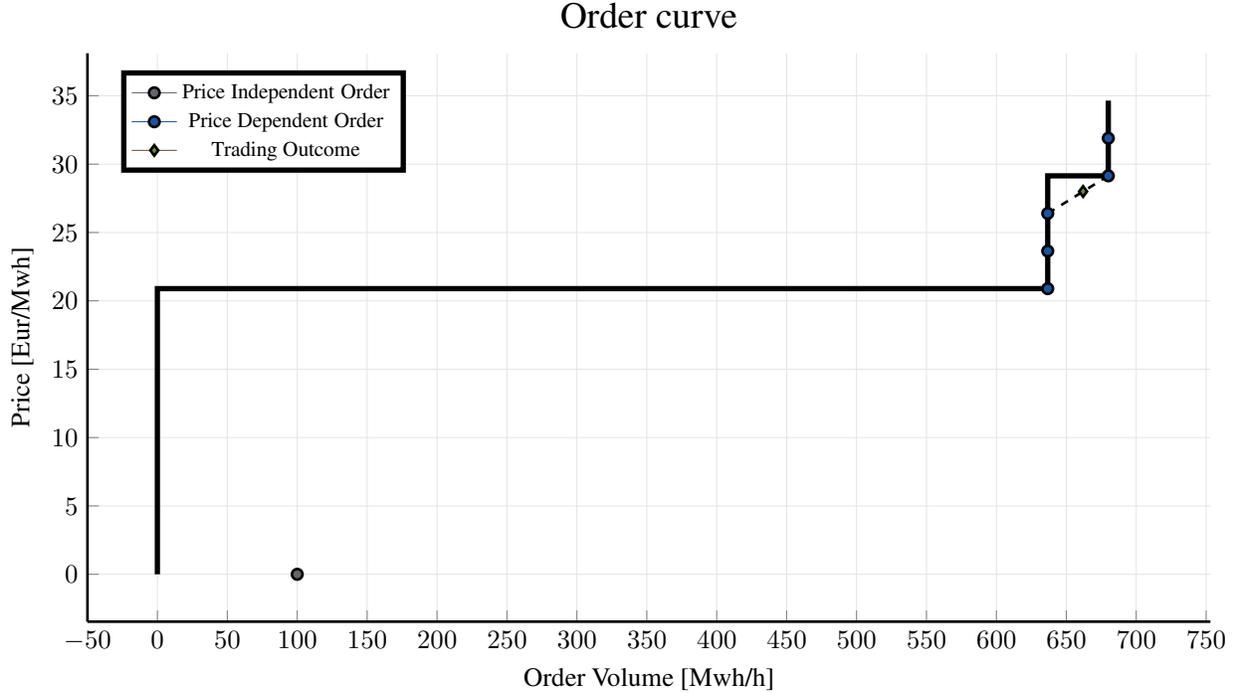
\begin{figure}
  \centering
  \input{single_order_ex.tex}
  \caption{Single hourly order example, showing volume interpolation after market price settlement. The price independent order is always accepted at market price.}
  \label{fig:hourlyorderex}
\end{figure}

Hourly orders can be placed in two ways. A price independent hourly order specifies an electricity volume that is to be purchased or sold at market price during a certain hour, independent of the market price. A price dependent order specifies electricity volumes at given price points. If the settled market price ends up between specified price steps a linear interpolation is performed between the adjacent volume orders to determine the order volume. A settled hourly order is illustrated in Figure~\ref{fig:hourlyorderex}.

Block orders span over an interval of consecutive hours. A regular block order is accepted in its entirety if the mean market price in the specified interval is higher or equal to the order price. The participants then become balance responsible for the order volume every hour of the specified interval, at the mean market price in the interval. Block orders where the price is higher than the mean market price in the given interval are rejected. Other block order variants exist, such as linked block orders and curtailable block orders. These all include conditional elements, and determining optimal orders would involve combinatorial optimization. This is also true for the remaining order types: exclusive group and flexible orders. We do not give further details into the conditional order types as only hourly orders and regular block orders will be used subsequently.

\subsection{Problem setting}
\label{sec:problem-setting}

We formulate the day-ahead planning problem from the perspective of a fictional hydropower producer that owns all 15 power stations in the Swedish river Skellefteälven. The producer is assumed to be price-taking so submitted orders do not influence the market price. The day-ahead model is limited to include only hourly orders and block orders. After market price realization, the producer optimizes the hydropower production with respect to the price and the future water value. Any imbalances are settled in a balancing market at a penalized prize. We assume that there are no fixed contracts to adhere to. In other words, all electricity production is sold for profit in the market. All submitted orders must adhere to the trade regulations specified by the Nordpool market. A general description of the day-ahead problem is given in~\eqref{eq:generaldayahead}.

\begin{equation} \label{eq:generaldayahead}
\begin{aligned}
 \maximize_{} & \quad \text{Profit} + \text{Water value} - \text{Imbalance penalty} \\
 \text{subject to} & \quad \text{Order strategy} \\
 & \quad \text{Physical limitations} \\
 & \quad \text{Economic/legal limitations}
\end{aligned}
\end{equation}

Because next-day market prices are unknown when placing orders, we formulate a two-stage stochastic program to generate optimal orders. The first-stage decisions are the orders submitted to the day-ahead market. A general description of the first-stage problem is given in~\eqref{eq:generalfirststage}.

\begin{equation} \label{eq:generalfirststage}
\begin{aligned}
 \maximize_{\mathclap{\text{Order strategy}}} & \quad \expect[]{\text{Revenue}\parentheses*{\text{Order strategy}, \text{Price}, \text{Inflow}}} \\
 \text{subject to} & \quad \text{Trade regulations} \\
\end{aligned}
\end{equation}

In each second-stage scenario, uncertain parameters are realized and the electricity production is optimized with respect to profits and water value while satisfying the settled order commitments. A general description of the second stage is given in~\eqref{eq:generalsecondstage}
\begin{equation} \label{eq:generalsecondstage}
\begin{aligned}
 \maximize_{\mathclap{\text{Production schedule}}} & \quad \text{Profit}\parentheses*{\text{Price}} + \text{Water value} - \text{Imbalance penalty} \\
 \text{subject to} & \quad \text{Commitments}\parentheses*{\text{Order strategy}, \text{Price}} \\
 & \quad \text{Hydrological balance}\parentheses*{\text{Inflow}} \\
 & \quad \text{Electricity production} \\
 & \quad \text{Load balance}\parentheses*{\text{Commitments}}
\end{aligned}
\end{equation}

In addition, the producer can take recourse decisions by trading surplus or shortage in a simple balancing market. Both market prices and local water inflows to the reservoirs are considered uncertain.

\subsection{Day-ahead model}
\label{sec:day-ahead-model}

The day-ahead model used in this work is similar to the model introduced in~\cite{Fleten2007}. We repeat the general structure and also highlight the key details of our model. Further, we sketch how the model is implemented in our software framework.

\subsubsection{General structure}
\label{sec:general-structure}

In the first stage, we model the day-ahead orders to be submitted to the market. Let $\mathcal{T} = \braces{t_1,\dots,t_{24}}$ denote indices for the $24$-hour horizon of the upcoming day. The set $\mathcal{B} = \braces{b_1,\dots,b_B}$ make up blocks of consecutive hours in the $24$-hour period. In order to avoid non-linear relations in the model, we fix a set of hourly price levels $\mathcal{P}_t = \braces{p_{1,t},\dots,p_{P,t}}$ to bid at beforehand. We explain how these prices are chosen in a following subsection. The block order prices are determined by calculating averages of the available prices levels over the given blocks. We introduce $x^I_t, x^D_{p,t}$, and $x^B_{p,b}$ to represent price independent orders, price dependent orders, and block orders respectively. As per NordPool regulations, the volumes in a price dependent sell order have to be constant or increasing with increasing prices. We enforce the following constraints:
\begin{equation*}
  x^D_{i,t} \leq x^D_{i+1,t}, \quad i \in \mathcal{P}, t \in \mathcal{T}.
\end{equation*}
In addition, we constrain the total volume offered to the market to not exceed $200\%$ of the production capacity, i.e.,
\begin{equation*}
  x^I_t + x^D_{P,t} + \sum_{i \in \mathcal{P}}\sum_{t \in b, b \in \mathcal{B}} x^B_{i,b} \leq 2 \sum_{h \in \mathcal{H}}\bar{P}_h,
\end{equation*}
where $\bar{P}_h$ is the maximum production of plant $h$ and can be obtained from Table~\ref{tab:skelleftealven}. Consequently, we allow imbalances in the order commitments, but limit the maximum imbalance already in the first stage.

In the second stage, we model the order commitments after price realization as well as the production schedule after inflow realization. We introduce the random variables $\rho^{\omega}_t$, that describe the hourly market price in scenario $\omega$. Let $y^{\omega}_t$ and $y^{\omega}_b$ represent the committed hourly volumes and the committed block volumes, in scenario $\omega$, respectively. Every hour $t$, the dispatched hourly volumes are determined through linear interpolation:
\begin{equation*}
  y^{\omega}_t = x^I_t + \frac{\rho^{\omega}_t-p_{i,t}}{p_{i+1,t} - p_{i,t}}x^D_{i+1,t} + \frac{p_{i+1,t}-\rho^{\omega}_t}{p_{i+1,t} - p_{i,t}}x^D_{i, t} \quad p_{i,t} \leq \rho^{\omega}_t \leq p_{i+1,t}.
\end{equation*}
The dispatched block volumes are given by
\begin{equation*}
  y^{\omega}_b = \sum_{\mathclap{p : \bar{p}(p,b) \leq \bar{\rho}^{\omega}_b}} x^B_{p,b}
\end{equation*}
where
\begin{equation*}
N  \bar{p}(i,b) = \frac{1}{|b|}\sum_{t \in b} p_{i,t}
\end{equation*}
and
\begin{equation*}
  \bar{\rho}^{\omega}_b = \frac{1}{|b|}\sum_{t \in b} \rho^{\omega}_t.
\end{equation*}
Next, we model the production. Let $\mathcal{H} = \braces{h_1,\dots,h_{15}}$ index the 15 hydroelectric power stations in Skellefteälven. For each plant and hour, let $Q^{\omega}_{h,s,t}$ and $S^{\omega}_{h,t}$ denote the water discharged and spilled in scenario $\omega$, respectively. The maximum discharge, $\bar{Q}_{h}$ is obtained from Table~\ref{tab:skelleftealven}. Further, let $P^{\omega}_t$ denote the total volume of electricity produced each hour in scenario $\omega$. We employ a piecewise linear approximation of the generation curve of each station. In other words,
\begin{equation*}
  P^{\omega}_t = \sum_{s \in \mathcal{S}} \mu_{h,s}Q^{\omega}_{h,s,t},
\end{equation*}
where $\mu_{h,s}$ is the marginal production equivalent of station $h$ and segment $s \in \mathcal{S}$. See~\cite{obel} for how to estimate the production curve in two piecewise linear segments based on the maximum discharge and capacity. In brief, we set
\begin{equation*}
  \begin{aligned}
    \mu_{h,1} &= \frac{\bar{P}_h}{\bar{Q}_{s}(0.75 + 0.95\cdot 0.25)} \\
    \mu_{h,2} &= 0.95\mu_{h,1}
  \end{aligned}
\end{equation*}
and let $\bar{Q}_{h,1} = 0.75 \bar{Q}_h$ in segment 1 and $\bar{Q}_{h,2} = 0.25 \bar{Q}_h$ in segment 2. The load balance is given by
\begin{equation*}
  y^{\omega} + \sum_{t \in b, b \in \mathcal{B}} y^{\omega}_b  - P^{\omega}_t = y^{\omega +}_t - y^{\omega -}_t.
\end{equation*}
In each hour in scenario $\omega$, any imbalance between committed volumes and produced volumes is equal to the difference between the imbalance variables $y^{\omega}_{t+}$ and $y^{\omega}_{t-}$. Any shortage $y^{\omega}_{t+}$ is bought from the balancing market, and any surplus $y^{\omega}_{t-}$ is sold to the balancing market. Finally, let $M^{\omega}_{h,t}$ denote the reservoir contents in plant $h$ during hour $t$. The maximum reservoir content, $\bar{M}_{h}$ is obtained from Table~\ref{tab:skelleftealven}. Flow conservation each hour is given by
\begin{equation*}
  \begin{aligned}
    M^{\omega}_{h,t} = \;&M^{\omega}_{h,t-1} \\
    &+ \sum_{i \in \mathcal{Q}_u(h)} \sum_{s \in \mathcal{S}} Q^{\omega}_{i,s,t-\tau_{ih}} + \sum_{i \in \mathcal{S}_u(h)} S^{\omega}_{i,t-\tau_{ih}} + V^{\omega}_h \\
    &- \sum _{s \in \mathcal{S}} Q^{\omega}_{h,s,t} - S^{\omega}_{h,t}
  \end{aligned}
\end{equation*}
Here, $V^{\omega}_h$ are random variables describing the local inflow to each plant in scenario $\omega$. The sets $\mathcal{Q}_u(h)$ and $\mathcal{S}_u(h)$ contain upstream plants where discharge and/or spillage can reach plant $h$ through connecting waterways. Note that the water travel times $\tau_{ih}$ between power stations are included in the incoming flow to each plant, and can be obtained from Table~\ref{tab:skelleftealven}. Internally, this is modeled by introducing auxiliary variables and constraints. Variable limits and the introduced parameters are all included in the deterministic data sets for Skellefteälven given in Table~\ref{tab:skelleftealven}. The revenue from a production schedule satisfying the above relations is given by
\begin{equation*}
  \sum_{t \in \mathcal{T}} \rho^{\omega}_ty^{\omega}_t + \sum_{b \in \mathcal{B}}\abs{b}\bar{\rho}^{\omega}_by^{\omega}_b + \sum_{t \in \mathcal{T}}\alpha_t\rho^{\omega}_ty^{\omega}_{t-} - \beta_t\rho^{\omega}_ty^{\omega}_{t+} + W^{\omega}(M^{\omega}_{1,24},\dots,M^{\omega}_{15,24}).
\end{equation*}
Note that, for any committed block order $y^{\omega}_b$, the order volume is dispatched every hour in the block at average market price. Hence, $\abs{b}\bar{\rho}^{\omega}_by^{\omega}_b$ is earned. The imbalance volumes are traded at penalized prices, using penalty factors $\alpha_t$ and $\beta_t$, for discouragement. It is hard to accurately model this penalty. Here, we use a $15\%$ penalty during peak hours, and $10\%$ otherwise. These values are based on observations of historic values, but can not be considered accurate. The final term in the revenue is the expected future value of water, which is a function of the water volumes that remain in the reservoirs after the period. In the following section, we introduce a polyhedral approximation of this function that can be modelled with linear terms. For now, we simply denote the water value in scenario $\omega$ by $W^{\omega}$. In summary, a stochastic program modeling the day-ahead problem is in essence given by
\begin{equation} \label{eq:dayahead}
  \begin{aligned}
   \maximize_{x^I_t, x^D_{i,h}, x^B_{i,b}} & \quad \expect[\xi]{\sum_{t \in \mathcal{T}} \rho^{\omega}_ty^{\omega}_t + \sum_{b \in \mathcal{B}}\abs{b}\bar{\rho}^{\omega}_by^{\omega}_b + \sum_{t \in \mathcal{T}} \parentheses*{\alpha_t\rho^{\omega}_ty^{\omega}_{t-} - \beta_t\rho^{\omega}_ty^{\omega}_{t+}} + W^{\omega}} \\
   \sbj & \quad x^D_{i,t} \leq x^D_{i+1,t}, \quad i \in \mathcal{P}, t \in \mathcal{T} \\
   & \quad x^I_t + x^D_{P,t} + \sum_{i \in \mathcal{P}}\sum_{t \in b, b \in \mathcal{B}} x^B_{i,b} \leq 2 \sum_{h \in \mathcal{H}}\bar{P}_h \\
   & \quad y^{\omega}_t = x^I_t + \frac{\rho^{\omega}_t-p_{i,t}}{p_{i+1,t} - p_{i,t}}x^D_{i+1,t} + \frac{p_{i+1,t}-\rho^{\omega}_t}{p_{i+1,t} - p_{i,t}}x^D_{i, t}, \quad t \in \mathcal{T} \\
   & \quad y^{\omega}_b = \sum_{p : \bar{p}(p,b) \leq \bar{\rho}^{\omega}_b} x^B_{p,b}, \quad b \in \mathcal{B} \\
   & \quad P^{\omega}_t = \sum_{s \in \mathcal{S}} \mu_{h,s}Q^{\omega}_{h,s,t}, \quad t \in  \mathcal{T} \\
   & \quad y^{\omega}_t + \sum_{\mathclap{t \in b, b \in \mathcal{B}}} y^{\omega}_b  - P^{\omega}_t = y^{\omega}_{t+} - y^{\omega}_{t-}, \quad t \in \mathcal{T} \\
   & \quad \begin{aligned}
     M^{\omega}_{h,t} = \;&M^{\omega}_{h,t-1} \\
     &+ \sum_{i \in \mathcal{Q}_u(h)} \sum_{s \in \mathcal{S}} Q^{\omega}_{i,s,t-\tau_{ih}} + \sum_{i \in \mathcal{S}_u(h)} S^{\omega}_{i,t-\tau_{ih}} \\
     &+ V^{\omega}_h \\
     &- \sum _{s \in \mathcal{S}} Q^{\omega}_{h,s,t} - S^{\omega}_{h,t}, \qquad h \in \mathcal{H}, t \in \mathcal{T} \\
   \end{aligned} \\
   & \quad 0 \leq Q^{\omega}_{h,s,t} \leq \bar{Q}^{\omega}_{h,s}, \quad h \in \mathcal{H}, s \in \mathcal{S}, t \in \mathcal{T} \\
   & \quad 0 \leq M^{\omega}_{h,t} \leq \bar{M}_h, \quad h \in \mathcal{H}, t \in \mathcal{T} \\
   & \quad y^{\omega}_t \geq 0,\; y^{\omega}_b \geq 0,\; y^{\omega}_{t+} \geq 0,\; y^{\omega}_{t-} \geq 0 \\
   & \quad P^{\omega}_t \geq 0,\; S^{\omega}_{h,t} \geq 0.
 \end{aligned}
\end{equation}

\subsubsection{Water evaluation}
\label{sec:water-evaluation}

The expected value of keeping water in the reservoirs must be accounted for in the production plan. If the water value is large, then it could be optimal to not produce, settle committed orders in the balancing market, and save water. Likewise, if the water value is small, it could be optimal to overproduce and sell the excess in the balancing market. Consequently, the water value will evidently also impact the optimal order strategy because the order commitments are instrumental in both scenarios. Thus, the accuracy of the water evaluation is critical for hydropower producers participating in the day-ahead market. If we assume that excess water can be used to produce and sell electricity at some expected future price, we get naive order strategies governed by price variations around the expected future price. We instead consider an auxilliary stochastic program, where the first-stage decisions determine the reservoir contents of every power station before the upcoming week. After realizing sequence of inflows and daily price curves, the second stage optimizes the weekly production of energy sold at market price. This simplified week-ahead problem is given by

\begin{equation} \label{eq:weekahead}
  \begin{aligned}
   \maximize_{M_{h,0}} & \quad \expect[\xi]{\sum_{t \in \tilde{\mathcal{T}}} \rho^{\omega}_tP^{\omega}_t} \\
   \st & \quad P^{\omega}_t = \sum_{s \in \mathcal{S}} \mu_{h,s}Q^{\omega}_{h,s,t}, \quad t \in \mathcal{T} \\
   & \quad \begin{aligned}
    M^{\omega}_{h,t} &= M^{\omega}_{h,t-1} \\
    &+ \sum_{i \in \mathcal{Q}_u(h)} \sum_{s \in \mathcal{S}} Q^{\omega}_{i,s,t-\tau_{ih}} + \sum_{i \in \mathcal{S}_u(h)} S^{\omega}_{i,t-\tau_{ih}} \\
    & + V^{\omega}_h \\
    &- \sum _{s \in \mathcal{S}} Q^{\omega}_{h,s,t} - S^{\omega}_{h,t}, \qquad h \in \mathcal{H}, t \in \mathcal{T} \\
  \end{aligned} \\
   & \quad 0 \leq Q^{\omega}_{h,s,t} \leq \bar{Q}_{h,s}, \quad h \in \mathcal{H}, t \in \mathcal{T} \\
   & \quad 0 \leq M^{\omega}_{h,t} \leq \bar{M}_h, \quad h \in \mathcal{H}, t \in \mathcal{T} \\
   & \quad P^{\omega}_t \geq 0,\; S^{\omega}_{h,t} \geq 0,
 \end{aligned}
\end{equation}
\noindent
\sloppy
where the time-horizon $\tilde{\mathcal{T}} = \braces{t_1,\dots,t_{168}}$ is now a week. The objective function $W^{\omega} = \expect[\xi]{W^{\omega}(M_{1,0},\dots,M_{15,0})}$ of this problem will be used as a water value function. The problem~\eqref{eq:weekahead} is trivial since the optimal decision will be to fill the reservoirs with enough water to be able to run at maximum capacity in the worst-case scenario. However, information about the water value can be extracted by solving~\eqref{eq:weekahead} with an L-shaped type method. The L-shaped method generates cutting planes of the form
\begin{equation} \label{eq:wcuts}
  \sum_{h \in \mathcal{H}}\partial W_{c,h}M_{h,0}+ W \geq w_c.
\end{equation}
\noindent
This form supports for the concave objective function $W$, which is a function of reservoir content in the system. Hence, after the algorithm has converged we have access to a polyhedral approximation of $W$ in the form of a collection of such cuts as~\eqref{eq:wcuts}. We can use these cuts to put an approximate value of the remaining volumes of water present in the reservoirs after meeting order commitments. The water value approximation enters the day-ahead problem~\eqref{eq:dayahead} in the following way:
\begin{equation*}
\begin{aligned}
 \maximize_{} & \quad \dots + W \\
 \st & \quad \vdots \\
 & \quad \sum_{h \in \mathcal{H}}\partial W_{c,h}M_{h,24}+ W \geq w_c \quad c \in \mathcal{C}. \\
\end{aligned}
\end{equation*}
In practice, we use a multiple-cut formulation
\begin{equation*}
  W = \sum_{i = 1}^{N} W_i
\end{equation*}
as the L-shaped method solves the week-ahead problem with a large number of scenarios $N$ more efficiently in this way. The end result is still a collection of cuts that approximate a polyhedral water value function of the final reservoir volumes.

\subsubsection{Price levels}
\label{sec:pricelevels}

\begin{figure}
  \centering
  \input{price_levels.tex}
  \caption{Expected daily electricity price out of $1000$ samples from the RNN forecaster. Two standard deviations above and below the expected price is shown each hour.}
  \label{fig:pricelevels}
\end{figure}
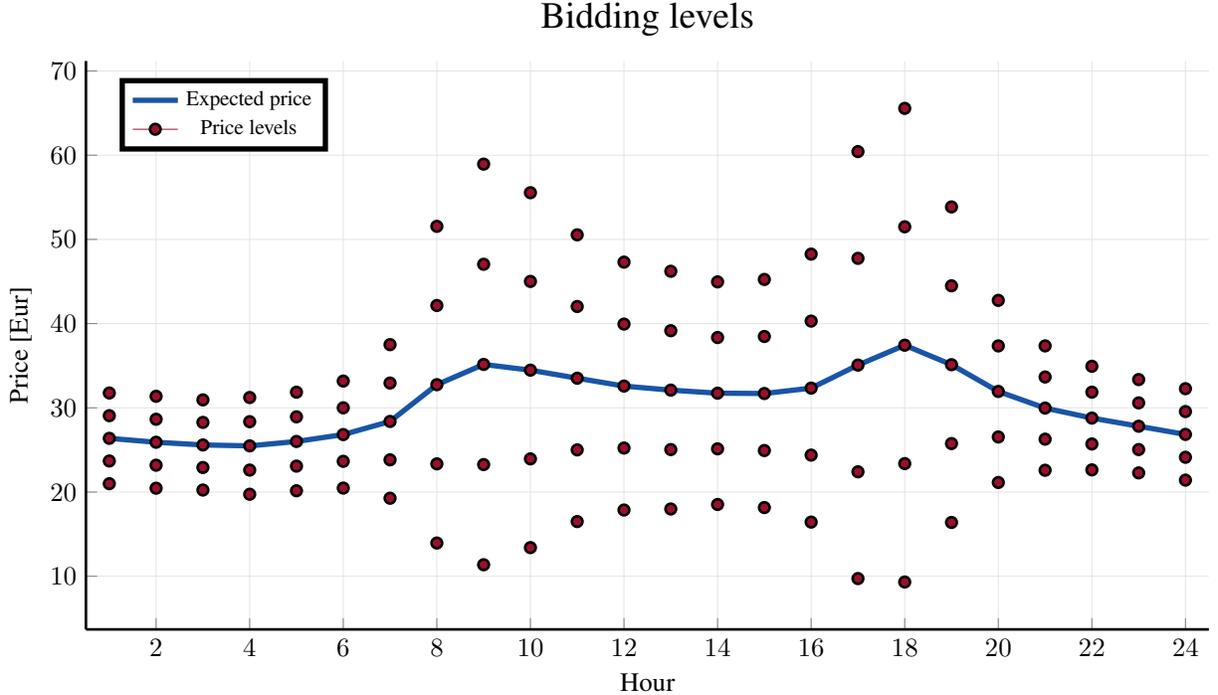

The price-dependent hourly orders and the block orders are specified at pre-chosen price levels. For flexibility, we allow these levels to vary with time. The set of price levels $\mathcal{P}_t$ for each hour is determined using the price forecaster. We sample a set of price scenarios and use the resulting hourly mean price and standard deviation as a baseline. In each hour $t$, we define price levels around the mean price using multiples of the standard deviation. A set of hourly price levels generated using this method is shown in Figure~\ref{fig:pricelevels}. For each block $b \in \mathcal{B}$, we define the possible block prices $\braces{p_{i,b}}_{i = 1}^{5}$ by computing mean price levels over the hours $t \in b$.

\subsubsection{Model implementation}
\label{sec:model-implementation}

We now outline how the day-ahead model~\eqref{eq:dayahead} is formulated in SPjl. To increase readability, we present an abridged version of the day-ahead model implementation in SPjl, by obfuscating parts of the code and making slight syntax changes. The full unabridged model is available at Github~\footnote{\url{https://github.com/martinbiel/HydroModels.jl}}.

First, we define a data structure to describe the uncertain parameters using the \jlinl{@define_scenario} command. We also create a sampler object, using the \jlinl{@sampler} command, which utilizes the noise-driven RNN forecasters introduced in~\cite{dayahead} to generate price curves and inflows. The code is shown in Listing~\ref{lst:dayahead_scenariodef}. Because we want to make use of the forecasters' seasonal capabilities, we also include a date field in the sampler object. The forecasters use the provided date to determine seasonal parameter inputs to the neural networks.

\begin{lstlisting}[language = julia, float, caption = {Day-ahead scenario definition in SPjl}, label = {lst:dayahead_scenariodef}]
@define_scenario DayAheadScenario = begin
    ρ::PriceCurve{24,Float64}
    Q̃::Inflows{typeof(Skellefteälven),Float64}
end

@sampler RecurrentDayAheadSampler = begin
    date::Date
    plants::PlantCollection
    price_forecaster::Forecaster
    flow_forecaster::Forecaster

    @sample DayAheadScenario begin
        price_curve = forecast(sampler.price_forecaster, month(sampler.date))
        flows = forecast(sampler.flow_forecaster, week(sampler.date))
        return DayAheadScenario(PriceCurve(price_curve), Inflows(flows))
    end
end
\end{lstlisting}

The day-ahead model definition in SPjl is presented in Listing~\ref{lst:dayaheaddef}.

\begin{lstlisting}[language = julia, float, caption = {Day-ahead problem definition in SPjl. The code has been condensed for readability.}, label = {lst:dayaheaddef}]
@stochastic_model begin
    @stage 1 begin
        @parameters horizon indices data
        @unpack hours, plants, bids, blockbids, blocks = indices
        @unpack hydrodata, regulations = data
        @decision(model, xᴵ[t in hours] >= 0) # Price-independent orders
        @decision(model, xᴰ[p in bids, t in hours] >= 0) # Price-dependent orders
        @decision(model, xᴮ[p in blockbids, b in blocks] >= 0) # Block orders
    end
    @stage 2 begin
        @parameters horizon indices data
        @unpack hours, plants, segments, blocks = indices
        @unpack hydrodata, water_value, regulations, bidlevels = data
        @uncertain ρ, V from ξ::DayAheadScenario
        @recourse(model, yᴴ[t in hours] >= 0) # Dispatched hourly volumes
        @recourse(model, yᴮ[b in blocks] >= 0) # Dispatched block volumes
        @recourse(model, y⁺[t in hours] >= 0) # Energy shortage
        @recourse(model, y⁻[t in hours] >= 0) # Energy surplus
        @recourse(model, 0 <= Q[h in plants,t in hours] <= Q_max[h]) # Discharge
        @recourse(model, S[h in plants,t in hours] >= 0) # Spillage
        @recourse(model, 0 <= M[h in plants,t in hours] <= M_max[p]) # Water volume
        @recourse(model, W[i in 1:nindices(water_value)]) # Water value
        @recourse(model, P[t in hours] >= 0) # Produced energy
        @expression(model, net_profit,
            sum(ρ[t]*yᴴ[t] for t in hours)
            + sum(|b|*(mean(ρ[hours_per_block[b]])*yᴮ[b] for b in blocks))
        @expression(model, intraday_trading,
            sum(penalty(ξ,t)*y⁺[t] - reward(ξ,t)*y⁻[t] for t in hours))
        @expression(model, value_of_stored_water,
            -sum(W[i] for i in 1:nindices(water_value)))
        @objective(model, Max, net_profit - intraday_trading + value_of_stored_water)
        # Bid-dispatch links
        @constraint(model, hourlybids[t in hours],
            yᴴ[t] == interpolate(ρ[t], bidlevels, xᴰ[t]) + xᴵ[t])
        @constraint(model, bidblocks[b in blocks],
            yᴮ[b] == sum(xᴮ[j,b] for j in accepted_blockorders(b)))
        # Hydrological balance
        @constraint(model, hydro_constraints[h in plants, t in hours],
            M[h,t] == (t > 1 ? M[h,t-1] : M₀[p]) #
            + sum(Q[i,t-τ] for i in intersect(Qu[p], plants)) # Inflows from
            + sum(S[i,t-τ] for i in intersect(Su[p], plants)) # upstream plants
            + V[p] # Local inflow
            - (Q[h,t] + S[h,t])) # Outflow
        # Production
        @constraint(model, production[t in hours],
            P[t] == sum(μ[p]*Q[h,t] for h in plants))
        # Load balance
        @constraint(model, loadbalance[t in hours],
            yᴴ[t] + sum(yᴮ[b] for b in active(t)) - P[t] == y⁺[t] - y⁻[t])
        # Water travel time ... (not shown)
        # Polyhedral water value
        @constraint(model, water_value_approximation[c in 1:ncuts(water_value)],
            sum(∂W[c,p]*M[h,T] for h in plants)
            + sum(W[i] for i in cut_indices(c)) >= w[c])
    end
end
\end{lstlisting}

\subsubsection{Algorithmic details}
\label{sec:algorithm-details}

We use the sample average approximation (SAA) scheme outlined in Section~\ref{sec:stoch-progr} to solve the day-ahead problem~\eqref{eq:dayahead}. During the SAA procedure, we solve many sampled instances of increasing size. We distribute the sampled instances on a 32-core compute node, using the parallel capabilities of SPjl. The instances are solved efficiently using a parallel L-shaped method accelerated using regularization~\cite{distlshaped} and cut aggregation~\cite{cutaggregation}.

\subsection{Numerical Experiments}
\label{sec:numerical_experiments}

\begin{figure}
  \centering
  \input{dayahead_confidence_intervals.tex}
  \caption{Confidence intervals around optimal value of the day-ahead problem as a function of sample size.}
  \label{fig:day-ahead_confidence_intervals}
\end{figure}
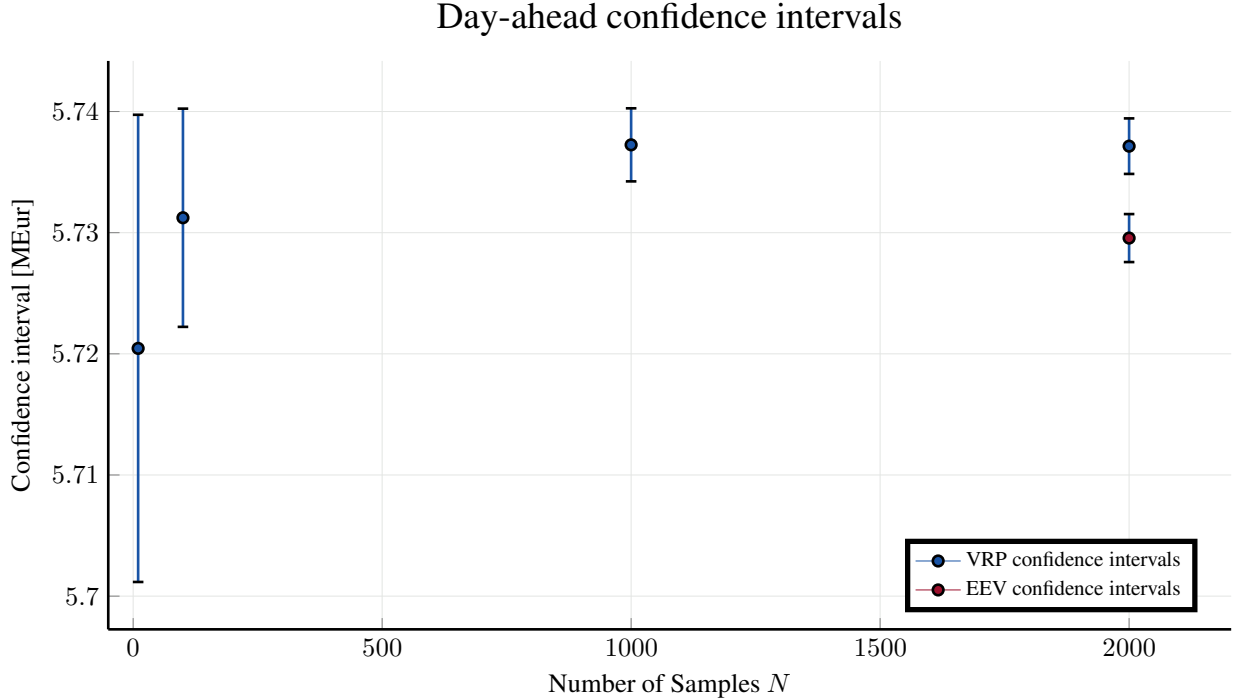

The results of the SAA algorithm is given in Figure~\ref{fig:day-ahead_confidence_intervals}. The confidence interval is stabilized at $2000$ samples. We compute a confidence interval around the EEV at this sample size as well. Because there is no overlap between the VRP and EEV, there is a statistically significant VSS within the interval $[0.058\% - 0.21\%]$. The total profit is however skewed by the future water evaluation. With respect to only the daily market profit, the relative VSS is about $1\%$. In addition, we note that these are daily marginal profits. Hence, the VSS accumulates and could be considered more significant.

\begin{figure}
  \input{day-ahead_evp_strategy.tex}
  \centering
  \caption{Optimal order strategy from a day-ahead expected value problem.}
  \label{fig:dayahead_evp_strategy}
\end{figure}
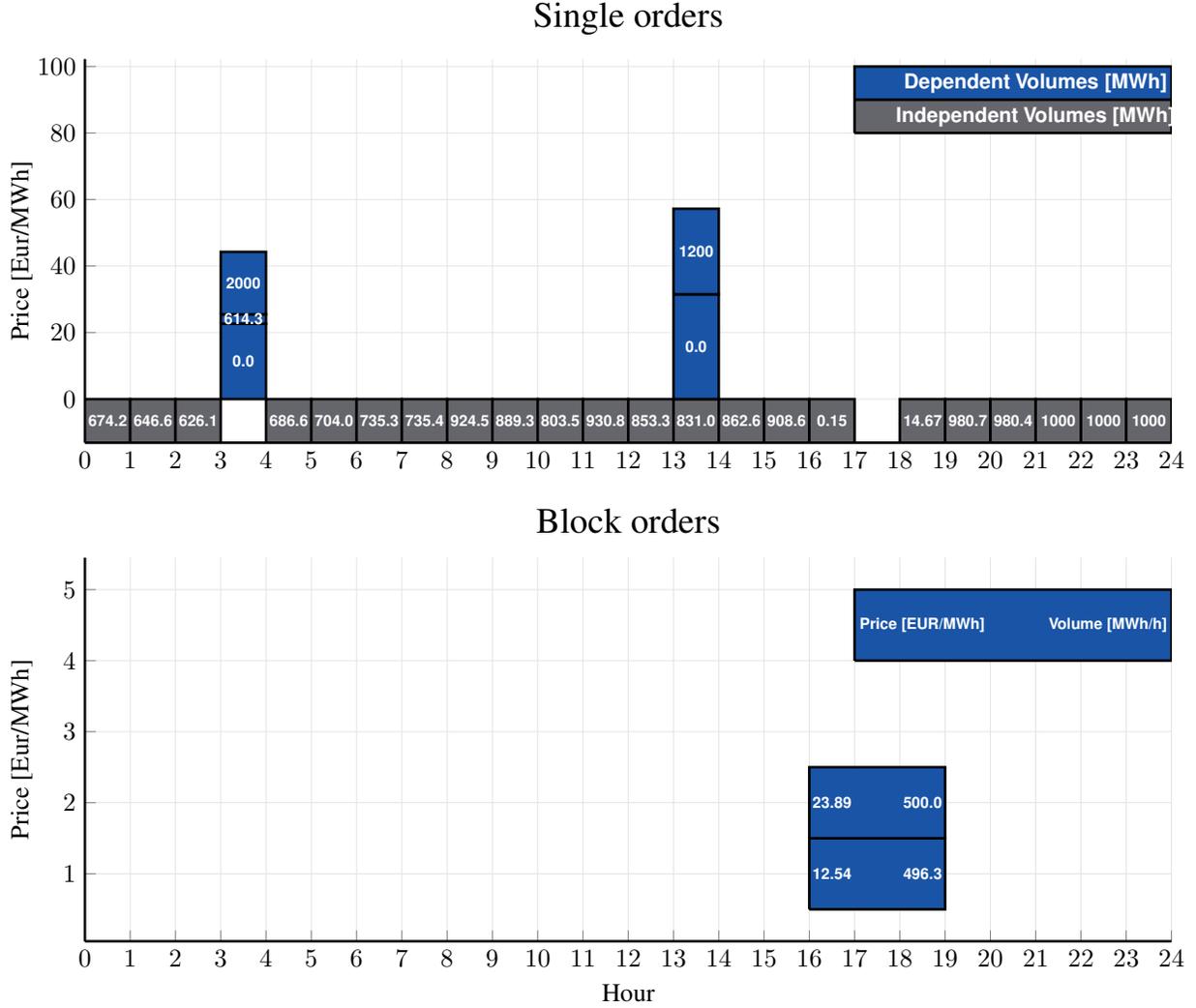

We solve a specific sampled instance of $2000$ scenarios and present the resulting strategy in Figure~\ref{fig:dayahead_strategy}. The stochastic solution uses a large block order in the afternoon where a large mean price is expected. In comparison, the deterministic strategy obtained by solving the expected value problem is shown in Figure~\ref{fig:dayahead_evp_strategy}. The deterministic strategy mostly utilizes price-independent orders, which is less flexible than the stochastic solution.

\begin{figure}
  \input{day-ahead_strategy.tex}
  \centering
  \caption{Optimal order strategy from a $2000$-scenario day-ahead instance.}
  \label{fig:dayahead_strategy}
\end{figure}
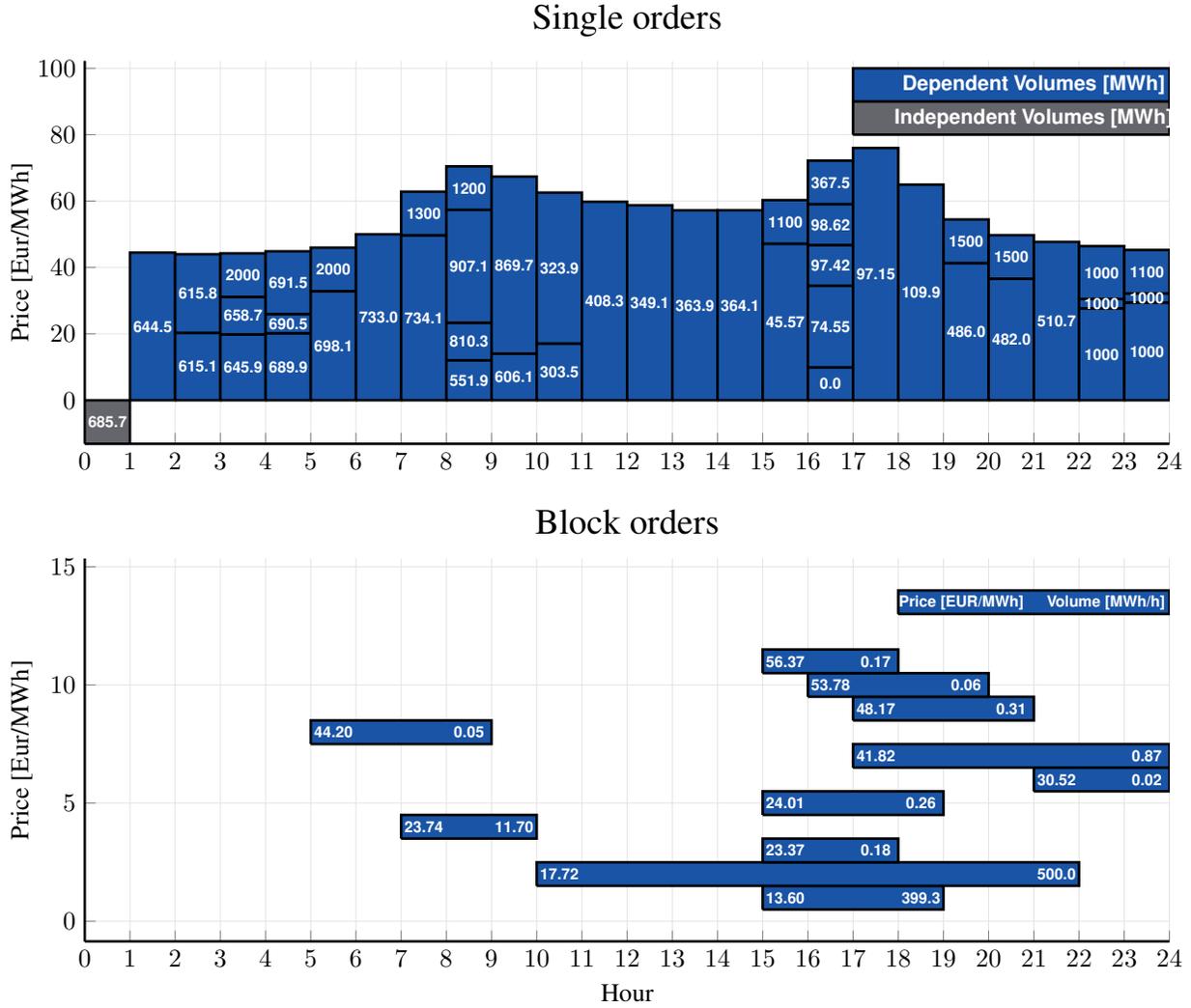

\section{Case study 2: Maintenance scheduling}
\label{sec:maintenance-scheduling}

Next, we consider a variation of the day-ahead problem that includes preventive power plant maintenance. We again employ the RNN forecasters to generate scenarios, formulate the maintenance scheduling problem in SPjl, and solve the problem in parallel using the special-purpose algorithms discussed in the second part of the thesis. We provide confidence intervals and note that the value of the stochastic solution is higher than the day-ahead problem.

\subsection{The maintenance scheduling problem}
\label{sec:maintenance-scheduling-problem_2}

Hydropower production is flexible due to the ability to store energy in water reservoirs. In Sweden, hydropower constitutes about $50\%$ of the total power production, and this flexibility is constantly utilized to ensure balance in the Swedish power system~\cite{balans}. It is therefore important that the hydroelectric power production can operate reliably to ensure stability in the total energy system. System failures should of course be corrected through maintenance, but it is also possible to perform planned preventive maintenance of the power stations to prolong their lifetime and increase reliability~\cite{maintenance}. By planning maintenance, it also becomes possible to minimize the loss of profit in the maintenance period.

We consider settings when electricity production is not possible during preventive maintenance. Because we are interested in the deregulated Nordic electricity market we also consider the orders submitted to the market the day before maintenance is planned. Large losses can be incurred if the power producer is balance responsible for large energy volumes during the maintenance period. It is also possible to miss profits if maintenance coincides with a period when the electricity price is large. Consequently, it is important to coordinate the maintenance scheduling with the order strategy.

Maintenance scheduling for hydroelectric systems have been considered in other contributions~\cite{maintenance,maintenance2}, but they do not include coordination with day-ahead order strategies. We will show that this inclusion leads to a significant value of the stochastic solution compared to a day-ahead formulation.

A hydroelectric maintenance scheduling problem involves specifying an optimal schedule for when to perform preventive maintenance on the hydropower power stations, and optimal order volumes to be submitted to the market, that maximizes the expected profit. We give a brief introduction to the problem in the following.

\subsubsection{Problem setting}
\label{sec:problem-setting_2}

We formulate a maintenance scheduling planning problem from the perspective of a fictional hydropower producer that owns all 15 power stations in the Swedish river Skellefteälven. The producer is again assumed to be price-taking. The day-ahead model is limited to include only hourly orders to ensure computational tractability. When the market price has been realized, the producer optimizes the power production with respect to the price while carrying out preventive maintenance according to a planned schedule. The time required for preventive maintenance of each plant is given in Table~\ref{tab:maintenance_times}. We allot 1-4 hours of maintenance for most power stations, and a longer period of 6-8 hours for the four stations connected to large reservoirs. The maintenance periods are fictional and are meant to resemble a preventive maintenance procedure that can be performed in a single day.

\begin{table}[htbp]
  \centering
  \begin{adjustbox}{max width = \textwidth}
    \begin{tabular}{c|c}
      \toprule
      \textbf{Plant} & {Maintenance time [hours]} \\
      \midrule
      Rebnis & 8 \\
      Sadva  & 6 \\
      Bergnas & 8 \\
      Slagnas & 8 \\
      Bastusel & 4 \\
      Grytfors & 4 \\
      Gallejaur & 3 \\
      Vargfors & 3 \\
      Rengard & 2 \\
      Batfors & 2 \\
      Finnfors & 2 \\
      Granfors & 2 \\
      Krangfors & 2 \\
      Selsfors & 1 \\
      Kvistforsen & 1 \\
      \bottomrule
    \end{tabular}
  \end{adjustbox}
  \caption{Time alloted for preventive maintenance of each power station in Skellefteälven.}
  \label{tab:maintenance_times}
\end{table}

We assume that a power station is not able to discharge water and produce power during the maintenance period. Any imbalances in setted market orders and produced power and are settled in a fictional balancing market at a penalized prize. We assume that there are no fixed contracts to adhere to. In other words, all electricity production is sold for profit in the market. All submitted orders must adhere to the trade regulations specified by the Nordpool market. A general description of the maintenance scheduling problem is given in~\eqref{eq:generalmaintenance}.

\begin{equation} \label{eq:generalmaintenance}
\begin{aligned}
 \maximize_{} & \quad \text{Profit} - \text{Imbalance penalty} \\
 \text{subject to} & \quad \text{Maintenance schedule} \\
 & \quad \text{Order strategy} \\
 & \quad \text{Physical limitations} \\
 & \quad \text{Economic/legal limitations}
\end{aligned}
\end{equation}

Because next-day market prices are unknown when placing orders and constructing the maintenance schedule, we formulate a two-stage stochastic program to generate optimal orders and schedule. The first-stage decisions are the orders submitted to the day-ahead market as well as the preventive maintenance schedule for the upcoming day. A general description of the first-stage problem is given in~\eqref{eq:generalfirststage_2}.

\begin{equation} \label{eq:generalfirststage_2}
\begin{aligned}
 \maximize_{\mathclap{\text{Orders + schedule}}} & \quad \expect[]{\text{Revenue}\parentheses*{\text{Orders}, \text{Schedule}, \text{Price}, \text{Inflow}}} \\
 \text{subject to} & \quad \text{Trade regulations} \\
 & \quad \text{Schedule restrictions}
\end{aligned}
\end{equation}

In each second-stage scenario, uncertain parameters are realized and the electricity production is optimized with respect to next-day profits while satisfying the settled order commitments and adhering to the maintenance schedule. A general description of the second stage is given in~\eqref{eq:generalsecondstage_2}
\begin{equation} \label{eq:generalsecondstage_2}
\begin{aligned}
 \maximize_{\mathclap{\text{Production schedule}}} & \quad \text{Profit}\parentheses*{\text{Price}} - \text{Imbalance penalty} \\
 \text{subject to} & \quad \text{Commitments}\parentheses*{\text{Order strategy}, \text{Price}} \\
 & \quad \text{Hydrological balance}\parentheses*{\text{Inflow}} \\
 & \quad \text{Electricity production}\parentheses*{\text{Maintenance schedule}} \\
 & \quad \text{Load balance}\parentheses*{\text{Commitments}}
\end{aligned}
\end{equation}

In addition, the producer can take recourse decisions by trading surplus or shortage in a simple balancing market. Both market prices and local water inflows to the reservoirs are considered uncertain.

\subsection{Maintenance scheduling model}
\label{sec:maintenance-scheduling-model_2}

The maintenance scheduling model is mostly derived from the day-ahead model. We also draw inspiration from~\cite{maintenance}. We outline the general structure and also highlight the key details of our model. Further, we sketch how the model is implemented in SPjl.

\subsubsection{General structure}
\label{sec:general-structure_2}

The maintenance scheduling model is defined as a variation of the day-ahead model presented in the previous chapter. We re-cap the essentials to keep this chapter self-contained. In the first stage, we model the day-ahead orders to be submitted to the market. Let $\mathcal{T} = \braces{t_1,\dots,t_{24}}$ denote indices for the $24$-hour horizon of the upcoming day. For simplicity, we refrain from using block orders in this model. We fix a set of hourly price levels $\mathcal{P}_t = \braces{p_{1,t},\dots,p_{P,t}}$ to bid at beforehand as shown in Section~\ref{sec:pricelevels}. We introduce $x^I_t$ and $x^D_{p,t}$ to represent price independent orders and price dependent orders respectively. As per NordPool regulations, the volumes in a price dependent sell order have to be constant or increasing with increasing prices. We enforce this using the following constraints:
\begin{equation*}
  x^D_{i,t} \leq x^D_{i+1,t}, \quad i \in \mathcal{P}, t \in \mathcal{T}.
\end{equation*}
In addition, we model the maintenance schedule by introducing binary decisions. Let $\mathcal{H} = \braces{h_1,\dots,h_{15}}$ index the 15 hydroelectric power stations in Skellefteälven. Next, define $s_{h,t}$ for each power station and hour, where, $s_{h,t} = 1$ indicates that plant $h$ is being maintained during hour $t$. To ensure that maintenance of each plant is finished we include the constraints
\begin{equation*}
  \sum_{t \in \mathcal{T}} s_{h,t} = D_h, \quad h \in \mathcal{H},
\end{equation*}
where $D_h$ is the number of hours required to perform preventive maintenance on plant $h$. Furthermore, we require that the maintenance of each plant is performed during a consecutive period, which can be modeled using
\begin{equation*}
  s_{h,t} - s_{h,t-1} \leq s_{h,t + D_h-1}, \quad h \in \mathcal{H}, t \in \mathcal{T}.
\end{equation*}

In the second stage, we model the order commitments after price realization as well as the production schedule after inflow realization. In contrast to the day-ahead model, we also adhere to the maintenance schedule. We introduce the random variables $\rho^{\omega}_t$, that describe the hourly market price. Let $y^{\omega}_t$ represent the committed hourly volumes. Every hour $t$, the dispatched hourly volumes are determined through linear interpolation:
\begin{equation*}
  y^{\omega}_t = x^I_t + \frac{\rho^{\omega}_t-p_{i,t}}{p_{i+1,t} - p_{i,t}}x^D_{i+1,t} + \frac{p_{i+1,t}-\rho^{\omega}_t}{p_{i+1,t} - p_{i,t}}x^D_{i, t} \quad p_{i,t} \leq \rho^{\omega}_t \leq p_{i+1,t}.
\end{equation*}
Next, we model the production. For each plant and hour, let $Q^{\omega}_{h,s,t}$ and $S^{\omega}_{h,t}$ denote the water discharged and spilled, respectively. The maximum discharge, $\bar{Q}_{h}$ is obtained from Table~\ref{tab:skelleftealven}. Further, let $P^{\omega}_t$ denote the total volume of electricity produced each hour. We employ a piecewise linear approximation of the generation curve of each station. In other words,
\begin{equation*}
  P^{\omega}_t = \sum_{s \in \mathcal{S}} \mu_{h,s}Q^{\omega}_{h,s,t},
\end{equation*}
where $\mu_{h,s}$ is the marginal production equivalent of station $h$ and segment $s \in \mathcal{S}$. See~\cite{obel} for how to estimate the production curve in two piecewise linear segments based on the maximum discharge and capacity. In brief, we set
\begin{equation*}
  \begin{aligned}
    \mu_{h,1} &= \frac{\bar{P}_h}{\bar{Q}_{s}(0.75 + 0.95\cdot 0.25)} \\
    \mu_{h,2} &= 0.95\mu_{h,1}
  \end{aligned}
\end{equation*}
and let $\bar{Q}_{h,1} = 0.75 \bar{Q}_h$ in segment 1 and $\bar{Q}_{h,2} = 0.25 \bar{Q}_h$ in segment 2. To adhere to the maintenance schedule, we require
\begin{equation*}
  Q^{\omega}_{h,s,t} \leq (1-s_{h,t})\bar{Q}^{\omega}_{h,s}, \quad h \in \mathcal{H}, s \in \mathcal{S}, t \in \mathcal{T},
\end{equation*}
so that discharging water from plant $h$ is not allowed if maintenance is taking place. The load balance is given by
\begin{equation*}
  y^{\omega} - P^{\omega}_t = y^{\omega +}_t - y^{\omega -}_t.
\end{equation*}
In each hour in scenario $\omega$, any imbalance between committed volumes and produced volumes is equal to the difference between the imbalance variables $y^{\omega}_{t+}$ and $y^{\omega}_{t-}$. Any shortage $y^{\omega}_{t+}$ is bought from the balancing market, and any surplus $y^{\omega}_{t-}$ is sold to the balancing market. Finally, let $M^{\omega}_{h,t}$ denote the reservoir contents in plant $h$ during hour $t$. The maximum reservoir content, $\bar{M}_{h}$ is obtained from Table~\ref{tab:skelleftealven}. Flow conservation each hour is given by
\begin{equation*}
  \begin{aligned}
    M^{\omega}_{h,t} = \;&M^{\omega}_{h,t-1} \\
    &+ \sum_{i \in \mathcal{Q}_u(h)} \sum_{s \in \mathcal{S}} Q^{\omega}_{i,s,t-\tau_{ih}} + \sum_{i \in \mathcal{S}_u(h)} S^{\omega}_{i,t-\tau_{ih}} + V^{\omega}_h \\
    &- \sum _{s \in \mathcal{S}} Q^{\omega}_{h,s,t} - S^{\omega}_{h,t}
  \end{aligned}
\end{equation*}
Here, $V^{\omega}_h$ are random variables describing the local inflow to each plant. The sets $\mathcal{Q}_u(h)$ and $\mathcal{S}_u(h)$ contain upstream plants where discharge and/or spillage can reach plant $h$ through connecting waterways. Note that the water travel times $\tau_{ih}$ between power stations are included in the incoming flow to each plant, and can be obtained from Table~\ref{tab:skelleftealven}. Internally, this is modeled by introducing auxiliary variables and constraints. Variable limits and the introduced parameters are all included in the deterministic data sets for Skellefteälven given in Table~\ref{tab:skelleftealven}. The revenue from a production schedule satisfying the above relations is given by
\begin{equation*}
  \sum_{t \in \mathcal{T}} \rho^{\omega}_ty^{\omega}_t + \sum_{t \in \mathcal{T}}\alpha_t\rho^{\omega}_ty^{\omega}_{t-} - \beta_t\rho^{\omega}_ty^{\omega}_{t+}.
\end{equation*}
We ignore the water value in this formulation and only optimize the next-day profits. The imbalance volumes are traded at penalized prices, using penalty factors $\alpha_t$ and $\beta_t$, for discouragement. We again use a $15\%$ penalty during peak hours, and $10\%$ otherwise. In summary, a stochastic program modeling the day-ahead problem is in essence given by
\begin{equation} \label{eq:maintenance}
  \begin{aligned}
   \maximize_{s_t,x^I_t, x^D_{i,h}} & \quad \expect[\xi]{\sum_{t \in \mathcal{T}} \rho^{\omega}_ty^{\omega}_t + \sum_{t \in \mathcal{T}} \parentheses*{\alpha_t\rho^{\omega}_ty^{\omega}_{t-} - \beta_t\rho^{\omega}_ty^{\omega}_{t+}}} \\
   \sbj & \quad x^D_{i,t} \leq x^D_{i+1,t}, \quad i \in \mathcal{P}, t \in \mathcal{T} \\
   & \quad \sum_{t \in \mathcal{T}} s_{h,t} = D_h, \quad h \in \mathcal{H} \\
   & \quad s_{h,t} - s_{h,t-1} \leq s_{h,t + D_h-1}, \quad h \in \mathcal{H}, t \in \mathcal{T} \\
   & \quad y^{\omega}_t = x^I_t + \frac{\rho^{\omega}_t-p_{i,t}}{p_{i+1,t} - p_{i,t}}x^D_{i+1,t} + \frac{p_{i+1,t}-\rho^{\omega}_t}{p_{i+1,t} - p_{i,t}}x^D_{i, t}, \quad t \in \mathcal{T} \\
   & \quad P^{\omega}_t = \sum_{s \in \mathcal{S}} \mu_{h,s}Q^{\omega}_{h,s,t}, \quad t \in  \mathcal{T} \\
   & \quad y^{\omega}_t - P^{\omega}_t = y^{\omega}_{t+} - y^{\omega}_{t-}, \quad t \in \mathcal{T} \\
   & \quad \begin{aligned}
     M^{\omega}_{h,t} = \;&M^{\omega}_{h,t-1} \\
     &+ \sum_{i \in \mathcal{Q}_u(h)} \sum_{s \in \mathcal{S}} Q^{\omega}_{i,s,t-\tau_{ih}} + \sum_{i \in \mathcal{S}_u(h)} S^{\omega}_{i,t-\tau_{ih}} \\
     &+ V^{\omega}_h \\
     &- \sum _{s \in \mathcal{S}} Q^{\omega}_{h,s,t} - S^{\omega}_{h,t}, \qquad h \in \mathcal{H}, t \in \mathcal{T} \\
   \end{aligned} \\
   & \quad 0 \leq Q^{\omega}_{h,s,t} \leq (1-s_{h,t})\bar{Q}^{\omega}_{h,s}, \quad h \in \mathcal{H}, s \in \mathcal{S}, t \in \mathcal{T} \\
   & \quad 0 \leq M^{\omega}_{h,t} \leq \bar{M}_h, \quad h \in \mathcal{H}, t \in \mathcal{T} \\
   & \quad y^{\omega}_t \geq 0,\; y^{\omega}_b \geq 0,\; y^{\omega}_{t+} \geq 0,\; y^{\omega}_{t-} \geq 0 \\
   & \quad P^{\omega}_t \geq 0,\; S^{\omega}_{h,t} \geq 0.
 \end{aligned}
\end{equation}

\subsubsection{Model implementation}
\label{sec:model-implementation_2}

We now outline how the maintenance scheduling model~\eqref{eq:maintenance} is formulated in SPjl. To increase readability, we again present an abridged version of the maintenance scheduling model implementation in SPjl and refer to Github~\footnote{\url{https://github.com/martinbiel/HydroModels.jl}} for the full unabridged model.

We re-use the \jlinl{DayAheadScenario} defined in the previous chapter because the maintenance scheduling problem has the same uncertainty model as the day-ahead problem. The maintenance scheduling model definition in SPjl is presented in Listing~\ref{lst:maintenancedef}.

\begin{lstlisting}[language = julia, float, caption = {Maintenance scheduling problem definition in SPjl. The code has been condensed for readability.}, label = {lst:maintenancedef}]
@stochastic_model begin
    @stage 1 begin
        @parameters horizon indices data
        @unpack hours, plants, bids = indices
        @unpack hydrodata, regulations = data
        @decision(model, xᴵ[t in hours] >= 0) # Price-independent orders
        @decision(model, xᴰ[p in bids, t in hours] >= 0) # Price-dependent orders
        # Maintenance schedule
        @decision(model, s[h in plants, t in hours], Bin)
        # Maintenance period should be consecutive
        @constraint(model, maintenance_times[h in plants, t in hours],
                    s[h,t] - s[h,t-1] <= s[h,t+D[h]-1]
        # Ensure maintenance is finished
        @constraint(model, maintenance_finished[h in plants],
                    sum(s[h,t] for t in hours) == D[h])
    end
    @stage 2 begin
        @parameters horizon indices data
        @unpack hours, plants, segments, blocks = indices
        @unpack hydrodata, bidlevels = data
        @uncertain ρ, V from ξ::DayAheadScenario
        @recourse(model, yᴴ[t in hours] >= 0) # Dispatched hourly volumes
        @recourse(model, y⁺[t in hours] >= 0) # Energy shortage
        @recourse(model, y⁻[t in hours] >= 0) # Energy surplus
        @recourse(model, 0 <= Q[h in plants,t in hours] <= Q_max[h]) # Discharge
        @recourse(model, S[h in plants,t in hours] >= 0) # Spillage
        @recourse(model, 0 <= M[h in plants,t in hours] <= M_max[p]) # Water volume
        @recourse(model, P[t in hours] >= 0) # Produced energy
        @expression(model, net_profit, sum(ρ[t]*yᴴ[t] for t in hours)
        @expression(model, intraday_trading,
            sum(penalty(ξ,t)*y⁺[t] - reward(ξ,t)*y⁻[t] for t in hours))
        @objective(model, Max, net_profit - intraday_trading)
        # Bid-dispatch links
        @constraint(model, hourlybids[t in hours],
                   yᴴ[t] == interpolate(ρ[t], bidlevels, xᴰ[t]) + xᴵ[t])
        # Pause production during maintenance hours
        @constraint(model, pause_production[h in plants, t in hours],
                    Q[h,t] <= (1 - s[h,t])*Q_max[h,s]
        # Hydrological balance
        @constraint(model, hydro_constraints[h in plants, t in hours],
            M[h,t] == (t > 1 ? M[h,t-1] : M₀[p]) #
            + sum(Q[i,t-τ] for i in intersect(Qu[p], plants)) # Inflows from
            + sum(S[i,t-τ] for i in intersect(Su[p], plants)) # upstream plants
            + V[p] # Local inflow
            - (Q[h,t] + S[h,t])) # Outflow
        # Production
        @constraint(model, production[t in hours],
            P[t] == sum(μ[p]*Q[h,t] for h in plants))
        # Load balance
        @constraint(model, loadbalance[t in hours],
            yᴴ[t] - P[t] == y⁺[t] - y⁻[t])
        # Water travel time ... (not shown)
    end
end
\end{lstlisting}

\subsubsection{Algorithm details}
\label{sec:algorithm-details_2}

We again employ the sample average approximation (SAA) scheme outlined in Section~\ref{sec:stoch-progr} to compute confidence intervals around the optimal value of the maintenance scheduling problem~\eqref{eq:maintenance}. The sampled instances are again distributed on the 32-core compute node and solved using a distributed L-shaped algorithm. In contrast to the day-ahead problem, the maintenance scheduling problem includes binary decisions in the first stage. Our initial experiments with solving this model reveal that load imbalance has a large impact on performance. As the binary constrained master problem grows in size iterations can be considerably prolonged because the master is not guaranteed to be solved in polynomial time. We therefore configure the L-shaped algorithm to use partial cut aggregation and aggressive cut consolidation to minimize the number of constraints in the master problem. In addition, we use trust-region regularization to speed up convergence. Trust-region regularization was shown to be effective on the day-ahead problem in our numerical benchmarks and the master problem is mixed-integer linear as opposed to mixed-integer quadratic compared to the other regularization schemes. We also warm-start the algorithm with an accelerated subgradient method. The mixed-integer quadratic problem solved in the projection step of the subgradient algorithm does not grow in size so load imbalance is not as prominent. We can therefore make some good initial progress before running the L-shaped procedure.

\subsection{Numerical Experiments}
\label{sec:numerical_experiments_2}

\begin{figure}
  \centering
  \input{maintenance_confidence_intervals.tex}
  \caption{Confidence intervals around the optimal value of the maintenance problem as a function of sample size.}
  \label{fig:maintenance_confidence_intervals}
\end{figure}
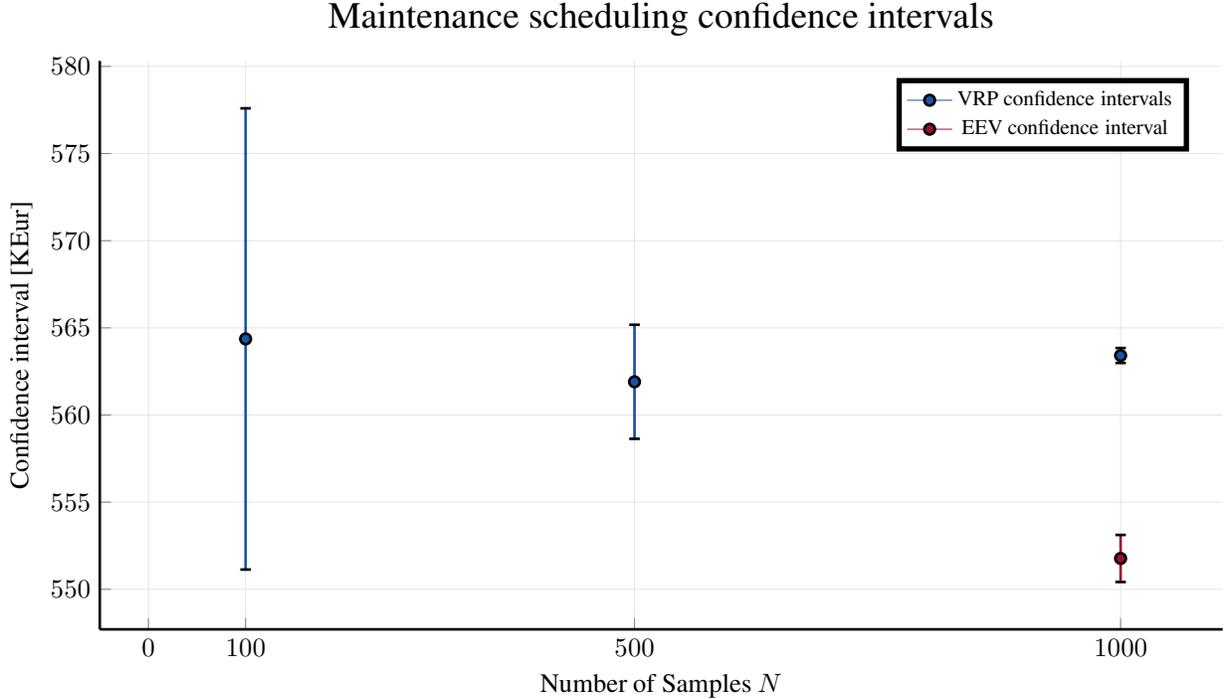

The results of the SAA algorithm is given in Figure~\ref{fig:maintenance_confidence_intervals}. The confidence interval is stabilized at $1000$ samples. We compute a confidence interval around the EEV at this sample size as well. Because there is no overlap between the VRP and EEV, there is a statistically significant VSS within the interval $[1.75\% - 2.38\%]$.

\begin{figure}
  \input{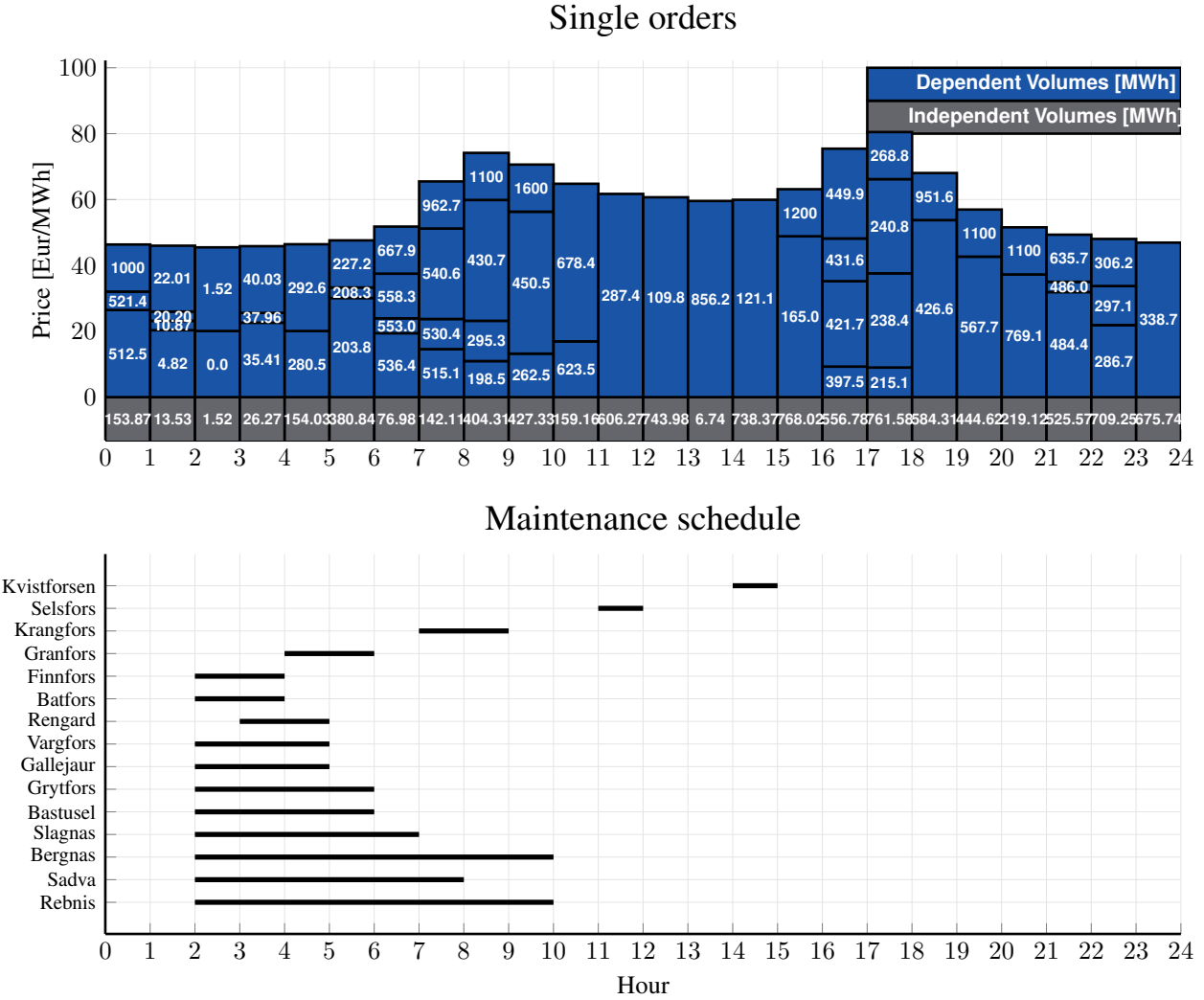}
  \centering
  \caption{Optimal order strategy and maintenance schedule from a $1000$-scenario maintenance scheduling instance.}
  \label{fig:maintenance_strategy}
\end{figure}

We solve a specific sampled instance of $1000$ scenarios and present the resulting strategy in Figure~\ref{fig:maintenance_strategy}. The strategy generated by solving the corresponding expected value problem is shown in Figure~\ref{fig:maintenance_evp_strategy}. The strategy corresponding to the stochastic solution results in a profit of $\num{567.42}$ thousand Euros, while the deterministic strategy yields $\num{537.12}$ thousand Euros. This corresponds to a VSS of $30.3$ thousand Euros, which is $5.3\%$ of the next-day profits. The deterministic strategy opts to start maintenance on all 15 plants in the morning when the electricy price is low. In contrast, the stochastic strategy schedule maintenance of the four downstream plants later in the day, which in expectation yields a larger profit. Similar to the day-ahead problem, the deterministic solution uses a crude set of price dependent orders compared to the more involved orders suggested by the stochastic strategy. The deterministic strategy on average incurs an energy shortage of $\num{2534.80}$ MWh, that must be settled in the intraday market at more expensive prices. The stochastic strategy only incurs an energy shortage of $\num{371.56}$ MWh on average through more conservative planning. This could explain the discrepancy in profits.

\begin{figure}
  \input{maintenance_evp_strategy.tex}
  \centering
  \caption{Optimal order strategy and maintenance schedule from a maintenance scheduling expected value instance.}
  \label{fig:maintenance_evp_strategy}
\end{figure}

\section{Case study 3: Capacity expansion}
\label{sec:capacity-expansion}

In this chapter, we consider capacity expansion of the power stations in Skellefteälven. We use a variation of the RNN forecasters to generate price scenarios over a longer horizon. The capacity expansion problem is then formulated and solved in parallel using SPjl.

\subsection{The capacity expansion problem}
\label{sec:capacity-expansion-problem_3}

Capacity expansion is another common problem formulation for power systems. Due to nuclear phase out by $2040$ and increased demand for electriciy there are predicted challenges in the Swedish power system~\cite{futureprice}. Investing in capacity expansions of existing hydroelectric plants could be a viable approach to meeting future increases in demand. A capacity expansion suggestion for Skellefteälven specifically has recently been proposed by the major owner of the power stations~\cite{skkraft}.

A hydroelectric capacity expansion problem involves specifying an optimal expansion plan that maximizes the expected future profit. We give a brief introduction to the problem in the following.

\subsubsection{Problem setting}
\label{sec:problem-setting_3}

We formulate a capacity expansion planning problem from the perspective of a fictional hydropower producer that owns all 15 power stations in the Swedish river Skellefteälven. The producer is again assumed to be price-taking. We do not consider order strategies to the market and just trade produced energy at market price. In contrast to the day-ahead problem, we schedule the hydropower production over a long planning horizon

When the future market prices and inflows have been realized, the producer optimizes the power production with respect to the price while making use of any extra capacity alloted in the first stage. We again assume that there are no fixed contracts to adhere to. In other words, all electricity production is sold for profit in the market. Any expansion is associated with an investment cost with interest payed back over the planning horizon. A general description of the day-ahead problem is given in~\eqref{eq:generalexpansion}.

\begin{equation} \label{eq:generalexpansion}
\begin{aligned}
 \maximize_{} & \quad \text{Profit} - \text{Expansion cost} \\
 \text{subject to} & \quad \text{Expansion strategy} \\
 & \quad \text{Physical limitations} \\
 & \quad \text{Economic/legal limitations}
\end{aligned}
\end{equation}

Because future market prices and inflows are unknown when planning the capacity expansion, we formulate a two-stage stochastic program to generate an optimal expansion plan. The first-stage decisions are the capacity expansion volumes distributed over the 15 power stations. A general description of the first-stage problem is given in~\eqref{eq:generalfirststage_3}.

\begin{equation} \label{eq:generalfirststage_3}
\begin{aligned}
 \maximize_{\mathclap{\text{Expansion strategy}}} & \quad \expect[]{\text{Revenue}\parentheses*{\text{Expansion strategy}, \text{Price}, \text{Inflow}}} - \text{Cost}\parentheses*{\text{Expansion strategy}} \\
 \text{subject to} & \quad \text{Maximum expansion}
\end{aligned}
\end{equation}

In each second-stage scenario, uncertain parameters are realized and the electricity production is optimized with respect to future profits. A general description of the second stage is given in~\eqref{eq:generalsecondstage_3}
\begin{equation} \label{eq:generalsecondstage_3}
\begin{aligned}
 \maximize_{\mathclap{\text{Production schedule}}} & \quad \text{Profit}\parentheses*{\text{Price}} \\
 \text{subject to} & \quad \text{Hydrological balance}\parentheses*{\text{Inflow}} \\
 & \quad \text{Electricity production}\parentheses*{\text{Expansion strategy}} \\
 & \quad \text{Load balance}
\end{aligned}
\end{equation}

\subsection{Capacity expansion model}
\label{sec:capacity-expansion-model_3}

The capacity expansion model includes the same hydropower scheduling elements as the previous two models. We also draw inspiration from~\cite{Sandstrom2019} for modeling the capacity expansion planning. We outline the general structure and also highlight the key details of our model. The details that coincide with the previous two models are again repeated to keep the chapter self-constained. Further, we sketch how the model is implemented in SPjl.

\subsubsection{General structure}
\label{sec:general-structure_3}

In the first stage, we model the capacity expansion decisions. Let $\mathcal{H} = \braces{h_1,\dots,h_{15}}$ index the 15 hydroelectric power stations in Skellefteälven and introduce $\Delta P_{h}$ to denote the increase in production capacity for each plant $h$. We limit the total expansion to $1000$ MW through
\begin{equation*}
  0 \leq \sum_{h \in \mathcal{H}} \Delta P_h \leq 1000.
\end{equation*}
We use an investment cost estimate of $0.79$ MEur/MW suggested in~\cite{skkraft}. We assume the investment is paid back over a 40 year horizon with a discount rate of $r = 5\%$. As suggested in~\cite{Sandstrom2019}, we can compute the equivalent cost for a shorter time horizon through
\begin{equation*}
  0.79\frac{r^E}{1 - (1-r^E)^{\frac{-40\cdot 365}{T}}} \sum_{h \in \mathcal{H}} \Delta P_h,
\end{equation*}
where $r^E$ is the equivalent interest rate
\begin{equation*}
  r^E = \parentheses*{1+r}^{\frac{T}{365}} - 1,
\end{equation*}
and $T$ is the considered time horizon measured in days.

In the second stage, we model the production schedule after price and inflow realization, while making use of the provisioned increase in capacity. Let $\mathcal{T} = \braces{t_1,\dots,t_{T}}$ denote period indices for the horizon of length $T$. We will consider horizons that are a year or longer, so hourly increments are not computationally tractable. First, we model the production. For each plant and period, let $Q^{\omega}_{h,s,t}$ and $S^{\omega}_{h,t}$ denote the water discharged and spilled, respectively. Further, let $P^{\omega}_t$ denote the total volume of electricity produced each period. We employ a piecewise linear approximation of the generation curve of each station. In other words,
\begin{equation*}
  P^{\omega}_t = \sum_{s \in \mathcal{S}} \mu_{h,s}Q^{\omega}_{h,s,t},
\end{equation*}
where $\mu_{h,s}$ is the marginal production equivalent of station $h$ and segment $s \in \mathcal{S}$. See~\cite{obel} for how to estimate the production curve in two piecewise linear segments based on the maximum discharge and capacity. In brief, we set
\begin{equation*}
  \begin{aligned}
    \mu_{h,1} &= \frac{\bar{P}_h}{\bar{Q}_{s}(0.75 + 0.95\cdot 0.25)} \\
    \mu_{h,2} &= 0.95\mu_{h,1}
  \end{aligned}
\end{equation*}
and let $\bar{Q}_{h,1} = 0.75 \bar{Q}_h$ in segment 1 and $\bar{Q}_{h,2} = 0.25 \bar{Q}_h$ in segment 2. As in~\cite{Sandstrom2019}, we assume that the ratio betweeen the maximum discharge and installed power of each plant is unchanged by the expansions. It follows that the production equivalents are unchanged, while the discharge limits are increased by
\begin{equation*}
  \begin{aligned}
    Q_{h,1,t} &\leq 0.75(\bar{Q}_h + \Delta Q_h) \\
    Q_{h,2,t} &\leq 0.25(\bar{Q}_h + \Delta Q_h),
  \end{aligned}
\end{equation*}
where
\begin{equation*}
  \Delta Q_h = \frac{\bar{Q}_h}{\bar{P}_h} \Delta P_h
\end{equation*}
Finally, let $M^{\omega}_{h,t}$ denote the reservoir contents in plant $h$ during period $t$. Flow conservation each period is given by
\begin{equation*}
  \begin{aligned}
    M^{\omega}_{h,t} = \;&M^{\omega}_{h,t-1} \\
    &+ \sum_{i \in \mathcal{Q}_u(h)} \sum_{s \in \mathcal{S}} Q^{\omega}_{i,s,t-\tau_{ih}} + \sum_{i \in \mathcal{S}_u(h)} S^{\omega}_{i,t-\tau_{ih}} + V^{\omega}_h \\
    &- \sum _{s \in \mathcal{S}} Q^{\omega}_{h,s,t} - S^{\omega}_{h,t}
  \end{aligned}
\end{equation*}
Here, $V^{\omega}_h$ are random variables describing the local inflow to each plant. The sets $\mathcal{Q}_u(h)$ and $\mathcal{S}_u(h)$ contain upstream plants where discharge and/or spillage can reach plant $h$ through connecting waterways. Note that the water travel times $\tau_{ih}$ between power stations are included in the incoming flow to each plant. Internally, this is modeled by introducing auxiliary variables and constraints. Variable limits and the introduced parameters are all included in the deterministic data sets for Skellefteälven given in Table~\ref{tab:skelleftealven} on~\cite{dayahead}. We introduce the random variables $\rho^{\omega}_t$, that describe the market price in period $t$. The revenue from a production schedule satisfying the above relations is then given by
\begin{equation*}
  \sum_{t \in \mathcal{T}} \rho^{\omega}_tP^{\omega}_t.
\end{equation*}
We ignore the water value in this formulation because we consider a long horizon. In summary, a stochastic program modeling the capacity expansion problem is in essence given by
\begin{equation} \label{eq:capacity}
  \begin{aligned}
   \maximize_{\Delta P_h} & \quad \expect[\xi]{\sum_{t \in \mathcal{T}} \rho^{\omega}_tP^{\omega}_t} - 0.79\frac{r^E}{1 - (1-r^E)^{\frac{-b40\cdot 365}{T}}} \sum_{h \in \mathcal{H}} \Delta P_h \\
   \sbj & \quad \Delta P_h \geq 0, \quad h \in \mathcal{H} \\
   & \quad P^{\omega}_t = \sum_{s \in \mathcal{S}} \mu_{h,s}Q^{\omega}_{h,s,t}, \quad t \in  \mathcal{T} \\
   & \quad \begin{aligned}
     M^{\omega}_{h,t} = \;&M^{\omega}_{h,t-1} \\
     &+ \sum_{i \in \mathcal{Q}_u(h)} \sum_{s \in \mathcal{S}} Q^{\omega}_{i,s,t-\tau_{ih}} + \sum_{i \in \mathcal{S}_u(h)} S^{\omega}_{i,t-\tau_{ih}} \\
     &+ V^{\omega}_h \\
     &- \sum _{s \in \mathcal{S}} Q^{\omega}_{h,s,t} - S^{\omega}_{h,t}, \qquad h \in \mathcal{H}, t \in \mathcal{T} \\
   \end{aligned} \\
   & \quad \Delta Q_h = \frac{\bar{Q}_h}{\bar{P}_h}\Delta P_h, \quad h \in \mathcal{H} \\
   & \quad 0 \leq Q^{\omega}_{h,1,t} \leq 0.75(\bar{Q}^{\omega}_{h} + \Delta Q_h), \quad h \in \mathcal{H}, t \in \mathcal{T} \\
   & \quad 0 \leq Q^{\omega}_{h,2,t} \leq 0.25(\bar{Q}^{\omega}_{h} + \Delta Q_h), \quad h \in \mathcal{H}, t \in \mathcal{T} \\
   & \quad 0 \leq M^{\omega}_{h,t} \leq \bar{M}_h, \quad h \in \mathcal{H}, t \in \mathcal{T} \\
   & \quad y^{\omega}_t \geq 0,\; y^{\omega}_b \geq 0,\; y^{\omega}_{t+} \geq 0,\; y^{\omega}_{t-} \geq 0 \\
   & \quad P^{\omega}_t \geq 0,\; S^{\omega}_{h,t} \geq 0.
 \end{aligned}
\end{equation}

\subsubsection{Time resolution}
\label{sec:time-resolution_3}

Because the planning horizon is considerably longer in this model, it is not feasible to consider hourly increments. We adopt the methodology used in~\cite{Sandstrom2019} to change the time resolution of planning problem. We implement an auxiliary $Resolution$ object that scales the relevant quantities in the model based on a predefined number of hours in each period. For example, if a period consists of $24$ hours, water volumes are scaled by $1/24$ and are then measured in $24$-hour equivalents as opposed to $1$-hour equivalents, the marginal production equivalents are scaled by $24$, and water flow times are re-calculated accordingly. After realization of the uncertain parameters, the values in a specific period are calculated through means. In the considered example of $24$ hours in each period, we would use the daily mean price and the daily water inflow.

\subsubsection{Future electricity price}
\label{sec:future-electr-price_3}

Due to challenges in the energy sector, the future electricity price is expected to be much higher than today. A study that considers a future scenario where nuclear power in Sweden is completely phased out in Sweden by 2040 predicts electricy prices $133\%$ above todays levels~\cite{futureprice}. Since our forecasters were trained historical data they do not encompass such effects. To pose a capacity expansion problem with a planning horizon longer than a year we add a rate factor in the price forecast. The predictions in~\cite{futureprice} corresponds to a yearly increase of $4\%$. We use this as a baseline.

\subsubsection{Model implementation}
\label{sec:model-implementation_3}

We now outline how the capacity expansion model~\eqref{eq:capacity} is formulated in SPjl. To increase readability, we again present an abridged version of the capacity expansion model implementation in SPjl and refer to Github~\footnote{\url{https://github.com/martinbiel/HydroModels.jl}} for the full unabridged model.

We first define a scenario data structure to describe the uncertain parameters. In contrast to the day-ahead scenario, the capacity expansion scenario has a variable time horizon. We also define a sampler object that forecasts price and inflow data using the noise-driven RNN forecasters introduced in~\cite{dayahead} over a predefined time horizon. The code is shown in Listing~\ref{lst:capacity_scenariodef}. The sampler begins from a chosen start date and then generates a price curve and inflow sequence stretching over the full horizon. During generation, the date is incremented to make use of the forecasters seasonal prediction capability. In addition, we pick a random yearly rate between $0-4\%$ that grows the forecast for each year in the planning horizon.

\begin{lstlisting}[language = julia, float, caption = {Capacity-expansion scenario definition in SPjl}, label = {lst:capacity_scenariodef}]
@define_scenario CapacityExpansionScenario{T₁,T₂} <: AbstractScenario
    probability::Probability
    ρ::PriceCurve{T₁,Float64}
    Q̃::InflowSequence{T₂,typeof(Skellefteälven),Float64}
end

@sampler RecurrentCapacityExpansionSampler{T₁,T₂}
    date::Date
    plants::PlantCollection
    price_forecaster::Forecaster
    flow_forecaster::Forecaster
    electricity_price_rate::Float64
    horizon::Horizon

    @sample CapacityExpansionScenario begin
        # Generate a random price rate
        price_rate = sampler.electricity_price_rate * rand()
        # Forecast the first price curve
        price_curve = forecast(sampler.price_forecaster, month(sampler.date))
        for d = 1:num_days(sampler.horizon)-1
            # Increment day
            date = sampler.date + Dates.Day(d)
            # Calculate years passed for the correct price ratee
            years_passed = Dates.value(Dates.Year(date) - Dates.Year(sampler.date))
            # Forecast price curve at current day and scale by the price rate
            daily_curve = forecast(sampler.price_forecaster, price_curve[end], month(date))
            daily_curve *= (1 + price_rate)^years_passed
            # Add to price curve
            append!(price_curve, daily_curve)
        end
        # Forecast the first inflows
        flows = forecast(sampler.flow_forecaster, week(sampler.date))
        for w = 1:num_weeks(sampler.horizon)-1
            # Increment week
            date = sampler.date + Dates.Week(w)
            # Forecast flow at current week
            weekly_flows = forecast(sampler.flow_forecaster, flows[:,end], week(date))
            # Add to inflows
            flows = hcat(flows, weekly_flows)
        end
        return CapacityExpansionScenario(PriceCurve(price_curve),
                                         InflowSequence(sampler.plants, flows))
    end
end
\end{lstlisting}

The capacity expansion model definition in SPjl is presented in Listing~\ref{lst:capacitydef}.

\begin{lstlisting}[language = julia, float, caption = {Maintenance scheduling problem definition in SPjl. The code has been condensed for readability.}, label = {lst:capacitydef}]
@stochastic_model begin
    @stage 1 begin
        @parameters horizon indices data
        @unpack hours, plants = indices
        @unpack hydrodata = data
        @decision(model, ΔP[p in plants] >= 0) # Capacity expansion
        # Limit expansion
        @constraint(model, limit_expansion,
                           sum(ΔH[p] for p in plants) <= 1000)
        # Minimize cost
        @objective(model, Max, -equivalent_cost(data, horizon) * ΔH̄)
    end
    @stage 2 begin
        @parameters horizon indices data
        @unpack periods, plants, segments, blocks = indices
        @unpack hydrodata, resolution = data
        @uncertain ρ, V from ξ::CapacityExpansionScenario
        @recourse(model, 0 <= Q[h in plants,t in periods] <= Q_max[h]) # Discharge
        @recourse(model, S[h in plants,t in periods] >= 0) # Spillage
        @recourse(model, 0 <= M[h in plants,t in periods] <= M_max[p]) # Water volume
        @recourse(model, P[t in periods] >= 0) # Produced energy
        @expression(model, net_profit,
            sum(mean_price(resolution, ρ, t)*P[t] for t in periods)
        @objective(model, Max, net_profit)
        # Capacity expansion
        @constraint(model, capacity_expansion[p in plants, s in segments, t in periods],
            Q[p,s,t] <= Q̄(hydrodata, p, s) + (Q̄[p]/P̄[p])ΔH[p])
        # Hydrological balance
        @constraint(model, hydro_constraints[h in plants, t in periods],
            M[h,t] == (t > 1 ? M[h,t-1] : water_volume(resolution, M₀)) #
            + sum(Q[i,t-τ] for i in intersect(Qu[p], plants)) # Inflows from
            + sum(S[i,t-τ] for i in intersect(Su[p], plants)) # upstream plants
            + mean_flow(resolution, V, t, p) # Local inflow
            - (Q[h,t] + S[h,t])) # Outflow
        # Production
        @constraint(model, production[t in periods],
            P[t] == sum(marginal_production(resolution, p, s)*Q[h,t] for h in plants))
        # Water travel time ... (not shown)
    end
end
\end{lstlisting}

\subsubsection{Algorithm details}
\label{sec:algorithm-details_3}

We again employ the sample average approximation (SAA) scheme outline in Section~\ref{sec:stoch-progr} to compute confidence intervals around the optimal value of the capacity expansion problem~\eqref{eq:capacity}. The sampled instances are again distributed on the 32-core compute node and solved using a distributed L-shaped algorithm. In contrast to the models in the previous chapters, solving the subproblems becomes the bottleneck due to the large planning horizon used in the capacity expansion problem. We therefore theorize that limiting the growth of the master using cut aggregation and consolidation will not lead to significant performance improvements. Instead, we elect to use the multi-cut formulation which should yield the best iteration complexity. Otherwise, we employ trust-region regularization to improve convergence.

\subsection{Numerical Experiments}
\label{sec:numerical_experiments_3}

\begin{figure}
  \centering
  \input{capacity_one_year_confidence_intervals.tex}
  \caption{Confidence intervals around optimal value of the capacity expansion problem, with a one year horizon, as a function of sample size.}
  \label{fig:capacity_one_year_confidence_intervals}
\end{figure}
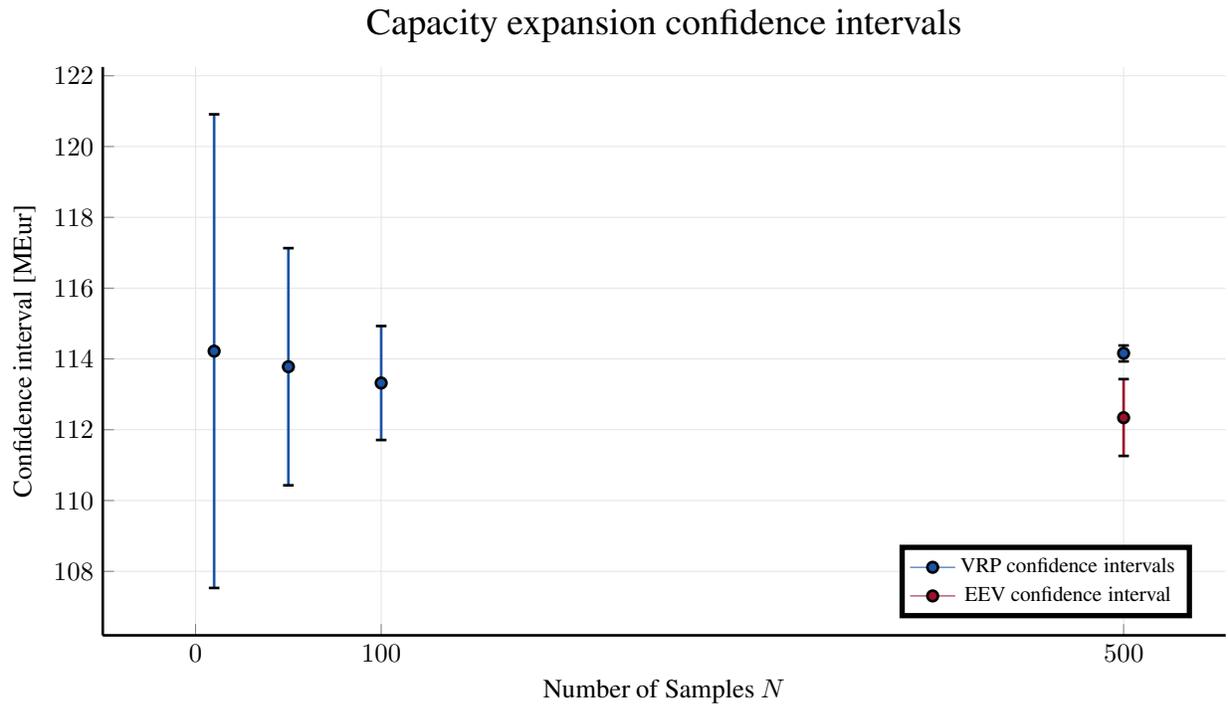

We first consider a planning horizon of one year, where the price forecaster is used as is. The time resolution is set to $24$ hours per period for computational tractability. The results of the SAA algorithm is given in Figure~\ref{fig:capacity_one_year_confidence_intervals}. The confidence interval is stabilized after $500$ samples. We compute a confidence interval around the EEV at this sample size as well. Because there is no overlap between the VRP and EEV, there is a statistically significant VSS within the interval $[0.44\% - 2.73\%]$.

The solution of a $500$-scenario instance is presented in Table~\ref{tab:capacity_one_year_expansion}. A total of $\num{149.29}$ MW is provisioned among the $15$ plants and results in a total profit of $114.39$ million Euros over the one-year planning horizon. If we fix all expansion decisions to zero and re-optimize then the profit without expansion is $111.53$ million Euros. Thus, about $2.86$ million Euros of extra profit is incurred through capacity expansion. The expansion plan generated by solving the deterministic problem instead yields $113.81$ million Euros, which implies a VSS of $0.58$ million Euros, or $0.5\%$. However, we note that the VSS is $20\%$ of the extra profit amassed from the capacity expansion.

\begin{table}[htbp]
  \centering
  \begin{adjustbox}{max width = \textwidth}
    \begin{tabular}{c|S}
      \toprule
      \textbf{Plant} & {Capacity expansion $\Delta P$ [MW]} \\
      \midrule
      Rebnis & 0.0 \\
      Sadva  & 0.0 \\
      Bergnas & 4.75 \\
      Slagnas & 4.21 \\
      Bastusel & 54.87 \\
      Grytfors & 20.15 \\
      Gallejaur & 0.0 \\
      Vargfors & 0.0 \\
      Rengard & 15.13 \\
      Batfors & 5.23 \\
      Finnfors & 2.78 \\
      Granfors & 12.53 \\
      Krangfors & 19.49 \\
      Selsfors & 3.17 \\
      Kvistforsen & 6.98 \\
      \bottomrule
    \end{tabular}
  \end{adjustbox}
  \caption{Capacity expansion strategy obtained from solving a $500$-scenario capacity expansion problem with a one-year planning horizon.}
  \label{tab:capacity_one_year_expansion}
\end{table}

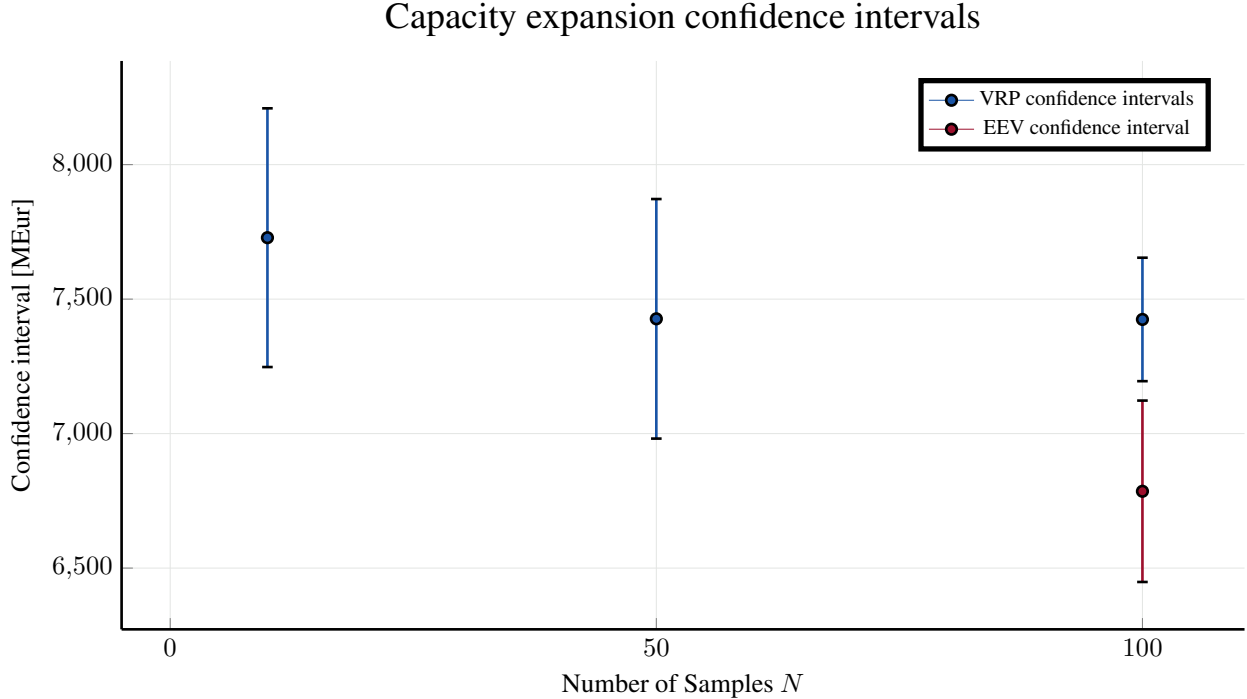
\begin{figure}
  \centering
  \input{capacity_20_year_confidence_intervals.tex}
  \caption{Confidence intervals around optimal value of the capacity expansion problem, with a 20-year horizon, as a function of sample size.}
  \label{fig:capacity_20_year_confidence_intervals}
\end{figure}

Next, we consider a long planning horizon of $20$ years. We set the time resolution to $120$ hours per period for computational tractability. In other words, the planning is constrained to intervals of five days. This decreases the accuracy of the model, but the model size will reach the memory limit of our hardware at $100$ scenarios already at this resolution. Confidence intervals around the optimal value are presented in Figure~\ref{fig:capacity_20_year_confidence_intervals}. The largest sample size we can use in our hardware setup is $N = 100$, which yields a fairly stable confidence interval. We compute a confidence interval around the EEV at this sample size as well. Because there is no overlap between the VRP and EEV, there is a statistically significant VSS within the interval $[0.97\% - 16.24\%]$. The VSS is considerably higher than when using the shorter planning horizon.

The solution of a $100$-scenario instance is presented in Table~\ref{tab:capacity_20_year_expansion}. A total of $\num{437.60}$ MW is provisioned among the $15$ plants and results in a total profit of $7649.18$ million Euros over the 20-year planning horizon. If we fix all expansion decisions to zero and re-optimize then the profit without expansion is $6199.90$ million Euros. Thus, about $1449.28$ million Euros of extra profit is incurred through capacity expansion. The expansion plan is about three times as aggressive as the one-year plan, and the gain from expansion is far greater. This is probably an effect of predicting a larger future electricity price. The expansion plan generated by solving the deterministic problem instead yields $7039.73$ million Euros, which implies a VSS of $609.45$ million Euros, or $7.97\%$. As a percentage of the expansion gain, the VSS is $42\%$. The results could be compared to the $363$ MW capacity expansion suggested in~\cite{skkraft}.

\begin{table}[htbp]
  \centering
  \begin{adjustbox}{max width = \textwidth}
    \begin{tabular}{c|S}
      \toprule
      \textbf{Plant} & {Capacity expansion $\Delta P$ [MW]} \\
      \midrule
      Rebnis & 0.0 \\
      Sadva  & 0.0 \\
      Bergnas & 0.0 \\
      Slagnas & 0.0 \\
      Bastusel & 0.0 \\
      Grytfors & 12.90 \\
      Gallejaur & 0.0 \\
      Vargfors & 0.0 \\
      Rengard & 34.86 \\
      Batfors & 37.53 \\
      Finnfors & 46.04 \\
      Granfors & 53.48 \\
      Krangfors & 82.76 \\
      Selsfors & 52.98 \\
      Kvistforsen & 117.04 \\
      \bottomrule
    \end{tabular}
  \end{adjustbox}
  \caption{Capacity expansion strategy obtained from solving a $100$-scenario capacity expansion problem with a 20-year planning horizon.}
  \label{tab:capacity_20_year_expansion}
\end{table}

\section{Concluding remarks}
\label{sec:concluding-remarks-1}

We have utilized the \jlinl{StochasticPrograms.jl} framework to pose and solve three large-scale planning problems related to hydropower production in the Swedish river Skellefteälven.

First, we formulated and solved a day-ahead planning problem in SPjl. The problem is formulated from the perspective of a hydropower producer participating in a day-ahead market. Both the electricity price and local water inflows are uncertain. We provided a detailed model formulation of the day-ahead problem and explained how to the model can be implemented in SPjl. We used the SAA algorithm to compute tight confidence intervals around the optimal solution of the day-ahead problems. This approach becomes tractable through the parallel capabilities of SPjl. The tight confidence intervals allow us to conclude that the gain from solving the stochastic program is statistically significant.

Next, we formulated and solved a stochastic program for determining optimal day-ahead order strategies in conjuction with a maintenance schedule. The SAA algorithm again yields tight confidence intervals around the stochastic solution, and the resulting VSS is larger than the day-ahead formulation. We argue that the intricacies involved with coordinating the maintenance schedule with the order strategies increases the value of a stochastic approach. The maintenance periods we have used are not necessarily accurate, but we believe that the strong experimental results indicate that the model could be a viable approach for optimally planning preventive maintenance while participating in a deregulated electricy market.

Finally, we considered a capacity expansion problem with a long planning horizon. The same methodology is applied as when solving the first two hydropower problems. However, the planning horizon is considerably longer, from one year up to 20 years compared to the 24 hour horizon used in the first two problems. We therefore utilized a coarser time resolution and smaller sample sizes to not exceed the memory capacity of our hardware setup. The SAA algorithm still produces relatively tight confidence intervals with a statistically significant VSS both when the planning horizon is one year and when it is 20 years. Also, the relative significance of the VSS is much greater when comparing to the extra profits incurred from the capacity expansion instead of the total profit. In brief, the VSS is more significant over longer planning horizons, even with a coarse time resolution. This suggests that stochastic programming is a viable tool for this type of planning problem.

In summary, the three case studies provide a solid proof-of-concept of the \jlinl{StochasticPrograms.jl} framework. Because the uncertainty modeling is decoupled from the optimization modeling, we can effortlessly re-use the forecasting technique presented~\cite{dayahead} to sample scenarios in all three applications. The distributed capabilities allow us to solve large-scale instances. Moreover, the specialized algorithms avaiable in SPjl are used to efficiently solve sampled instances. Consequently, we can afford to run SAA algorithms until tight confidence intervals are obtained. This allows us to be precise when discussing the generated solutions and prove that the value of the stochastic solution is statistically significant.

\bibliographystyle{unsrt}
\bibliography{references}

\end{document}

%% file: rnnarch.tex
\begin{tikzpicture}
  \node at (0,0) (anchor) {};

  \node[left of=anchor,xshift=-100pt] (gaussian) {$\tilde{w}$};
  \node[draw,right of = gaussian, node distance = 60pt] (init) {Initializer};
  \draw[->] (gaussian) -- (init);
  \node[below of = init, node distance = 20pt] (seasonal) {$s$};
  \draw[->] (seasonal) -- (init);
  \node[draw, right of = init, node distance = 60pt] (seq) {RNN};
  \draw[->] (init) -- node[pos=0.5, above] {$x_0$} (seq);
  \node[below of = seq, node distance = 20pt] (seasonalseq) {$s$};
  \draw[->] (seasonalseq) -- (seq);
  \node[right of = seq, node distance = 60pt] (output) {$x_t$};
  \draw[->] (seq) -- node[pos = 0.5] (sr0) {} (output);
  \node[above of = sr0] (sr1) {};
  \node[above of = seq] (sr2) {};
  \draw[->] (sr0.center) -- (sr1.south) -- (sr2.south) -- (seq);
\end{tikzpicture}

%% file: single_order_ex.tex
\begin{tikzpicture}[]
  \begin{axis}[
    width = {\textwidth},
    height = {0.4\textheight},
    xlabel = {Order Volume [Mwh/h]},
    ylabel = {Price [Eur/Mwh]},
    title = {Order curve},
    unbounded coords=jump,
    scaled x ticks = false,
    xmajorgrids = true,
    xmin = -50,
    xtick align = inside,
    xticklabel style = {font = {\fontsize{10 pt}{10.4 pt}\selectfont}, rotate = 0.0},
    x grid style = {color = kth-lightgray,
      line width = 0.25,
      solid},
    axis x line* = left,
    x axis line style = {line width = 1,
      solid},
    scaled y ticks = false,
    ylabel style = {font = {\fontsize{10 pt}{14.3 pt}\selectfont}, rotate = 0.0},
    ymajorgrids = true,
    ytick align = inside,
    yticklabel style = {font = {\fontsize{10 pt}{10.4 pt}\selectfont}, rotate = 0.0},
    y grid style = {color = kth-lightgray,
      line width = 0.25,
      solid},axis y line* = left,y axis line style = {line width = 1,
      solid},
    xshift = 0mm,
    yshift = {0.2\textheight},
    title style = {font = {\fontsize{14 pt}{18.2 pt}\selectfont}, rotate = 0.0},
    legend pos = {north west},
    legend style = {line width = 2, solid,font = {\fontsize{8 pt}{10.4 pt}\selectfont}},colorbar style={title=}]
    \addplot+ [line width = 2, color = black,
    solid,mark = none,
    mark size = 2.0,
    mark options = {
      color = black,
      fill = black,
      line width = 1,
      rotate = 0,
      solid
    },forget plot]coordinates {
      (0.0, 0.0)
      (0.0, 0.0)
      (0.0, 20.89356016358743)
      (636.7057290307187, 20.89356016358743)
      (636.7057290307187, 20.89356016358743)
      (636.7057290307187, 23.644073132559704)
      (636.7057290307187, 23.644073132559704)
      (636.7057290307187, 26.39458610153198)
      (636.7057290307187, 29.14509907050426)
      (680.0388233393326, 29.14509907050426)
      (680.0388233393326, 29.14509907050426)
      (680.0388233393326, 31.895612039476532)
      (680.0388233393326, 31.895612039476532)
      (680.0388233393326, 34.646125008448806)
    };
    \addplot+[draw=none,
    color = kth-darkgray,
    line width = 0,
    solid,mark = *,
    mark size = 2.0,
    mark options = {
      color = black,
      fill = kth-darkgray,
      line width = 1,
      rotate = 0,
      solid
    }] coordinates {
      (100.0, 0.0)
    };
    \addlegendentry{Price Independent Order}
    \addplot+[draw=none,
    line width = 0,
    color = kth-blue,
    solid,mark = *,
    mark size = 2.0,
    mark options = {
      color = black,
      fill = kth-blue,
      line width = 1,
      rotate = 0,
      solid
    }] coordinates {
      (636.7057290307187, 20.89356016358743)
      (636.7057290307187, 23.644073132559704)
      (636.7057290307187, 26.39458610153198)
      (680.0388233393326, 29.14509907050426)
      (680.0388233393326, 31.895612039476532)
    };
    \addlegendentry{Price Dependent Order}
    \addplot+ [line width = 1,
    color = black,
    dashed,mark = none,
    forget plot]coordinates {
      (636.7057290307187, 26.39458610153198)
      (680.0388233393326, 29.14509907050426)
    };
    \addplot+[draw=none,
    line width = 0,
    solid,mark = diamond*,
    mark size = 2.0,
    mark options = {
      color = black,
      fill = kth-green,
      line width = 1,
      rotate = 0,
      solid
    }] coordinates {
      (661.9983026893217, 28.0)
    };
    \addlegendentry{Trading Outcome}
  \end{axis}
\end{tikzpicture}

%% file: price_levels.tex
\begin{tikzpicture}[]
  \begin{axis}[
    width = {\textwidth},
    height = {0.4\textheight},
    xlabel = {Hour},
    ylabel = {Price [Eur]},
    title = {Bidding levels},
    unbounded coords=jump,
    scaled x ticks = false,
    xmajorgrids = true,
    xmin = 0.5,
    xmax = 24.5,
    xtick align = inside,
    xticklabel style = {font = {\fontsize{10 pt}{10.4 pt}\selectfont}, rotate = 0.0},
    x grid style = {color = kth-lightgray,
      line width = 0.25,
      solid},
    axis x line* = left,
    x axis line style = {line width = 1,
      solid},
    scaled y ticks = false,
    ylabel style = {font = {\fontsize{10 pt}{14.3 pt}\selectfont}, rotate = 0.0},
    ymajorgrids = true,
    ytick align = inside,
    yticklabel style = {font = {\fontsize{10 pt}{10.4 pt}\selectfont}, rotate = 0.0},
    y grid style = {color = kth-lightgray,
      line width = 0.25,
      solid},axis y line* = left,y axis line style = {line width = 1,
      solid},
    xshift = 0mm,
    yshift = {0.2\textheight},
    title style = {font = {\fontsize{14 pt}{18.2 pt}\selectfont}, rotate = 0.0},
    legend pos = {north west},
    legend style = {line width = 2, solid,font = {\fontsize{8 pt}{10.4 pt}\selectfont}},colorbar style={title=}]
    \addplot+ [color = kth-blue,
    line width = 2,
    solid,mark = none,
    mark size = 2.0,
    mark options = {
      color = black,
      fill = kth-blue,
      line width = 1,
      rotate = 0,
      solid
    }]coordinates {
      (1.0, 26.377586431503296)
      (2.0, 25.91467671966553)
      (3.0, 25.592450896263124)
      (4.0, 25.481782537460326)
      (5.0, 26.004895790100097)
      (6.0, 26.825305461883545)
      (7.0, 28.385489170074464)
      (8.0, 32.749079216003416)
      (9.0, 35.15360792827606)
      (10.0, 34.47467072868347)
      (11.0, 33.52035451698303)
      (12.0, 32.583070892333986)
      (13.0, 32.102346574783326)
      (14.0, 31.740470558166503)
      (15.0, 31.701566961288453)
      (16.0, 32.342629470825194)
      (17.0, 35.07680814361572)
      (18.0, 37.43803865814209)
      (19.0, 35.123736673355104)
      (20.0, 31.9463049659729)
      (21.0, 29.970932359695436)
      (22.0, 28.784440793991088)
      (23.0, 27.81591905593872)
      (24.0, 26.84460446166992)
    };
    \addlegendentry{Expected price}
    \addplot+[draw=none, color = kth-red,
    line width = 0,
    solid,mark = *,
    mark size = 2.0,
    mark options = {
      color = black,
      fill = kth-red,
      line width = 1,
      rotate = 0,
      solid
    }] coordinates {
      (1.0, 20.99341396253017)
      (2.0, 20.45733220765932)
      (3.0, 20.239489859538217)
      (4.0, 19.737804863058187)
      (5.0, 20.152802829257215)
      (6.0, 20.470774425740707)
      (7.0, 19.26656043356128)
      (8.0, 13.945209308787053)
      (9.0, 11.369126198144443)
      (10.0, 13.403605200008595)
      (11.0, 16.49702675391818)
      (12.0, 17.865773269678996)
      (13.0, 17.988791440942336)
      (14.0, 18.523584451985485)
      (15.0, 18.15928997909782)
      (16.0, 16.431537024584017)
      (17.0, 9.728789733695685)
      (18.0, 9.315750056061752)
      (19.0, 16.391890863433396)
      (20.0, 21.133446449558875)
      (21.0, 22.586666222843522)
      (22.0, 22.63350740400468)
      (23.0, 22.273971688418406)
      (24.0, 21.415767819469757)
    };
    \addlegendentry{Price levels}
    \addplot+[draw=none, color = kth-red,
    line width = 0,
    solid,mark = *,
    mark size = 2.0,
    mark options = {
      color = black,
      fill = kth-red,
      line width = 1,
      rotate = 0,
      solid
    },forget plot] coordinates {
      (1.0, 23.685500197016733)
      (2.0, 23.186004463662425)
      (3.0, 22.91597037790067)
      (4.0, 22.609793700259257)
      (5.0, 23.078849309678656)
      (6.0, 23.648039943812126)
      (7.0, 23.826024801817873)
      (8.0, 23.347144262395233)
      (9.0, 23.26136706321025)
      (10.0, 23.93913796434603)
      (11.0, 25.008690635450606)
      (12.0, 25.22442208100649)
      (13.0, 25.04556900786283)
      (14.0, 25.132027505075996)
      (15.0, 24.930428470193135)
      (16.0, 24.387083247704606)
      (17.0, 22.402798938655703)
      (18.0, 23.37689435710192)
      (19.0, 25.75781376839425)
      (20.0, 26.539875707765887)
      (21.0, 26.27879929126948)
      (22.0, 25.708974098997885)
      (23.0, 25.044945372178564)
      (24.0, 24.130186140569837)
    };
    \addplot+[draw=none, color = kth-red,
    line width = 0,
    solid,mark = *,
    mark size = 2.0,
    mark options = {
      color = black,
      fill = kth-red,
      line width = 1,
      rotate = 0,
      solid
    },forget plot] coordinates {
      (1.0, 26.377586431503296)
      (2.0, 25.91467671966553)
      (3.0, 25.592450896263124)
      (4.0, 25.481782537460326)
      (5.0, 26.004895790100097)
      (6.0, 26.825305461883545)
      (7.0, 28.385489170074464)
      (8.0, 32.749079216003416)
      (9.0, 35.15360792827606)
      (10.0, 34.47467072868347)
      (11.0, 33.52035451698303)
      (12.0, 32.583070892333986)
      (13.0, 32.102346574783326)
      (14.0, 31.740470558166503)
      (15.0, 31.701566961288453)
      (16.0, 32.342629470825194)
      (17.0, 35.07680814361572)
      (18.0, 37.43803865814209)
      (19.0, 35.123736673355104)
      (20.0, 31.9463049659729)
      (21.0, 29.970932359695436)
      (22.0, 28.784440793991088)
      (23.0, 27.81591905593872)
      (24.0, 26.84460446166992)
    };
    \addplot+[draw=none, color = kth-red,
    line width = 0,
    solid,mark = *,
    mark size = 2.0,
    mark options = {
      color = black,
      fill = kth-red,
      line width = 1,
      rotate = 0,
      solid
    },forget plot] coordinates {
      (1.0, 29.06967266598986)
      (2.0, 28.643348975668633)
      (3.0, 28.268931414625577)
      (4.0, 28.353771374661395)
      (5.0, 28.930942270521538)
      (6.0, 30.002570979954964)
      (7.0, 32.944953538331056)
      (8.0, 42.1510141696116)
      (9.0, 47.04584879334187)
      (10.0, 45.01020349302091)
      (11.0, 42.03201839851546)
      (12.0, 39.94171970366148)
      (13.0, 39.15912414170382)
      (14.0, 38.34891361125701)
      (15.0, 38.47270545238377)
      (16.0, 40.29817569394578)
      (17.0, 47.75081734857574)
      (18.0, 51.499182959182264)
      (19.0, 44.48965957831596)
      (20.0, 37.35273422417991)
      (21.0, 33.66306542812139)
      (22.0, 31.85990748898429)
      (23.0, 30.58689273969888)
      (24.0, 29.559022782770004)
    };
    \addplot+[draw=none, color = kth-red,
    line width = 0,
    solid,mark = *,
    mark size = 2.0,
    mark options = {
      color = black,
      fill = kth-red,
      line width = 1,
      rotate = 0,
      solid
    },forget plot] coordinates {
      (1.0, 31.761758900476423)
      (2.0, 31.372021231671738)
      (3.0, 30.94541193298803)
      (4.0, 31.225760211862465)
      (5.0, 31.85698875094298)
      (6.0, 33.17983649802638)
      (7.0, 37.50441790658765)
      (8.0, 51.552949123219776)
      (9.0, 58.93808965840768)
      (10.0, 55.54573625735834)
      (11.0, 50.54368228004789)
      (12.0, 47.300368514988975)
      (13.0, 46.215901708624315)
      (14.0, 44.95735666434752)
      (15.0, 45.243843943479085)
      (16.0, 48.25372191706637)
      (17.0, 60.42482655353576)
      (18.0, 65.56032726022244)
      (19.0, 53.85558248327681)
      (20.0, 42.759163482386924)
      (21.0, 37.35519849654735)
      (22.0, 34.9353741839775)
      (23.0, 33.357866423459036)
      (24.0, 32.273441103870084)
    };
  \end{axis}
\end{tikzpicture}

%% file: dayahead_confidence_intervals.tex
\begin{tikzpicture}[]
  \begin{axis}[
    width = {\textwidth},
    height = {0.4\textheight},
    xlabel = {Number of Samples $N$},
    ylabel = {Confidence interval [MEur]},
    title = {Day-ahead confidence intervals},
    unbounded coords=jump,
    scaled x ticks = false,
    xmajorgrids = true,
    xmin = -50,
    xtick = {0.0,500.0,1000.0,1500.0,2000.0},
    xticklabels = {$0$,$500$,$1000$,$1500$,$2000$},
    xtick align = inside,
    xticklabel style = {font = {\fontsize{10 pt}{10.4 pt}\selectfont}, rotate = 0.0},
    x grid style = {color = kth-lightgray,
      line width = 0.25,
      solid},
    axis x line* = left,
    x axis line style = {line width = 1,
      solid},
    scaled y ticks = false,
    ylabel style = {font = {\fontsize{10 pt}{14.3 pt}\selectfont}, rotate = 0.0},
    ymajorgrids = true,
    ytick align = inside,
    yticklabel style = {font = {\fontsize{10 pt}{10.4 pt}\selectfont}, rotate = 0.0},
    y grid style = {color = kth-lightgray,
      line width = 0.25,
      solid},axis y line* = left,y axis line style = {line width = 1,
      solid},
    xshift = 0mm,
    title style = {font = {\fontsize{14 pt}{18.2 pt}\selectfont}, rotate = 0.0},
    legend pos = {south east},
    legend style = {line width = 2, solid,font = {\fontsize{8 pt}{10.4 pt}\selectfont}},colorbar style={title=}]
    \addplot+[draw=none, color = kth-blue,
    line width = 0,
    solid,mark = *,
    mark size = 2.0,
    mark options = {
      color = black,
      fill = kth-blue,
      line width = 1,
      rotate = 0,
      solid
    }] coordinates {
      (10.0, 5.7204498019063715)
      (100.0, 5.731231656548085)
      (1000.0, 5.737249943856671)
      (2000.0, 5.737140981226794)
    };
    \addlegendentry{VRP confidence intervals}
    \addplot+[draw=none, color = kth-red,
    line width = 0,
    solid,mark = *,
    mark size = 2.0,
    mark options = {
      color = black,
      fill = kth-red,
      line width = 1,
      rotate = 0,
      solid
    }] coordinates {
      (2000.0, 5.729554235)
    };
    \addlegendentry{EEV confidence intervals}
    \addplot+ [color = kth-blue,
    line width = 1,
    solid,mark = -,
    mark size = 2.0,
    mark options = {
      color = black,
      line width = 1,
      rotate = 0,
      solid
    },forget plot]coordinates {
      (10.0, 5.701168201217919)
      (10.0, 5.739731402594825)
    };
    \addplot+ [color = kth-blue,
    line width = 1,
    solid,mark = -,
    mark size = 2.0,
    mark options = {
      color = black,
      line width = 1,
      rotate = 0,
      solid
    },forget plot]coordinates {
      (100.0, 5.722225596706217)
      (100.0, 5.740237716389952)
    };
    \addplot+ [color = kth-blue,
    line width = 1,
    solid,mark = -,
    mark size = 2.0,
    mark options = {
      color = black,
      line width = 1,
      rotate = 0,
      solid
    },forget plot]coordinates {
      (1000.0, 5.734234395853527)
      (1000.0, 5.740265491859815)
    };
    \addplot+ [color = kth-blue,
    line width = 1,
    solid,mark = -,
    mark size = 2.0,
    mark options = {
      color = black,
      line width = 1,
      rotate = 0,
      solid
    },forget plot]coordinates {
      (2000.0, 5.7348487136523975)
      (2000.0, 5.739433248801191)
    };
    \addplot+ [color = kth-blue,
    line width = 1,
    solid,mark = -,
    mark size = 2.0,
    mark options = {
      color = black,
      line width = 1,
      rotate = 0,
      solid
    },forget plot]coordinates {
      (2000.0, 5.72757388)
      (2000.0, 5.73153459)
    };
  \end{axis}
\end{tikzpicture}

%% file: day-ahead_strategy.tex
\begin{tikzpicture}[]
  \begin{axis}[
    width = {\textwidth},
    height = {0.3\textheight},
    xlabel = {},
    ylabel = {Price [Eur/MWh]},
    title = {Single orders},
    unbounded coords=jump,
    scaled x ticks = false,
    xmajorgrids = true,
    xmin = 0,
    xmax = 24,
    ymin = {-13.112593618627468},
    ymax = {102.22033269617995},
    xtick align = inside,
    xtick = {0.0,1.0,2.0,3.0,4.0,5.0,6.0,7.0,8.0,9.0,10.0,11.0,12.0,13.0,14.0,15.0,16.0,17.0,18.0,19.0,20.0,21.0,22.0,23.0,24.0},
    xticklabels = {$0$,$1$,$2$,$3$,$4$,$5$,$6$,$7$,$8$,$9$,$10$,$11$,$12$,$13$,$14$,$15$,$16$,$17$,$18$,$19$,$20$,$21$,$22$,$23$,$24$},
    xticklabel style = {font = {\fontsize{10 pt}{10.4 pt}\selectfont}, rotate = 0.0},
    x grid style = {color = kth-lightgray,
      line width = 0.25,
      solid},
    axis x line* = left,
    x axis line style = {line width = 1,
      solid},
    scaled y ticks = false,
    ylabel style = {yshift = -5pt, font = {\fontsize{10 pt}{14.3 pt}\selectfont}, rotate = 0.0},
    ymajorgrids = true,
    ytick align = inside,
    yticklabel style = {font = {\fontsize{10 pt}{10.4 pt}\selectfont}, rotate = 0.0},
    y grid style = {color = kth-lightgray,
      line width = 0.25,
      solid},axis y line* = left,y axis line style = {line width = 1,
      solid},
    xshift = 0mm,
    yshift = {0.3\textheight},
    title style = {font = {\fontsize{14 pt}{18.2 pt}\selectfont}, rotate = 0.0},
    legend pos = {north west},
    legend style = {line width = 2, solid,font = {\fontsize{8 pt}{10.4 pt}\selectfont}},colorbar style={title=}]
    \addplot+ [line width = 1,
    solid,mark = none,
    color = black,
    mark size = 2.0,
    mark options = {
      color = black,
      fill = kth-darkgray,
      line width = 1,
      rotate = 0,
      solid
    },fill = kth-darkgray, forget plot]coordinates {
      (0.0, -13.112593618627468)
      (1.0, -13.112593618627468)
      (1.0, 0.0)
      (0.0, 0.0)
      (0.0, -13.112593618627468)
    };
    \addplot+ [line width = 1,
    solid,mark = none,
    color = black,
    mark size = 2.0,
    mark options = {
      fill = kth-darkgray,
      line width = 1,
      rotate = 0,
      solid
    },fill = kth-darkgray, forget plot]coordinates {
      (17.0, 80)
      (24.0, 80)
      (24.0, 90)
      (17.0, 90)
      (17.0, 80)
    };
    \addplot+ [line width = 1,
    solid,mark = none,
    color = black,
    mark size = 2.0,
    mark options = {
      color = black,
      fill = kth-blue,
      line width = 1,
      rotate = 0,
      solid
    },fill = kth-blue, forget plot]coordinates {
      (1.0, 0.0)
      (2.0, 0.0)
      (2.0, 44.491096562153785)
      (1.0, 44.491096562153785)
      (1.0, 0.0)
    };
    \addplot+ [line width = 1,
    solid,mark = none,
    color = black,
    mark size = 2.0,
    mark options = {
      color = black,
      fill = kth-blue,
      line width = 1,
      rotate = 0,
      solid
    },fill = kth-blue, forget plot]coordinates {
      (2.0, 0.0)
      (3.0, 0.0)
      (3.0, 20.304737464236357)
      (2.0, 20.304737464236357)
      (2.0, 0.0)
    };
    \addplot+ [line width = 1,
    solid,mark = none,
    color = black,
    mark size = 2.0,
    mark options = {
      color = black,
      fill = kth-blue,
      line width = 1,
      rotate = 0,
      solid
    },fill = kth-blue, forget plot]coordinates {
      (2.0, 20.304737464236357)
      (3.0, 20.304737464236357)
      (3.0, 43.986893132159054)
      (2.0, 43.986893132159054)
      (2.0, 20.304737464236357)
    };
    \addplot+ [line width = 1,
    solid,mark = none,
    color = black,
    mark size = 2.0,
    mark options = {
      color = black,
      fill = kth-blue,
      line width = 1,
      rotate = 0,
      solid
    },fill = kth-blue, forget plot]coordinates {
      (3.0, 0.0)
      (4.0, 0.0)
      (4.0, 19.817582353470804)
      (3.0, 19.817582353470804)
      (3.0, 0.0)
    };
    \addplot+ [line width = 1,
    solid,mark = none,
    color = black,
    mark size = 2.0,
    mark options = {
      color = black,
      fill = kth-blue,
      line width = 1,
      rotate = 0,
      solid
    },fill = kth-blue, forget plot]coordinates {
      (3.0, 19.817582353470804)
      (4.0, 19.817582353470804)
      (4.0, 31.14522571003628)
      (3.0, 31.14522571003628)
      (3.0, 19.817582353470804)
    };
    \addplot+ [line width = 1,
    solid,mark = none,
    color = black,
    mark size = 2.0,
    mark options = {
      color = black,
      fill = kth-blue,
      line width = 1,
      rotate = 0,
      solid
    },fill = kth-blue, forget plot]coordinates {
      (3.0, 31.14522571003628)
      (4.0, 31.14522571003628)
      (4.0, 44.25781932866374)
      (3.0, 44.25781932866374)
      (3.0, 31.14522571003628)
    };
    \addplot+ [line width = 1,
    solid,mark = none,
    color = black,
    mark size = 2.0,
    mark options = {
      color = black,
      fill = kth-blue,
      line width = 1,
      rotate = 0,
      solid
    },fill = kth-blue, forget plot]coordinates {
      (4.0, 0.0)
      (5.0, 0.0)
      (5.0, 20.203211391501167)
      (4.0, 20.203211391501167)
      (4.0, 0.0)
    };
    \addplot+ [line width = 1,
    solid,mark = none,
    color = black,
    mark size = 2.0,
    mark options = {
      color = black,
      fill = kth-blue,
      line width = 1,
      rotate = 0,
      solid
    },fill = kth-blue, forget plot]coordinates {
      (4.0, 20.203211391501167)
      (5.0, 20.203211391501167)
      (5.0, 25.97668363380432)
      (4.0, 25.97668363380432)
      (4.0, 20.203211391501167)
    };
    \addplot+ [line width = 1,
    solid,mark = none,
    color = black,
    mark size = 2.0,
    mark options = {
      color = black,
      fill = kth-blue,
      line width = 1,
      rotate = 0,
      solid
    },fill = kth-blue, forget plot]coordinates {
      (4.0, 25.97668363380432)
      (5.0, 25.97668363380432)
      (5.0, 44.86274949473494)
      (4.0, 44.86274949473494)
      (4.0, 25.97668363380432)
    };
    \addplot+ [line width = 1,
    solid,mark = none,
    color = black,
    mark size = 2.0,
    mark options = {
      color = black,
      fill = kth-blue,
      line width = 1,
      rotate = 0,
      solid
    },fill = kth-blue, forget plot]coordinates {
      (5.0, 0.0)
      (6.0, 0.0)
      (6.0, 32.86009143933869)
      (5.0, 32.86009143933869)
      (5.0, 0.0)
    };
    \addplot+ [line width = 1,
    solid,mark = none,
    color = black,
    mark size = 2.0,
    mark options = {
      color = black,
      fill = kth-blue,
      line width = 1,
      rotate = 0,
      solid
    },fill = kth-blue, forget plot]coordinates {
      (5.0, 32.86009143933869)
      (6.0, 32.86009143933869)
      (6.0, 45.97268505796616)
      (5.0, 45.97268505796616)
      (5.0, 32.86009143933869)
    };
    \addplot+ [line width = 1,
    solid,mark = none,
    color = black,
    mark size = 2.0,
    mark options = {
      color = black,
      fill = kth-blue,
      line width = 1,
      rotate = 0,
      solid
    },fill = kth-blue, forget plot]coordinates {
      (6.0, 0.0)
      (7.0, 0.0)
      (7.0, 50.00373002427962)
      (6.0, 50.00373002427962)
      (6.0, 0.0)
    };
    \addplot+ [line width = 1,
    solid,mark = none,
    color = black,
    mark size = 2.0,
    mark options = {
      color = black,
      fill = kth-blue,
      line width = 1,
      rotate = 0,
      solid
    },fill = kth-blue, forget plot]coordinates {
      (7.0, 0.0)
      (8.0, 0.0)
      (8.0, 49.712140792857895)
      (7.0, 49.712140792857895)
      (7.0, 0.0)
    };
    \addplot+ [line width = 1,
    solid,mark = none,
    color = black,
    mark size = 2.0,
    mark options = {
      color = black,
      fill = kth-blue,
      line width = 1,
      rotate = 0,
      solid
    },fill = kth-blue, forget plot]coordinates {
      (7.0, 49.712140792857895)
      (8.0, 49.712140792857895)
      (8.0, 62.824734411485366)
      (7.0, 62.824734411485366)
      (7.0, 49.712140792857895)
    };
    \addplot+ [line width = 1,
    solid,mark = none,
    color = black,
    mark size = 2.0,
    mark options = {
      color = black,
      fill = kth-blue,
      line width = 1,
      rotate = 0,
      solid
    },fill = kth-blue, forget plot]coordinates {
      (8.0, 0.0)
      (9.0, 0.0)
      (9.0, 12.030588153077126)
      (8.0, 12.030588153077126)
      (8.0, 0.0)
    };
    \addplot+ [line width = 1,
    solid,mark = none,
    color = black,
    mark size = 2.0,
    mark options = {
      color = black,
      fill = kth-blue,
      line width = 1,
      rotate = 0,
      solid
    },fill = kth-blue, forget plot]coordinates {
      (8.0, 12.030588153077126)
      (9.0, 12.030588153077126)
      (9.0, 23.359214580154894)
      (8.0, 23.359214580154894)
      (8.0, 12.030588153077126)
    };
    \addplot+ [line width = 1,
    solid,mark = none,
    color = black,
    mark size = 2.0,
    mark options = {
      color = black,
      fill = kth-blue,
      line width = 1,
      rotate = 0,
      solid
    },fill = kth-blue, forget plot]coordinates {
      (8.0, 23.359214580154894)
      (9.0, 23.359214580154894)
      (9.0, 57.34509386138821)
      (8.0, 57.34509386138821)
      (8.0, 23.359214580154894)
    };
    \addplot+ [line width = 1,
    solid,mark = none,
    color = black,
    mark size = 2.0,
    mark options = {
      color = black,
      fill = kth-blue,
      line width = 1,
      rotate = 0,
      solid
    },fill = kth-blue, forget plot]coordinates {
      (8.0, 57.3450938613882)
      (9.0, 57.3450938613882)
      (9.0, 70.45768748001566)
      (8.0, 70.45768748001566)
      (8.0, 57.3450938613882)
    };
    \addplot+ [line width = 1,
    solid,mark = none,
    color = black,
    mark size = 2.0,
    mark options = {
      color = black,
      fill = kth-blue,
      line width = 1,
      rotate = 0,
      solid
    },fill = kth-blue, forget plot]coordinates {
      (9.0, 0.0)
      (10.0, 0.0)
      (10.0, 14.073677296654779)
      (9.0, 14.073677296654779)
      (9.0, 0.0)
    };
    \addplot+ [line width = 1,
    solid,mark = none,
    color = black,
    mark size = 2.0,
    mark options = {
      color = black,
      fill = kth-blue,
      line width = 1,
      rotate = 0,
      solid
    },fill = kth-blue, forget plot]coordinates {
      (9.0, 14.073677296654779)
      (10.0, 14.073677296654779)
      (10.0, 67.37438359568485)
      (9.0, 67.37438359568485)
      (9.0, 14.073677296654779)
    };
    \addplot+ [line width = 1,
    solid,mark = none,
    color = black,
    mark size = 2.0,
    mark options = {
      color = black,
      fill = kth-blue,
      line width = 1,
      rotate = 0,
      solid
    },fill = kth-blue, forget plot]coordinates {
      (10.0, 0.0)
      (11.0, 0.0)
      (11.0, 17.093752087298178)
      (10.0, 17.093752087298178)
      (10.0, 0.0)
    };
    \addplot+ [line width = 1,
    solid,mark = none,
    color = black,
    mark size = 2.0,
    mark options = {
      color = black,
      fill = kth-blue,
      line width = 1,
      rotate = 0,
      solid
    },fill = kth-blue, forget plot]coordinates {
      (10.0, 17.093752087298178)
      (11.0, 17.093752087298178)
      (11.0, 62.54236682849421)
      (10.0, 62.54236682849421)
      (10.0, 17.093752087298178)
    };
    \addplot+ [line width = 1,
    solid,mark = none,
    color = black,
    mark size = 2.0,
    mark options = {
      color = black,
      fill = kth-blue,
      line width = 1,
      rotate = 0,
      solid
    },fill = kth-blue, forget plot]coordinates {
      (11.0, 0.0)
      (12.0, 0.0)
      (12.0, 59.80312977614818)
      (11.0, 59.80312977614818)
      (11.0, 0.0)
    };
    \addplot+ [line width = 1,
    solid,mark = none,
    color = black,
    mark size = 2.0,
    mark options = {
      color = black,
      fill = kth-blue,
      line width = 1,
      rotate = 0,
      solid
    },fill = kth-blue, forget plot]coordinates {
      (12.0, 0.0)
      (13.0, 0.0)
      (13.0, 58.74093007310725)
      (12.0, 58.74093007310725)
      (12.0, 0.0)
    };
    \addplot+ [line width = 1,
    solid,mark = none,
    color = black,
    mark size = 2.0,
    mark options = {
      color = black,
      fill = kth-blue,
      line width = 1,
      rotate = 0,
      solid
    },fill = kth-blue, forget plot]coordinates {
      (13.0, 0.0)
      (14.0, 0.0)
      (14.0, 57.23391661133232)
      (13.0, 57.23391661133232)
      (13.0, 0.0)
    };
    \addplot+ [line width = 1,
    solid,mark = none,
    color = black,
    mark size = 2.0,
    mark options = {
      color = black,
      fill = kth-blue,
      line width = 1,
      rotate = 0,
      solid
    },fill = kth-blue, forget plot]coordinates {
      (14.0, 0.0)
      (15.0, 0.0)
      (15.0, 57.25436618645419)
      (14.0, 57.25436618645419)
      (14.0, 0.0)
    };
    \addplot+ [line width = 1,
    solid,mark = none,
    color = black,
    mark size = 2.0,
    mark options = {
      color = black,
      fill = kth-blue,
      line width = 1,
      rotate = 0,
      solid
    },fill = kth-blue, forget plot]coordinates {
      (15.0, 0.0)
      (16.0, 0.0)
      (16.0, 47.161460775019506)
      (15.0, 47.161460775019506)
      (15.0, 0.0)
    };
    \addplot+ [line width = 1,
    solid,mark = none,
    color = black,
    mark size = 2.0,
    mark options = {
      color = black,
      fill = kth-blue,
      line width = 1,
      rotate = 0,
      solid
    },fill = kth-blue, forget plot]coordinates {
      (15.0, 47.161460775019506)
      (16.0, 47.161460775019506)
      (16.0, 60.27405439364698)
      (15.0, 60.27405439364698)
      (15.0, 47.161460775019506)
    };
    \addplot+ [line width = 1,
    solid,mark = none,
    color = black,
    mark size = 2.0,
    mark options = {
      color = black,
      fill = kth-blue,
      line width = 1,
      rotate = 0,
      solid
    },fill = kth-blue, forget plot]coordinates {
      (16.0, 0.0)
      (17.0, 0.0)
      (17.0, 9.899841711437812)
      (16.0, 9.899841711437812)
      (16.0, 0.0)
    };
    \addplot+ [line width = 1,
    solid,mark = none,
    color = black,
    mark size = 2.0,
    mark options = {
      color = black,
      fill = kth-blue,
      line width = 1,
      rotate = 0,
      solid
    },fill = kth-blue, forget plot]coordinates {
      (16.0, 9.899841711437812)
      (17.0, 9.899841711437812)
      (17.0, 34.483305475234985)
      (16.0, 34.483305475234985)
      (16.0, 9.899841711437812)
    };
    \addplot+ [line width = 1,
    solid,mark = none,
    color = black,
    mark size = 2.0,
    mark options = {
      color = black,
      fill = kth-blue,
      line width = 1,
      rotate = 0,
      solid
    },fill = kth-blue, forget plot]coordinates {
      (16.0, 34.483305475234985)
      (17.0, 34.483305475234985)
      (17.0, 46.77503735713357)
      (16.0, 46.77503735713357)
      (16.0, 34.483305475234985)
    };
    \addplot+ [line width = 1,
    solid,mark = none,
    color = black,
    mark size = 2.0,
    mark options = {
      color = black,
      fill = kth-blue,
      line width = 1,
      rotate = 0,
      solid
    },fill = kth-blue, forget plot]coordinates {
      (16.0, 46.77503735713357)
      (17.0, 46.77503735713357)
      (17.0, 59.066769239032155)
      (16.0, 59.066769239032155)
      (16.0, 46.77503735713357)
    };
    \addplot+ [line width = 1,
    solid,mark = none,
    color = black,
    mark size = 2.0,
    mark options = {
      color = black,
      fill = kth-blue,
      line width = 1,
      rotate = 0,
      solid
    },fill = kth-blue, forget plot]coordinates {
      (16.0, 59.066769239032155)
      (17.0, 59.066769239032155)
      (17.0, 72.17936285765963)
      (16.0, 72.17936285765963)
      (16.0, 59.066769239032155)
    };
    \addplot+ [line width = 1,
    solid,mark = none,
    color = black,
    mark size = 2.0,
    mark options = {
      color = black,
      fill = kth-blue,
      line width = 1,
      rotate = 0,
      solid
    },fill = kth-blue, forget plot]coordinates {
      (17.0, 0.0)
      (18.0, 0.0)
      (18.0, 75.99514545892501)
      (17.0, 75.99514545892501)
      (17.0, 0.0)
    };
    \addplot+ [line width = 1,
    solid,mark = none,
    color = black,
    mark size = 2.0,
    mark options = {
      color = black,
      fill = kth-blue,
      line width = 1,
      rotate = 0,
      solid
    },fill = kth-blue, forget plot]coordinates {
      (18.0, 0.0)
      (19.0, 0.0)
      (19.0, 64.94566713642983)
      (18.0, 64.94566713642983)
      (18.0, 0.0)
    };
    \addplot+ [line width = 1,
    solid,mark = none,
    color = black,
    mark size = 2.0,
    mark options = {
      color = black,
      fill = kth-blue,
      line width = 1,
      rotate = 0,
      solid
    },fill = kth-blue, forget plot]coordinates {
      (19.0, 0.0)
      (20.0, 0.0)
      (20.0, 41.33787737220982)
      (19.0, 41.33787737220982)
      (19.0, 0.0)
    };
    \addplot+ [line width = 1,
    solid,mark = none,
    color = black,
    mark size = 2.0,
    mark options = {
      color = black,
      fill = kth-blue,
      line width = 1,
      rotate = 0,
      solid
    },fill = kth-blue, forget plot]coordinates {
      (19.0, 41.33787737220982)
      (20.0, 41.33787737220982)
      (20.0, 54.45047099083729)
      (19.0, 54.45047099083729)
      (19.0, 41.33787737220982)
    };
    \addplot+ [line width = 1,
    solid,mark = none,
    color = black,
    mark size = 2.0,
    mark options = {
      color = black,
      fill = kth-blue,
      line width = 1,
      rotate = 0,
      solid
    },fill = kth-blue, forget plot]coordinates {
      (20.0, 0.0)
      (21.0, 0.0)
      (21.0, 36.61898147571877)
      (20.0, 36.61898147571877)
      (20.0, 0.0)
    };
    \addplot+ [line width = 1,
    solid,mark = none,
    color = black,
    mark size = 2.0,
    mark options = {
      color = black,
      fill = kth-blue,
      line width = 1,
      rotate = 0,
      solid
    },fill = kth-blue, forget plot]coordinates {
      (20.0, 36.61898147571877)
      (21.0, 36.61898147571877)
      (21.0, 49.73157509434624)
      (20.0, 49.73157509434624)
      (20.0, 36.61898147571877)
    };
    \addplot+ [line width = 1,
    solid,mark = none,
    color = black,
    mark size = 2.0,
    mark options = {
      color = black,
      fill = kth-blue,
      line width = 1,
      rotate = 0,
      solid
    },fill = kth-blue, forget plot]coordinates {
      (21.0, 0.0)
      (22.0, 0.0)
      (22.0, 47.70970733207092)
      (21.0, 47.70970733207092)
      (21.0, 0.0)
    };
    \addplot+ [line width = 1,
    solid,mark = none,
    color = black,
    mark size = 2.0,
    mark options = {
      color = black,
      fill = kth-blue,
      line width = 1,
      rotate = 0,
      solid
    },fill = kth-blue, forget plot]coordinates {
      (22.0, 0.0)
      (23.0, 0.0)
      (23.0, 27.690850854873656)
      (22.0, 27.690850854873656)
      (22.0, 0.0)
    };
    \addplot+ [line width = 1,
    solid,mark = none,
    color = black,
    mark size = 2.0,
    mark options = {
      color = black,
      fill = kth-blue,
      line width = 1,
      rotate = 0,
      solid
    },fill = kth-blue, forget plot]coordinates {
      (22.0, 27.690850854873656)
      (23.0, 27.690850854873656)
      (23.0, 30.49850052089886)
      (22.0, 30.49850052089886)
      (22.0, 27.690850854873656)
    };
    \addplot+ [line width = 1,
    solid,mark = none,
    color = black,
    mark size = 2.0,
    mark options = {
      color = black,
      fill = kth-blue,
      line width = 1,
      rotate = 0,
      solid
    },fill = kth-blue, forget plot]coordinates {
      (22.0, 30.49850052089886)
      (23.0, 30.49850052089886)
      (23.0, 46.41874380555153)
      (22.0, 46.41874380555153)
      (22.0, 30.49850052089886)
    };
    \addplot+ [line width = 1,
    solid,mark = none,
    color = black,
    mark size = 2.0,
    mark options = {
      color = black,
      fill = kth-blue,
      line width = 1,
      rotate = 0,
      solid
    },fill = kth-blue, forget plot]coordinates {
      (23.0, 0.0)
      (24.0, 0.0)
      (24.0, 29.440919082076665)
      (23.0, 29.440919082076665)
      (23.0, 0.0)
    };
    \addplot+ [line width = 1,
    solid,mark = none,
    color = black,
    mark size = 2.0,
    mark options = {
      color = black,
      fill = kth-blue,
      line width = 1,
      rotate = 0,
      solid
    },fill = kth-blue, forget plot]coordinates {
      (23.0, 29.440919082076665)
      (24.0, 29.440919082076665)
      (24.0, 32.15981687909626)
      (23.0, 32.15981687909626)
      (23.0, 29.440919082076665)
    };
    \addplot+ [line width = 1,
    solid,mark = none,
    color = black,
    mark size = 2.0,
    mark options = {
      color = black,
      fill = kth-blue,
      line width = 1,
      rotate = 0,
      solid
    },fill = kth-blue, forget plot]coordinates {
      (23.0, 32.15981687909626)
      (24.0, 32.15981687909626)
      (24.0, 45.272410497723726)
      (23.0, 45.272410497723726)
      (23.0, 32.15981687909626)
    };
    \addplot+ [line width = 1,
    solid,mark = none,
    color = black,
    mark size = 2.0,
    mark options = {
      fill = kth-blue,
      line width = 1,
      rotate = 0,
      solid
    },fill = kth-blue, forget plot]coordinates {
      (17.0, 90)
      (24.0, 90)
      (24.0, 100)
      (17.0, 100)
      (17.0, 90)
    };
    \node at (axis cs:0.5, -6.556296809313734) [,
    color=white,
    rotate=0.0,
    font={\fontsize{6 pt}{15.600000000000001 pt}\bfseries\sffamily}
    ] {685.7};
    \node at (axis cs:1.5, 22.245548281076893) [,
    color=white,
    rotate=0.0,
    font={\fontsize{6 pt}{15.600000000000001 pt}\bfseries\sffamily}
    ] {644.5};
    \node at (axis cs:2.5, 10.152368732118179) [,
    color=white,
    rotate=0.0,
    font={\fontsize{6 pt}{15.600000000000001 pt}\bfseries\sffamily}
    ] {615.1};
    \node at (axis cs:2.5, 32.145815298197704) [,
    color=white,
    rotate=0.0,
    font={\fontsize{6 pt}{15.600000000000001 pt}\bfseries\sffamily}
    ] {615.8};
    \node at (axis cs:3.5, 9.908791176735402) [,
    color=white,
    rotate=0.0,
    font={\fontsize{6 pt}{15.600000000000001 pt}\bfseries\sffamily}
    ] {645.9};
    \node at (axis cs:3.5, 25.48140403175354) [,
    color=white,
    rotate=0.0,
    font={\fontsize{6 pt}{15.600000000000001 pt}\bfseries\sffamily}
    ] {658.7};
    \node at (axis cs:3.5, 37.70152251935001) [,
    color=white,
    rotate=0.0,
    font={\fontsize{6 pt}{15.600000000000001 pt}\bfseries\sffamily}
    ] {2000};
    \node at (axis cs:4.5, 10.101605695750584) [,
    color=white,
    rotate=0.0,
    font={\fontsize{6 pt}{15.600000000000001 pt}\bfseries\sffamily}
    ] {689.9};
    \node at (axis cs:4.5, 23.089947512652742) [,
    color=white,
    rotate=0.0,
    font={\fontsize{6 pt}{15.600000000000001 pt}\bfseries\sffamily}
    ] {690.5};
    \node at (axis cs:4.5, 35.419716564269635) [,
    color=white,
    rotate=0.0,
    font={\fontsize{6 pt}{15.600000000000001 pt}\bfseries\sffamily}
    ] {691.5};
    \node at (axis cs:5.5, 16.430045719669344) [,
    color=white,
    rotate=0.0,
    font={\fontsize{6 pt}{15.600000000000001 pt}\bfseries\sffamily}
    ] {698.1};
    \node at (axis cs:5.5, 39.416388248652424) [,
    color=white,
    rotate=0.0,
    font={\fontsize{6 pt}{15.600000000000001 pt}\bfseries\sffamily}
    ] {2000};
    \node at (axis cs:6.5, 25.00186501213981) [,
    color=white,
    rotate=0.0,
    font={\fontsize{6 pt}{15.600000000000001 pt}\bfseries\sffamily}
    ] {733.0};
    \node at (axis cs:7.5, 24.856070396428947) [,
    color=white,
    rotate=0.0,
    font={\fontsize{6 pt}{15.600000000000001 pt}\bfseries\sffamily}
    ] {734.1};
    \node at (axis cs:7.5, 56.26843760217163) [,
    color=white,
    rotate=0.0,
    font={\fontsize{6 pt}{15.600000000000001 pt}\bfseries\sffamily}
    ] {1300};
    \node at (axis cs:8.5, 6.015294076538563) [,
    color=white,
    rotate=0.0,
    font={\fontsize{6 pt}{15.600000000000001 pt}\bfseries\sffamily}
    ] {551.9};
    \node at (axis cs:8.5, 17.69490136661601) [,
    color=white,
    rotate=0.0,
    font={\fontsize{6 pt}{15.600000000000001 pt}\bfseries\sffamily}
    ] {810.3};
    \node at (axis cs:8.5, 40.352154220771546) [,
    color=white,
    rotate=0.0,
    font={\fontsize{6 pt}{15.600000000000001 pt}\bfseries\sffamily}
    ] {907.1};
    \node at (axis cs:8.5, 63.90139067070193) [,
    color=white,
    rotate=0.0,
    font={\fontsize{6 pt}{15.600000000000001 pt}\bfseries\sffamily}
    ] {1200};
    \node at (axis cs:9.5, 7.0368386483273895) [,
    color=white,
    rotate=0.0,
    font={\fontsize{6 pt}{15.600000000000001 pt}\bfseries\sffamily}
    ] {606.1};
    \node at (axis cs:9.5, 40.72403044616981) [,
    color=white,
    rotate=0.0,
    font={\fontsize{6 pt}{15.600000000000001 pt}\bfseries\sffamily}
    ] {869.7};
    \node at (axis cs:10.5, 8.546876043649089) [,
    color=white,
    rotate=0.0,
    font={\fontsize{6 pt}{15.600000000000001 pt}\bfseries\sffamily}
    ] {303.5};
    \node at (axis cs:10.5, 39.8180594578962) [,
    color=white,
    rotate=0.0,
    font={\fontsize{6 pt}{15.600000000000001 pt}\bfseries\sffamily}
    ] {323.9};
    \node at (axis cs:11.5, 29.90156488807409) [,
    color=white,
    rotate=0.0,
    font={\fontsize{6 pt}{15.600000000000001 pt}\bfseries\sffamily}
    ] {408.3};
    \node at (axis cs:12.5, 29.370465036553625) [,
    color=white,
    rotate=0.0,
    font={\fontsize{6 pt}{15.600000000000001 pt}\bfseries\sffamily}
    ] {349.1};
    \node at (axis cs:13.5, 28.61695830566616) [,
    color=white,
    rotate=0.0,
    font={\fontsize{6 pt}{15.600000000000001 pt}\bfseries\sffamily}
    ] {363.9};
    \node at (axis cs:14.5, 28.627183093227096) [,
    color=white,
    rotate=0.0,
    font={\fontsize{6 pt}{15.600000000000001 pt}\bfseries\sffamily}
    ] {364.1};
    \node at (axis cs:15.5, 23.580730387509753) [,
    color=white,
    rotate=0.0,
    font={\fontsize{6 pt}{15.600000000000001 pt}\bfseries\sffamily}
    ] {45.57};
    \node at (axis cs:15.5, 53.71775758433324) [,
    color=white,
    rotate=0.0,
    font={\fontsize{6 pt}{15.600000000000001 pt}\bfseries\sffamily}
    ] {1100};
    \node at (axis cs:16.5, 4.949920855718906) [,
    color=white,
    rotate=0.0,
    font={\fontsize{6 pt}{15.600000000000001 pt}\bfseries\sffamily}
    ] {0.0};
    \node at (axis cs:16.5, 22.1915735933364) [,
    color=white,
    rotate=0.0,
    font={\fontsize{6 pt}{15.600000000000001 pt}\bfseries\sffamily}
    ] {74.55};
    \node at (axis cs:16.5, 40.62917141618428) [,
    color=white,
    rotate=0.0,
    font={\fontsize{6 pt}{15.600000000000001 pt}\bfseries\sffamily}
    ] {97.42};
    \node at (axis cs:16.5, 52.92090329808286) [,
    color=white,
    rotate=0.0,
    font={\fontsize{6 pt}{15.600000000000001 pt}\bfseries\sffamily}
    ] {98.62};
    \node at (axis cs:16.5, 65.62306604834589) [,
    color=white,
    rotate=0.0,
    font={\fontsize{6 pt}{15.600000000000001 pt}\bfseries\sffamily}
    ] {367.5};
    \node at (axis cs:17.5, 37.997572729462505) [,
    color=white,
    rotate=0.0,
    font={\fontsize{6 pt}{15.600000000000001 pt}\bfseries\sffamily}
    ] {97.15};
    \node at (axis cs:18.5, 32.47283356821492) [,
    color=white,
    rotate=0.0,
    font={\fontsize{6 pt}{15.600000000000001 pt}\bfseries\sffamily}
    ] {109.9};
    \node at (axis cs:19.5, 20.66893868610491) [,
    color=white,
    rotate=0.0,
    font={\fontsize{6 pt}{15.600000000000001 pt}\bfseries\sffamily}
    ] {486.0};
    \node at (axis cs:19.5, 47.89417418152355) [,
    color=white,
    rotate=0.0,
    font={\fontsize{6 pt}{15.600000000000001 pt}\bfseries\sffamily}
    ] {1500};
    \node at (axis cs:20.5, 18.309490737859385) [,
    color=white,
    rotate=0.0,
    font={\fontsize{6 pt}{15.600000000000001 pt}\bfseries\sffamily}
    ] {482.0};
    \node at (axis cs:20.5, 43.175278285032505) [,
    color=white,
    rotate=0.0,
    font={\fontsize{6 pt}{15.600000000000001 pt}\bfseries\sffamily}
    ] {1500};
    \node at (axis cs:21.5, 23.85485366603546) [,
    color=white,
    rotate=0.0,
    font={\fontsize{6 pt}{15.600000000000001 pt}\bfseries\sffamily}
    ] {510.7};
    \node at (axis cs:22.5, 13.845425427436828) [,
    color=white,
    rotate=0.0,
    font={\fontsize{6 pt}{15.600000000000001 pt}\bfseries\sffamily}
    ] {1000};
    \node at (axis cs:22.5, 29.09467568788626) [,
    color=white,
    rotate=0.0,
    font={\fontsize{6 pt}{15.600000000000001 pt}\bfseries\sffamily}
    ] {1000};
    \node at (axis cs:22.5, 38.45862216322519) [,
    color=white,
    rotate=0.0,
    font={\fontsize{6 pt}{15.600000000000001 pt}\bfseries\sffamily}
    ] {1000};
    \node at (axis cs:23.5, 14.720459541038332) [,
    color=white,
    rotate=0.0,
    font={\fontsize{6 pt}{15.600000000000001 pt}\bfseries\sffamily}
    ] {1000};
    \node at (axis cs:23.5, 30.800367980586465) [,
    color=white,
    rotate=0.0,
    font={\fontsize{6 pt}{15.600000000000001 pt}\bfseries\sffamily}
    ] {1000};
    \node at (axis cs:23.5, 38.71611368840999) [,
    color=white,
    rotate=0.0,
    font={\fontsize{6 pt}{15.600000000000001 pt}\bfseries\sffamily}
    ] {1100};
    \node at (axis cs:21.0, 85.0) [,
    color=white,
    rotate=0.0,
    font={\fontsize{8 pt}{15.600000000000001 pt}\bfseries\sffamily}
    ] {Independent Volumes [MWh]};
    \node at (axis cs:21.0, 95.0) [,
    color=white,
    rotate=0.0,
    font={\fontsize{8 pt}{15.600000000000001 pt}\bfseries\sffamily}
    ] {Dependent Volumes [MWh]};
  \end{axis}
  \begin{axis}[
    width = {\textwidth},
    height = {0.3\textheight},
    xlabel = {Hour},
    ylabel = {Price [Eur/MWh]},
    title = {Block orders},
    unbounded coords=jump,
    scaled x ticks = false,
    xmajorgrids = true,
    xmin = 0,
    xmax = 24,
    xtick align = inside,
    xtick = {0.0,1.0,2.0,3.0,4.0,5.0,6.0,7.0,8.0,9.0,10.0,11.0,12.0,13.0,14.0,15.0,16.0,17.0,18.0,19.0,20.0,21.0,22.0,23.0,24.0},
    xticklabels = {$0$,$1$,$2$,$3$,$4$,$5$,$6$,$7$,$8$,$9$,$10$,$11$,$12$,$13$,$14$,$15$,$16$,$17$,$18$,$19$,$20$,$21$,$22$,$23$,$24$},
    xticklabel style = {font = {\fontsize{10 pt}{10.4 pt}\selectfont}, rotate = 0.0},
    x grid style = {color = kth-lightgray,
      line width = 0.25,
      solid},
    axis x line* = left,
    x axis line style = {line width = 1,
      solid},
    scaled y ticks = false,
    ylabel style = {font = {\fontsize{10 pt}{14.3 pt}\selectfont}, rotate = 0.0},
    ymajorgrids = true,
    ytick align = inside,
    yticklabel style = {font = {\fontsize{10 pt}{10.4 pt}\selectfont}, rotate = 0.0},
    y grid style = {color = kth-lightgray,
      line width = 0.25,
      solid},axis y line* = left,y axis line style = {line width = 1,
      solid},
    xshift = 0mm,
    title style = {font = {\fontsize{14 pt}{18.2 pt}\selectfont}, rotate = 0.0},
    legend pos = {north west},
    legend style = {line width = 2, solid,font = {\fontsize{8 pt}{10.4 pt}\selectfont}},colorbar style={title=}]
    \addplot+ [line width = 1,
    solid,mark = none,
    color = black,
    mark size = 2.0,
    mark options = {
      color = black,
      fill = kth-blue,
      line width = 1,
      rotate = 0,
      solid
    },fill = kth-blue, forget plot]coordinates {
      (18.0, 13.0)
      (24.0, 13.0)
      (24.0, 14.0)
      (18.0, 14.0)
      (18.0, 13.0)
    };
    \addplot+ [line width = 1,
    solid,mark = none,
    color = black,
    mark size = 2.0,
    mark options = {
      color = black,
      fill = kth-blue,
      line width = 1,
      rotate = 0,
      solid
    },fill = kth-blue, forget plot]coordinates {
      (15.0, 0.5)
      (19.0, 0.5)
      (19.0, 1.5)
      (15.0, 1.5)
      (15.0, 0.5)
    };
    \addplot+ [line width = 1,
    solid,mark = none,
    color = black,
    mark size = 2.0,
    mark options = {
      color = black,
      fill = kth-blue,
      line width = 1,
      rotate = 0,
      solid
    },fill = kth-blue, forget plot]coordinates {
      (10.0, 1.5)
      (22.0, 1.5)
      (22.0, 2.5)
      (10.0, 2.5)
      (10.0, 1.5)
    };
    \addplot+ [line width = 1,
    solid,mark = none,
    color = black,
    mark size = 2.0,
    mark options = {
      color = black,
      fill = kth-blue,
      line width = 1,
      rotate = 0,
      solid
    },fill = kth-blue, forget plot]coordinates {
      (15.0, 2.5)
      (18.0, 2.5)
      (18.0, 3.5)
      (15.0, 3.5)
      (15.0, 2.5)
    };
    \addplot+ [line width = 1,
    solid,mark = none,
    color = black,
    mark size = 2.0,
    mark options = {
      color = black,
      fill = kth-blue,
      line width = 1,
      rotate = 0,
      solid
    },fill = kth-blue, forget plot]coordinates {
      (7.0, 3.5)
      (10.0, 3.5)
      (10.0, 4.5)
      (7.0, 4.5)
      (7.0, 3.5)
    };
    \addplot+ [line width = 1,
    solid,mark = none,
    color = black,
    mark size = 2.0,
    mark options = {
      color = black,
      fill = kth-blue,
      line width = 1,
      rotate = 0,
      solid
    },fill = kth-blue, forget plot]coordinates {
      (15.0, 4.5)
      (19.0, 4.5)
      (19.0, 5.5)
      (15.0, 5.5)
      (15.0, 4.5)
    };
    \addplot+ [line width = 1,
    solid,mark = none,
    color = black,
    mark size = 2.0,
    mark options = {
      color = black,
      fill = kth-blue,
      line width = 1,
      rotate = 0,
      solid
    },fill = kth-blue, forget plot]coordinates {
      (21.0, 5.5)
      (24.0, 5.5)
      (24.0, 6.5)
      (21.0, 6.5)
      (21.0, 5.5)
    };
    \addplot+ [line width = 1,
    solid,mark = none,
    color = black,
    mark size = 2.0,
    mark options = {
      color = black,
      fill = kth-blue,
      line width = 1,
      rotate = 0,
      solid
    },fill = kth-blue, forget plot]coordinates {
      (17.0, 6.5)
      (24.0, 6.5)
      (24.0, 7.5)
      (17.0, 7.5)
      (17.0, 6.5)
    };
    \addplot+ [line width = 1,
    solid,mark = none,
    color = black,
    mark size = 2.0,
    mark options = {
      color = black,
      fill = kth-blue,
      line width = 1,
      rotate = 0,
      solid
    },fill = kth-blue, forget plot]coordinates {
      (5.0, 7.5)
      (9.0, 7.5)
      (9.0, 8.5)
      (5.0, 8.5)
      (5.0, 7.5)
    };
    \addplot+ [line width = 1,
    solid,mark = none,
    color = black,
    mark size = 2.0,
    mark options = {
      color = black,
      fill = kth-blue,
      line width = 1,
      rotate = 0,
      solid
    },fill = kth-blue, forget plot]coordinates {
      (17.0, 8.5)
      (21.0, 8.5)
      (21.0, 9.5)
      (17.0, 9.5)
      (17.0, 8.5)
    };
    \addplot+ [line width = 1,
    solid,mark = none,
    color = black,
    mark size = 2.0,
    mark options = {
      color = black,
      fill = kth-blue,
      line width = 1,
      rotate = 0,
      solid
    },fill = kth-blue, forget plot]coordinates {
      (16.0, 9.5)
      (20.0, 9.5)
      (20.0, 10.5)
      (16.0, 10.5)
      (16.0, 9.5)
    };
    \addplot+ [line width = 1,
    solid,mark = none,
    color = black,
    mark size = 2.0,
    mark options = {
      color = black,
      fill = kth-blue,
      line width = 1,
      rotate = 0,
      solid
    },fill = kth-blue, forget plot]coordinates {
      (15.0, 10.5)
      (18.0, 10.5)
      (18.0, 11.5)
      (15.0, 11.5)
      (15.0, 10.5)
    };
    \node at (axis cs:15.5, 1) [,
    color=white,
    rotate=0.0,
    font={\fontsize{6 pt}{18.2 pt}\bfseries\sffamily}
    ] {13.60};
    \node at (axis cs:18.5, 1) [,
    color=white,
    rotate=0.0,
    font={\fontsize{6 pt}{18.2 pt}\bfseries\sffamily}
    ] {399.3};
    \node at (axis cs:10.5, 2) [,
    color=white,
    rotate=0.0,
    font={\fontsize{6 pt}{18.2 pt}\bfseries\sffamily}
    ] {17.72};
    \node at (axis cs:21.5, 2) [,
    color=white,
    rotate=0.0,
    font={\fontsize{6 pt}{18.2 pt}\bfseries\sffamily}
    ] {500.0};
    \node at (axis cs:15.5, 3) [,
    color=white,
    rotate=0.0,
    font={\fontsize{6 pt}{18.2 pt}\bfseries\sffamily}
    ] {23.37};
    \node at (axis cs:17.5, 3) [,
    color=white,
    rotate=0.0,
    font={\fontsize{6 pt}{18.2 pt}\bfseries\sffamily}
    ] {0.18};
    \node at (axis cs:7.5, 4) [,
    color=white,
    rotate=0.0,
    font={\fontsize{6 pt}{18.2 pt}\bfseries\sffamily}
    ] {23.74};
    \node at (axis cs:9.5, 4) [,
    color=white,
    rotate=0.0,
    font={\fontsize{6 pt}{18.2 pt}\bfseries\sffamily}
    ] {11.70};
    \node at (axis cs:15.5, 5) [,
    color=white,
    rotate=0.0,
    font={\fontsize{6 pt}{18.2 pt}\bfseries\sffamily}
    ] {24.01};
    \node at (axis cs:18.5, 5) [,
    color=white,
    rotate=0.0,
    font={\fontsize{6 pt}{18.2 pt}\bfseries\sffamily}
    ] {0.26};
    \node at (axis cs:21.5, 6) [,
    color=white,
    rotate=0.0,
    font={\fontsize{6 pt}{18.2 pt}\bfseries\sffamily}
    ] {30.52};
    \node at (axis cs:23.5, 6) [,
    color=white,
    rotate=0.0,
    font={\fontsize{6 pt}{18.2 pt}\bfseries\sffamily}
    ] {0.02};
    \node at (axis cs:17.5, 7) [,
    color=white,
    rotate=0.0,
    font={\fontsize{6 pt}{18.2 pt}\bfseries\sffamily}
    ] {41.82};
    \node at (axis cs:23.5, 7) [,
    color=white,
    rotate=0.0,
    font={\fontsize{6 pt}{18.2 pt}\bfseries\sffamily}
    ] {0.87};
    \node at (axis cs:5.5, 8) [,
    color=white,
    rotate=0.0,
    font={\fontsize{6 pt}{18.2 pt}\bfseries\sffamily}
    ] {44.20};
    \node at (axis cs:8.5, 8) [,
    color=white,
    rotate=0.0,
    font={\fontsize{6 pt}{18.2 pt}\bfseries\sffamily}
    ] {0.05};
    \node at (axis cs:17.5, 9) [,
    color=white,
    rotate=0.0,
    font={\fontsize{6 pt}{18.2 pt}\bfseries\sffamily}
    ] {48.17};
    \node at (axis cs:20.5, 9) [,
    color=white,
    rotate=0.0,
    font={\fontsize{6 pt}{18.2 pt}\bfseries\sffamily}
    ] {0.31};
    \node at (axis cs:16.5, 10) [,
    color=white,
    rotate=0.0,
    font={\fontsize{6 pt}{18.2 pt}\bfseries\sffamily}
    ] {53.78};
    \node at (axis cs:19.5, 10) [,
    color=white,
    rotate=0.0,
    font={\fontsize{6 pt}{18.2 pt}\bfseries\sffamily}
    ] {0.06};
    \node at (axis cs:15.5, 11) [,
    color=white,
    rotate=0.0,
    font={\fontsize{6 pt}{18.2 pt}\bfseries\sffamily}
    ] {56.37};
    \node at (axis cs:17.5, 11) [,
    color=white,
    rotate=0.0,
    font={\fontsize{6 pt}{18.2 pt}\bfseries\sffamily}
    ] {0.17};
    \node at (axis cs:19.4, 13.5) [,
    color=white,
    rotate=0.0,
    font={\fontsize{6 pt}{18.2 pt}\bfseries\sffamily}
    ] {Price [EUR/MWh]};
    \node at (axis cs:22.6, 13.5) [,
    color=white,
    rotate=0.0,
    font={\fontsize{6 pt}{18.2 pt}\bfseries\sffamily}
    ] {Volume [MWh/h]};
  \end{axis}
\end{tikzpicture}

%% file: day-ahead_evp_strategy.tex
\begin{tikzpicture}[]
  \begin{axis}[
    width = {\textwidth},
    height = {0.3\textheight},
    xlabel = {},
    ylabel = {Price [Eur/MWh]},
    title = {Single orders},
    unbounded coords=jump,
    scaled x ticks = false,
    xmajorgrids = true,
    xmin = 0,
    xmax = 24,
    ymin = {-13.112593618627468},
    ymax = {102.22033269617995},
    xtick align = inside,
    xtick = {0.0,1.0,2.0,3.0,4.0,5.0,6.0,7.0,8.0,9.0,10.0,11.0,12.0,13.0,14.0,15.0,16.0,17.0,18.0,19.0,20.0,21.0,22.0,23.0,24.0},
    xticklabels = {$0$,$1$,$2$,$3$,$4$,$5$,$6$,$7$,$8$,$9$,$10$,$11$,$12$,$13$,$14$,$15$,$16$,$17$,$18$,$19$,$20$,$21$,$22$,$23$,$24$},
    xticklabel style = {font = {\fontsize{10 pt}{10.4 pt}\selectfont}, rotate = 0.0},
    x grid style = {color = kth-lightgray,
      line width = 0.25,
      solid},
    axis x line* = left,
    x axis line style = {line width = 1,
      solid},
    scaled y ticks = false,
    ylabel style = {yshift = -5pt, font = {\fontsize{10 pt}{14.3 pt}\selectfont}, rotate = 0.0},
    ymajorgrids = true,
    ytick align = inside,
    yticklabel style = {font = {\fontsize{10 pt}{10.4 pt}\selectfont}, rotate = 0.0},
    y grid style = {color = kth-lightgray,
      line width = 0.25,
      solid},axis y line* = left,y axis line style = {line width = 1,
      solid},
    xshift = 0mm,
    yshift = {0.3\textheight},
    title style = {font = {\fontsize{14 pt}{18.2 pt}\selectfont}, rotate = 0.0},
    legend pos = {north west},
    legend style = {line width = 2, solid,font = {\fontsize{8 pt}{10.4 pt}\selectfont}},colorbar style={title=}]
    \addplot+ [line width = 1,
    solid,mark = none,
    color = black,
    mark size = 2.0,
    mark options = {
      fill = kth-darkgray,
      line width = 1,
      rotate = 0,
      solid
    },fill = kth-darkgray, forget plot]coordinates {
      (0.0, -13.112593618627468)
      (1.0, -13.112593618627468)
      (1.0, 0.0)
      (0.0, 0.0)
      (0.0, -13.112593618627468)
    };
    \addplot+ [line width = 1,
    solid,mark = none,
    color = black,
    mark size = 2.0,
    mark options = {
      fill = kth-darkgray,
      line width = 1,
      rotate = 0,
      solid
    },fill = kth-darkgray, forget plot]coordinates {
      (1.0, -13.112593618627468)
      (2.0, -13.112593618627468)
      (2.0, 0.0)
      (1.0, 0.0)
      (1.0, -13.112593618627468)
    };
    \addplot+ [line width = 1,
    solid,mark = none,
    color = black,
    mark size = 2.0,
    mark options = {
      fill = kth-darkgray,
      line width = 1,
      rotate = 0,
      solid
    },fill = kth-darkgray, forget plot]coordinates {
      (2.0, -13.112593618627468)
      (3.0, -13.112593618627468)
      (3.0, 0.0)
      (2.0, 0.0)
      (2.0, -13.112593618627468)
    };
    \addplot+ [line width = 1,
    solid,mark = none,
    color = black,
    mark size = 2.0,
    mark options = {
      fill = kth-darkgray,
      line width = 1,
      rotate = 0,
      solid
    },fill = kth-darkgray, forget plot]coordinates {
      (4.0, -13.112593618627468)
      (5.0, -13.112593618627468)
      (5.0, 0.0)
      (4.0, 0.0)
      (4.0, -13.112593618627468)
    };
    \addplot+ [line width = 1,
    solid,mark = none,
    color = black,
    mark size = 2.0,
    mark options = {
      fill = kth-darkgray,
      line width = 1,
      rotate = 0,
      solid
    },fill = kth-darkgray, forget plot]coordinates {
      (5.0, -13.112593618627468)
      (6.0, -13.112593618627468)
      (6.0, 0.0)
      (5.0, 0.0)
      (5.0, -13.112593618627468)
    };
    \addplot+ [line width = 1,
    solid,mark = none,
    color = black,
    mark size = 2.0,
    mark options = {
      fill = kth-darkgray,
      line width = 1,
      rotate = 0,
      solid
    },fill = kth-darkgray, forget plot]coordinates {
      (6.0, -13.112593618627468)
      (7.0, -13.112593618627468)
      (7.0, 0.0)
      (6.0, 0.0)
      (6.0, -13.112593618627468)
    };
    \addplot+ [line width = 1,
    solid,mark = none,
    color = black,
    mark size = 2.0,
    mark options = {
      fill = kth-darkgray,
      line width = 1,
      rotate = 0,
      solid
    },fill = kth-darkgray, forget plot]coordinates {
      (7.0, -13.112593618627468)
      (8.0, -13.112593618627468)
      (8.0, 0.0)
      (7.0, 0.0)
      (7.0, -13.112593618627468)
    };
    \addplot+ [line width = 1,
    solid,mark = none,
    color = black,
    mark size = 2.0,
    mark options = {
      fill = kth-darkgray,
      line width = 1,
      rotate = 0,
      solid
    },fill = kth-darkgray, forget plot]coordinates {
      (8.0, -13.112593618627468)
      (9.0, -13.112593618627468)
      (9.0, 0.0)
      (8.0, 0.0)
      (8.0, -13.112593618627468)
    };
    \addplot+ [line width = 1,
    solid,mark = none,
    color = black,
    mark size = 2.0,
    mark options = {
      fill = kth-darkgray,
      line width = 1,
      rotate = 0,
      solid
    },fill = kth-darkgray, forget plot]coordinates {
      (9.0, -13.112593618627468)
      (10.0, -13.112593618627468)
      (10.0, 0.0)
      (9.0, 0.0)
      (9.0, -13.112593618627468)
    };
    \addplot+ [line width = 1,
    solid,mark = none,
    color = black,
    mark size = 2.0,
    mark options = {
      fill = kth-darkgray,
      line width = 1,
      rotate = 0,
      solid
    },fill = kth-darkgray, forget plot]coordinates {
      (10.0, -13.112593618627468)
      (11.0, -13.112593618627468)
      (11.0, 0.0)
      (10.0, 0.0)
      (10.0, -13.112593618627468)
    };
    \addplot+ [line width = 1,
    solid,mark = none,
    color = black,
    mark size = 2.0,
    mark options = {
      fill = kth-darkgray,
      line width = 1,
      rotate = 0,
      solid
    },fill = kth-darkgray, forget plot]coordinates {
      (11.0, -13.112593618627468)
      (12.0, -13.112593618627468)
      (12.0, 0.0)
      (11.0, 0.0)
      (11.0, -13.112593618627468)
    };
    \addplot+ [line width = 1,
    solid,mark = none,
    color = black,
    mark size = 2.0,
    mark options = {
      fill = kth-darkgray,
      line width = 1,
      rotate = 0,
      solid
    },fill = kth-darkgray, forget plot]coordinates {
      (12.0, -13.112593618627468)
      (13.0, -13.112593618627468)
      (13.0, 0.0)
      (12.0, 0.0)
      (12.0, -13.112593618627468)
    };
    \addplot+ [line width = 1,
    solid,mark = none,
    color = black,
    mark size = 2.0,
    mark options = {
      fill = kth-darkgray,
      line width = 1,
      rotate = 0,
      solid
    },fill = kth-darkgray, forget plot]coordinates {
      (13.0, -13.112593618627468)
      (14.0, -13.112593618627468)
      (14.0, 0.0)
      (13.0, 0.0)
      (13.0, -13.112593618627468)
    };
    \addplot+ [line width = 1,
    solid,mark = none,
    color = black,
    mark size = 2.0,
    mark options = {
      fill = kth-darkgray,
      line width = 1,
      rotate = 0,
      solid
    },fill = kth-darkgray, forget plot]coordinates {
      (14.0, -13.112593618627468)
      (15.0, -13.112593618627468)
      (15.0, 0.0)
      (14.0, 0.0)
      (14.0, -13.112593618627468)
    };
    \addplot+ [line width = 1,
    solid,mark = none,
    color = black,
    mark size = 2.0,
    mark options = {
      fill = kth-darkgray,
      line width = 1,
      rotate = 0,
      solid
    },fill = kth-darkgray, forget plot]coordinates {
      (15.0, -13.112593618627468)
      (16.0, -13.112593618627468)
      (16.0, 0.0)
      (15.0, 0.0)
      (15.0, -13.112593618627468)
    };
    \addplot+ [line width = 1,
    solid,mark = none,
    color = black,
    mark size = 2.0,
    mark options = {
      fill = kth-darkgray,
      line width = 1,
      rotate = 0,
      solid
    },fill = kth-darkgray, forget plot]coordinates {
      (16.0, -13.112593618627468)
      (17.0, -13.112593618627468)
      (17.0, 0.0)
      (16.0, 0.0)
      (16.0, -13.112593618627468)
    };
    \addplot+ [line width = 1,
    solid,mark = none,
    color = black,
    mark size = 2.0,
    mark options = {
      fill = kth-darkgray,
      line width = 1,
      rotate = 0,
      solid
    },fill = kth-darkgray, forget plot]coordinates {
      (18.0, -13.112593618627468)
      (19.0, -13.112593618627468)
      (19.0, 0.0)
      (18.0, 0.0)
      (18.0, -13.112593618627468)
    };
    \addplot+ [line width = 1,
    solid,mark = none,
    color = black,
    mark size = 2.0,
    mark options = {
      fill = kth-darkgray,
      line width = 1,
      rotate = 0,
      solid
    },fill = kth-darkgray, forget plot]coordinates {
      (19.0, -13.112593618627468)
      (20.0, -13.112593618627468)
      (20.0, 0.0)
      (19.0, 0.0)
      (19.0, -13.112593618627468)
    };
    \addplot+ [line width = 1,
    solid,mark = none,
    color = black,
    mark size = 2.0,
    mark options = {
      fill = kth-darkgray,
      line width = 1,
      rotate = 0,
      solid
    },fill = kth-darkgray, forget plot]coordinates {
      (20.0, -13.112593618627468)
      (21.0, -13.112593618627468)
      (21.0, 0.0)
      (20.0, 0.0)
      (20.0, -13.112593618627468)
    };
    \addplot+ [line width = 1,
    solid,mark = none,
    color = black,
    mark size = 2.0,
    mark options = {
      fill = kth-darkgray,
      line width = 1,
      rotate = 0,
      solid
    },fill = kth-darkgray, forget plot]coordinates {
      (21.0, -13.112593618627468)
      (22.0, -13.112593618627468)
      (22.0, 0.0)
      (21.0, 0.0)
      (21.0, -13.112593618627468)
    };
    \addplot+ [line width = 1,
    solid,mark = none,
    color = black,
    mark size = 2.0,
    mark options = {
      fill = kth-darkgray,
      line width = 1,
      rotate = 0,
      solid
    },fill = kth-darkgray, forget plot]coordinates {
      (22.0, -13.112593618627468)
      (23.0, -13.112593618627468)
      (23.0, 0.0)
      (22.0, 0.0)
      (22.0, -13.112593618627468)
    };
    \addplot+ [line width = 1,
    solid,mark = none,
    color = black,
    mark size = 2.0,
    mark options = {
      fill = kth-darkgray,
      line width = 1,
      rotate = 0,
      solid
    },fill = kth-darkgray, forget plot]coordinates {
      (23.0, -13.112593618627468)
      (24.0, -13.112593618627468)
      (24.0, 0.0)
      (23.0, 0.0)
      (23.0, -13.112593618627468)
    };
    \addplot+ [line width = 1,
    solid,mark = none,
    color = black,
    mark size = 2.0,
    mark options = {
      fill = kth-darkgray,
      line width = 1,
      rotate = 0,
      solid
    },fill = kth-darkgray, forget plot]coordinates {
      (17.0, 80)
      (24.0, 80)
      (24.0, 90)
      (17.0, 90)
      (17.0, 80)
    };
    \addplot+ [line width = 1,
    solid,mark = none,
    color = black,
    mark size = 2.0,
    mark options = {
      fill = kth-blue,
      line width = 1,
      rotate = 0,
      solid
    },fill = kth-blue, forget plot]coordinates {
      (3.0, 0.0)
      (4.0, 0.0)
      (4.0, 22.649493192612173)
      (3.0, 22.649493192612173)
      (3.0, 0.0)
    };
    \addplot+ [line width = 1,
    solid,mark = none,
    color = black,
    mark size = 2.0,
    mark options = {
      fill = kth-blue,
      line width = 1,
      rotate = 0,
      solid
    },fill = kth-blue, forget plot]coordinates {
      (3.0, 22.649493192612173)
      (4.0, 22.649493192612173)
      (4.0, 25.48140403175354)
      (3.0, 25.48140403175354)
      (3.0, 22.649493192612173)
    };
    \addplot+ [line width = 1,
    solid,mark = none,
    color = black,
    mark size = 2.0,
    mark options = {
      fill = kth-blue,
      line width = 1,
      rotate = 0,
      solid
    },fill = kth-blue, forget plot]coordinates {
      (3.0, 25.48140403175354)
      (4.0, 25.48140403175354)
      (4.0, 44.25781932866374)
      (3.0, 44.25781932866374)
      (3.0, 25.48140403175354)
    };
    \addplot+ [line width = 1,
    solid,mark = none,
    color = black,
    mark size = 2.0,
    mark options = {
      fill = kth-blue,
      line width = 1,
      rotate = 0,
      solid
    },fill = kth-blue, forget plot]coordinates {
      (13.0, 0.0)
      (14.0, 0.0)
      (14.0, 31.51378666114807)
      (13.0, 31.51378666114807)
      (13.0, 0.0)
    };
    \addplot+ [line width = 1,
    solid,mark = none,
    color = black,
    mark size = 2.0,
    mark options = {
      fill = kth-blue,
      line width = 1,
      rotate = 0,
      solid
    },fill = kth-blue, forget plot]coordinates {
      (13.0, 31.51378666114807)
      (14.0, 31.51378666114807)
      (14.0, 57.23391661133232)
      (13.0, 57.23391661133232)
      (13.0, 31.51378666114807)
    };
    \addplot+ [line width = 1,
    solid,mark = none,
    color = black,
    mark size = 2.0,
    mark options = {
      fill = kth-blue,
      line width = 1,
      rotate = 0,
      solid
    },fill = kth-blue, forget plot]coordinates {
      (17.0, 90)
      (24.0, 90)
      (24.0, 100)
      (17.0, 100)
      (17.0, 90)
    };
    \node at (axis cs:0.5, -6.556296809313734) [,
    color=white,
    rotate=0.0,
    font={\fontsize{6 pt}{15.600000000000001 pt}\bfseries\sffamily}
    ] {674.2};
    \node at (axis cs:1.5, -6.556296809313734) [,
    color=white,
    rotate=0.0,
    font={\fontsize{6 pt}{15.600000000000001 pt}\bfseries\sffamily}
    ] {646.6};
    \node at (axis cs:2.5, -6.556296809313734) [,
    color=white,
    rotate=0.0,
    font={\fontsize{6 pt}{15.600000000000001 pt}\bfseries\sffamily}
    ] {626.1};
    \node at (axis cs:3.5, 11.324746596306086) [,
    color=white,
    rotate=0.0,
    font={\fontsize{6 pt}{15.600000000000001 pt}\bfseries\sffamily}
    ] {0.0};
    \node at (axis cs:3.5, 24.065448612182855) [,
    color=white,
    rotate=0.0,
    font={\fontsize{6 pt}{15.600000000000001 pt}\bfseries\sffamily}
    ] {614.3};
    \node at (axis cs:3.5, 34.86961168020864) [,
    color=white,
    rotate=0.0,
    font={\fontsize{6 pt}{15.600000000000001 pt}\bfseries\sffamily}
    ] {2000};
    \node at (axis cs:4.5, -6.556296809313734) [,
    color=white,
    rotate=0.0,
    font={\fontsize{6 pt}{15.600000000000001 pt}\bfseries\sffamily}
    ] {686.6};
    \node at (axis cs:5.5, -6.556296809313734) [,
    color=white,
    rotate=0.0,
    font={\fontsize{6 pt}{15.600000000000001 pt}\bfseries\sffamily}
    ] {704.0};
    \node at (axis cs:6.5, -6.556296809313734) [,
    color=white,
    rotate=0.0,
    font={\fontsize{6 pt}{15.600000000000001 pt}\bfseries\sffamily}
    ] {735.3};
    \node at (axis cs:7.5, -6.556296809313734) [,
    color=white,
    rotate=0.0,
    font={\fontsize{6 pt}{15.600000000000001 pt}\bfseries\sffamily}
    ] {735.4};
    \node at (axis cs:8.5, -6.556296809313734) [,
    color=white,
    rotate=0.0,
    font={\fontsize{6 pt}{15.600000000000001 pt}\bfseries\sffamily}
    ] {924.5};
    \node at (axis cs:9.5, -6.556296809313734) [,
    color=white,
    rotate=0.0,
    font={\fontsize{6 pt}{15.600000000000001 pt}\bfseries\sffamily}
    ] {889.3};
    \node at (axis cs:10.5, -6.556296809313734) [,
    color=white,
    rotate=0.0,
    font={\fontsize{6 pt}{15.600000000000001 pt}\bfseries\sffamily}
    ] {803.5};
    \node at (axis cs:11.5, -6.556296809313734) [,
    color=white,
    rotate=0.0,
    font={\fontsize{6 pt}{15.600000000000001 pt}\bfseries\sffamily}
    ] {930.8};
    \node at (axis cs:12.5, -6.556296809313734) [,
    color=white,
    rotate=0.0,
    font={\fontsize{6 pt}{15.600000000000001 pt}\bfseries\sffamily}
    ] {853.3};
    \node at (axis cs:13.5, -6.556296809313734) [,
    color=white,
    rotate=0.0,
    font={\fontsize{6 pt}{15.600000000000001 pt}\bfseries\sffamily}
    ] {831.0};
    \node at (axis cs:13.5, 15.756893330574036) [,
    color=white,
    rotate=0.0,
    font={\fontsize{6 pt}{15.600000000000001 pt}\bfseries\sffamily}
    ] {0.0};
    \node at (axis cs:13.5, 44.3738516362402) [,
    color=white,
    rotate=0.0,
    font={\fontsize{6 pt}{15.600000000000001 pt}\bfseries\sffamily}
    ] {1200};
    \node at (axis cs:14.5, -6.556296809313734) [,
    color=white,
    rotate=0.0,
    font={\fontsize{6 pt}{15.600000000000001 pt}\bfseries\sffamily}
    ] {862.6};
    \node at (axis cs:15.5, -6.556296809313734) [,
    color=white,
    rotate=0.0,
    font={\fontsize{6 pt}{15.600000000000001 pt}\bfseries\sffamily}
    ] {908.6};
    \node at (axis cs:16.5, -6.556296809313734) [,
    color=white,
    rotate=0.0,
    font={\fontsize{6 pt}{15.600000000000001 pt}\bfseries\sffamily}
    ] {0.15};
    \node at (axis cs:18.5, -6.556296809313734) [,
    color=white,
    rotate=0.0,
    font={\fontsize{6 pt}{15.600000000000001 pt}\bfseries\sffamily}
    ] {14.67};
    \node at (axis cs:19.5, -6.556296809313734) [,
    color=white,
    rotate=0.0,
    font={\fontsize{6 pt}{15.600000000000001 pt}\bfseries\sffamily}
    ] {980.7};
    \node at (axis cs:20.5, -6.556296809313734) [,
    color=white,
    rotate=0.0,
    font={\fontsize{6 pt}{15.600000000000001 pt}\bfseries\sffamily}
    ] {980.4};
    \node at (axis cs:21.5, -6.556296809313734) [,
    color=white,
    rotate=0.0,
    font={\fontsize{6 pt}{15.600000000000001 pt}\bfseries\sffamily}
    ] {1000};
    \node at (axis cs:22.5, -6.556296809313734) [,
    color=white,
    rotate=0.0,
    font={\fontsize{6 pt}{15.600000000000001 pt}\bfseries\sffamily}
    ] {1000};
    \node at (axis cs:23.5, -6.556296809313734) [,
    color=white,
    rotate=0.0,
    font={\fontsize{6 pt}{15.600000000000001 pt}\bfseries\sffamily}
    ] {1000};
    \node at (axis cs:21.0, 85.0) [,
    color=white,
    rotate=0.0,
    font={\fontsize{8 pt}{15.600000000000001 pt}\bfseries\sffamily}
    ] {Independent Volumes [MWh]};
    \node at (axis cs:21.0, 95.0) [,
    color=white,
    rotate=0.0,
    font={\fontsize{8 pt}{15.600000000000001 pt}\bfseries\sffamily}
    ] {Dependent Volumes [MWh]};
  \end{axis}
  \begin{axis}[
    width = {\textwidth},
    height = {0.3\textheight},
    xlabel = {Hour},
    ylabel = {Price [Eur/MWh]},
    title = {Block orders},
    unbounded coords=jump,
    scaled x ticks = false,
    xmajorgrids = true,
    xmin = 0,
    xmax = 24,
    xtick align = inside,
    xtick = {0.0,1.0,2.0,3.0,4.0,5.0,6.0,7.0,8.0,9.0,10.0,11.0,12.0,13.0,14.0,15.0,16.0,17.0,18.0,19.0,20.0,21.0,22.0,23.0,24.0},
    xticklabels = {$0$,$1$,$2$,$3$,$4$,$5$,$6$,$7$,$8$,$9$,$10$,$11$,$12$,$13$,$14$,$15$,$16$,$17$,$18$,$19$,$20$,$21$,$22$,$23$,$24$},
    xticklabel style = {font = {\fontsize{10 pt}{10.4 pt}\selectfont}, rotate = 0.0},
    x grid style = {color = kth-lightgray,
      line width = 0.25,
      solid},
    axis x line* = left,
    x axis line style = {line width = 1,
      solid},
    scaled y ticks = false,
    ylabel style = {yshift = 5pt, font = {\fontsize{10 pt}{14.3 pt}\selectfont}, rotate = 0.0},
    ymajorgrids = true,
    ytick align = inside,
    yticklabel style = {font = {\fontsize{10 pt}{10.4 pt}\selectfont}, rotate = 0.0},
    y grid style = {color = kth-lightgray,
      line width = 0.25,
      solid},axis y line* = left,y axis line style = {line width = 1,
      solid},
    xshift = 0mm,
    title style = {font = {\fontsize{14 pt}{18.2 pt}\selectfont}, rotate = 0.0},
    legend pos = {north west},
    legend style = {line width = 2, solid,font = {\fontsize{8 pt}{10.4 pt}\selectfont}},colorbar style={title=}]
    \addplot+ [line width = 1,
    solid,mark = none,
    color = black,
    mark size = 2.0,
    mark options = {
      fill = kth-blue,
      line width = 1,
      rotate = 0,
      solid
    },fill = kth-blue, forget plot]coordinates {
      (17.0, 4.0)
      (24.0, 4.0)
      (24.0, 5.0)
      (17.0, 5.0)
      (17.0, 4.0)
    };
    \addplot+ [line width = 1,
    solid,mark = none,
    color = black,
    mark size = 2.0,
    mark options = {
      fill = kth-blue,
      line width = 1,
      rotate = 0,
      solid
    },fill = kth-blue, forget plot]coordinates {
      (16.0, 0.5)
      (19.0, 0.5)
      (19.0, 1.5)
      (16.0, 1.5)
      (16.0, 0.5)
    };
    \addplot+ [line width = 1,
    solid,mark = none,
    color = black,
    mark size = 2.0,
    mark options = {
      fill = kth-blue,
      line width = 1,
      rotate = 0,
      solid
    },fill = kth-blue, forget plot]coordinates {
      (16.0, 1.5)
      (19.0, 1.5)
      (19.0, 2.5)
      (16.0, 2.5)
      (16.0, 1.5)
    };
    \node at (axis cs:16.5, 1) [,
    color=white,
    rotate=0.0,
    font={\fontsize{6 pt}{18.2 pt}\bfseries\sffamily}
    ] {12.54};
    \node at (axis cs:18.5, 1) [,
    color=white,
    rotate=0.0,
    font={\fontsize{6 pt}{18.2 pt}\bfseries\sffamily}
    ] {496.3};
    \node at (axis cs:16.5, 2) [,
    color=white,
    rotate=0.0,
    font={\fontsize{6 pt}{18.2 pt}\bfseries\sffamily}
    ] {23.89};
    \node at (axis cs:18.5, 2) [,
    color=white,
    rotate=0.0,
    font={\fontsize{6 pt}{18.2 pt}\bfseries\sffamily}
    ] {500.0};
    \node at (axis cs:18.5, 4.5) [,
    color=white,
    rotate=0.0,
    font={\fontsize{6 pt}{18.2 pt}\bfseries\sffamily}
    ] {Price [EUR/MWh]};
    \node at (axis cs:22.6, 4.5) [,
    color=white,
    rotate=0.0,
    font={\fontsize{6 pt}{18.2 pt}\bfseries\sffamily}
    ] {Volume [MWh/h]};
  \end{axis}
\end{tikzpicture}

%% file: maintenance_confidence_intervals.tex
\begin{tikzpicture}[]
  \begin{axis}[
    width = {\textwidth},
    height = {0.4\textheight},
    xlabel = {Number of Samples $N$},
    ylabel = {Confidence interval [KEur]},
    title = {Maintenance scheduling confidence intervals},
    unbounded coords=jump,
    scaled x ticks = false,
    xmajorgrids = true,
    xmin = -50,
    xtick = {0.0,100,500.0,1000.0},
    xticklabels = {$0$,$100$,$500$,$1000$},
    xtick align = inside,
    xticklabel style = {font = {\fontsize{10 pt}{10.4 pt}\selectfont}, rotate = 0.0},
    x grid style = {color = kth-lightgray,
      line width = 0.25,
      solid},
    axis x line* = left,
    x axis line style = {line width = 1,
      solid},
    scaled y ticks = false,
    ylabel style = {font = {\fontsize{10 pt}{14.3 pt}\selectfont}, rotate = 0.0},
    ymajorgrids = true,
    ytick align = inside,
    yticklabel style = {font = {\fontsize{10 pt}{10.4 pt}\selectfont}, rotate = 0.0},
    y grid style = {color = kth-lightgray,
      line width = 0.25,
      solid},axis y line* = left,y axis line style = {line width = 1,
      solid},
    xshift = 0mm,
    title style = {font = {\fontsize{14 pt}{18.2 pt}\selectfont}, rotate = 0.0},
    legend pos = {north east},
    legend style = {line width = 2, solid,font = {\fontsize{8 pt}{10.4 pt}\selectfont}},colorbar style={title=}]
    \addplot+[draw=none, color = kth-blue,
    line width = 0,
    solid,mark = *,
    mark size = 2.0,
    mark options = {
      color = black,
      fill = kth-blue,
      line width = 1,
      rotate = 0,
      solid
    }] coordinates {
      (100.0, 564.36446)
      (500.0, 561.905815)
      (1000.0, 563.41468)
    };
    \addlegendentry{VRP confidence intervals}
    \addplot+[draw=none, color = kth-red,
    line width = 0,
    solid,mark = *,
    mark size = 2.0,
    mark options = {
      color = black,
      fill = kth-red,
      line width = 1,
      rotate = 0,
      solid
    }] coordinates {
      (1000.0, 551.770485)
    };
    \addlegendentry{EEV confidence interval}
    \addplot+ [color = kth-blue,
    line width = 1,
    solid,mark = -,
    mark size = 2.0,
    mark options = {
      color = black,
      line width = 1,
      rotate = 0,
      solid
    },forget plot]coordinates {
      (100.0, 551.13346)
      (100.0, 577.59546)
    };
    \addplot+ [color = kth-blue,
    line width = 1,
    solid,mark = -,
    mark size = 2.0,
    mark options = {
      color = black,
      line width = 1,
      rotate = 0,
      solid
    },forget plot]coordinates {
      (500.0, 558.63425)
      (500.0, 565.17738)
    };
    \addplot+ [color = kth-blue,
    line width = 1,
    solid,mark = -,
    mark size = 2.0,
    mark options = {
      color = black,
      line width = 1,
      rotate = 0,
      solid
    },forget plot]coordinates {
      (1000.0, 562.98567)
      (1000.0, 563.84369)
    };
    \addplot+ [color = kth-red,
    line width = 1,
    solid,mark = -,
    mark size = 2.0,
    mark options = {
      color = black,
      line width = 1,
      rotate = 0,
      solid
    },forget plot]coordinates {
      (1000.0, 550.42071)
      (1000.0, 553.12026)
    };
  \end{axis}
\end{tikzpicture}

%% file: maintenance_evp_strategy.tex
\begin{tikzpicture}[]
  \begin{axis}[
    width = {\textwidth},
    height = {0.3\textheight},
    xlabel = {},
    ylabel = {Price [Eur/MWh]},
    title = {Single orders},
    unbounded coords=jump,
    scaled x ticks = false,
    xmajorgrids = true,
    xmin = 0,
    xmax = 24,
    ymin = {-13.112593618627468},
    ymax = {102.22033269617995},
    xtick align = inside,
    xtick = {0.0,1.0,2.0,3.0,4.0,5.0,6.0,7.0,8.0,9.0,10.0,11.0,12.0,13.0,14.0,15.0,16.0,17.0,18.0,19.0,20.0,21.0,22.0,23.0,24.0},
    xticklabels = {$0$,$1$,$2$,$3$,$4$,$5$,$6$,$7$,$8$,$9$,$10$,$11$,$12$,$13$,$14$,$15$,$16$,$17$,$18$,$19$,$20$,$21$,$22$,$23$,$24$},
    xticklabel style = {font = {\fontsize{10 pt}{10.4 pt}\selectfont}, rotate = 0.0},
    x grid style = {color = kth-lightgray,
      line width = 0.25,
      solid},
    axis x line* = left,
    x axis line style = {line width = 1,
      solid},
    scaled y ticks = false,
    ylabel style = {yshift = -5pt, font = {\fontsize{10 pt}{14.3 pt}\selectfont}, rotate = 0.0},
    ymajorgrids = true,
    ytick align = inside,
    yticklabel style = {font = {\fontsize{10 pt}{10.4 pt}\selectfont}, rotate = 0.0},
    y grid style = {color = kth-lightgray,
      line width = 0.25,
      solid},axis y line* = left,y axis line style = {line width = 1,
      solid},
    xshift = 0mm,
    yshift = {0.3\textheight},
    title style = {font = {\fontsize{14 pt}{18.2 pt}\selectfont}, rotate = 0.0},
    legend pos = {north west},
    legend style = {line width = 2, solid,font = {\fontsize{8 pt}{10.4 pt}\selectfont}},colorbar style={title=}]
    \addplot+ [line width = 1,
    solid,mark = none,
    color = black,
    mark size = 2.0,
    mark options = {
      color = black,
      fill = kth-darkgray,
      line width = 1,
      rotate = 0,
      solid
    },fill = kth-darkgray, forget plot]coordinates {
      (0.0, -14.298777702158802)
      (1.0, -14.298777702158802)
      (1.0, 0.0)
      (0.0, 0.0)
      (0.0, -14.298777702158802)
    };
    \addplot+ [line width = 1,
    solid,mark = none,
    color = black,
    mark size = 2.0,
    mark options = {
      color = black,
      fill = kth-darkgray,
      line width = 1,
      rotate = 0,
      solid
    },fill = kth-darkgray, forget plot]coordinates {
      (1.0, -14.298777702158802)
      (2.0, -14.298777702158802)
      (2.0, 0.0)
      (1.0, 0.0)
      (1.0, -14.298777702158802)
    };
    \addplot+ [line width = 1,
    solid,mark = none,
    color = black,
    mark size = 2.0,
    mark options = {
      color = black,
      fill = kth-darkgray,
      line width = 1,
      rotate = 0,
      solid
    },fill = kth-darkgray, forget plot]coordinates {
      (3.0, -14.298777702158802)
      (4.0, -14.298777702158802)
      (4.0, 0.0)
      (3.0, 0.0)
      (3.0, -14.298777702158802)
    };
    \addplot+ [line width = 1,
    solid,mark = none,
    color = black,
    mark size = 2.0,
    mark options = {
      color = black,
      fill = kth-darkgray,
      line width = 1,
      rotate = 0,
      solid
    },fill = kth-darkgray, forget plot]coordinates {
      (4.0, -14.298777702158802)
      (5.0, -14.298777702158802)
      (5.0, 0.0)
      (4.0, 0.0)
      (4.0, -14.298777702158802)
    };
    \addplot+ [line width = 1,
    solid,mark = none,
    color = black,
    mark size = 2.0,
    mark options = {
      color = black,
      fill = kth-darkgray,
      line width = 1,
      rotate = 0,
      solid
    },fill = kth-darkgray, forget plot]coordinates {
      (5.0, -14.298777702158802)
      (6.0, -14.298777702158802)
      (6.0, 0.0)
      (5.0, 0.0)
      (5.0, -14.298777702158802)
    };
    \addplot+ [line width = 1,
    solid,mark = none,
    color = black,
    mark size = 2.0,
    mark options = {
      color = black,
      fill = kth-darkgray,
      line width = 1,
      rotate = 0,
      solid
    },fill = kth-darkgray, forget plot]coordinates {
      (6.0, -14.298777702158802)
      (7.0, -14.298777702158802)
      (7.0, 0.0)
      (6.0, 0.0)
      (6.0, -14.298777702158802)
    };
    \addplot+ [line width = 1,
    solid,mark = none,
    color = black,
    mark size = 2.0,
    mark options = {
      color = black,
      fill = kth-darkgray,
      line width = 1,
      rotate = 0,
      solid
    },fill = kth-darkgray, forget plot]coordinates {
      (7.0, -14.298777702158802)
      (8.0, -14.298777702158802)
      (8.0, 0.0)
      (7.0, 0.0)
      (7.0, -14.298777702158802)
    };
    \addplot+ [line width = 1,
    solid,mark = none,
    color = black,
    mark size = 2.0,
    mark options = {
      color = black,
      fill = kth-darkgray,
      line width = 1,
      rotate = 0,
      solid
    },fill = kth-darkgray, forget plot]coordinates {
      (8.0, -14.298777702158802)
      (9.0, -14.298777702158802)
      (9.0, 0.0)
      (8.0, 0.0)
      (8.0, -14.298777702158802)
    };
    \addplot+ [line width = 1,
    solid,mark = none,
    color = black,
    mark size = 2.0,
    mark options = {
      color = black,
      fill = kth-darkgray,
      line width = 1,
      rotate = 0,
      solid
    },fill = kth-darkgray, forget plot]coordinates {
      (9.0, -14.298777702158802)
      (10.0, -14.298777702158802)
      (10.0, 0.0)
      (9.0, 0.0)
      (9.0, -14.298777702158802)
    };
    \addplot+ [line width = 1,
    solid,mark = none,
    color = black,
    mark size = 2.0,
    mark options = {
      color = black,
      fill = kth-darkgray,
      line width = 1,
      rotate = 0,
      solid
    },fill = kth-darkgray, forget plot]coordinates {
      (10.0, -14.298777702158802)
      (11.0, -14.298777702158802)
      (11.0, 0.0)
      (10.0, 0.0)
      (10.0, -14.298777702158802)
    };
    \addplot+ [line width = 1,
    solid,mark = none,
    color = black,
    mark size = 2.0,
    mark options = {
      color = black,
      fill = kth-darkgray,
      line width = 1,
      rotate = 0,
      solid
    },fill = kth-darkgray, forget plot]coordinates {
      (11.0, -14.298777702158802)
      (12.0, -14.298777702158802)
      (12.0, 0.0)
      (11.0, 0.0)
      (11.0, -14.298777702158802)
    };
    \addplot+ [line width = 1,
    solid,mark = none,
    color = black,
    mark size = 2.0,
    mark options = {
      color = black,
      fill = kth-darkgray,
      line width = 1,
      rotate = 0,
      solid
    },fill = kth-darkgray, forget plot]coordinates {
      (12.0, -14.298777702158802)
      (13.0, -14.298777702158802)
      (13.0, 0.0)
      (12.0, 0.0)
      (12.0, -14.298777702158802)
    };
    \addplot+ [line width = 1,
    solid,mark = none,
    color = black,
    mark size = 2.0,
    mark options = {
      color = black,
      fill = kth-darkgray,
      line width = 1,
      rotate = 0,
      solid
    },fill = kth-darkgray, forget plot]coordinates {
      (13.0, -14.298777702158802)
      (14.0, -14.298777702158802)
      (14.0, 0.0)
      (13.0, 0.0)
      (13.0, -14.298777702158802)
    };
    \addplot+ [line width = 1,
    solid,mark = none,
    color = black,
    mark size = 2.0,
    mark options = {
      color = black,
      fill = kth-darkgray,
      line width = 1,
      rotate = 0,
      solid
    },fill = kth-darkgray, forget plot]coordinates {
      (14.0, -14.298777702158802)
      (15.0, -14.298777702158802)
      (15.0, 0.0)
      (14.0, 0.0)
      (14.0, -14.298777702158802)
    };
    \addplot+ [line width = 1,
    solid,mark = none,
    color = black,
    mark size = 2.0,
    mark options = {
      color = black,
      fill = kth-darkgray,
      line width = 1,
      rotate = 0,
      solid
    },fill = kth-darkgray, forget plot]coordinates {
      (15.0, -14.298777702158802)
      (16.0, -14.298777702158802)
      (16.0, 0.0)
      (15.0, 0.0)
      (15.0, -14.298777702158802)
    };
    \addplot+ [line width = 1,
    solid,mark = none,
    color = black,
    mark size = 2.0,
    mark options = {
      color = black,
      fill = kth-darkgray,
      line width = 1,
      rotate = 0,
      solid
    },fill = kth-darkgray, forget plot]coordinates {
      (16.0, -14.298777702158802)
      (17.0, -14.298777702158802)
      (17.0, 0.0)
      (16.0, 0.0)
      (16.0, -14.298777702158802)
    };
    \addplot+ [line width = 1,
    solid,mark = none,
    color = black,
    mark size = 2.0,
    mark options = {
      color = black,
      fill = kth-darkgray,
      line width = 1,
      rotate = 0,
      solid
    },fill = kth-darkgray, forget plot]coordinates {
      (18.0, -14.298777702158802)
      (19.0, -14.298777702158802)
      (19.0, 0.0)
      (18.0, 0.0)
      (18.0, -14.298777702158802)
    };
    \addplot+ [line width = 1,
    solid,mark = none,
    color = black,
    mark size = 2.0,
    mark options = {
      color = black,
      fill = kth-darkgray,
      line width = 1,
      rotate = 0,
      solid
    },fill = kth-darkgray, forget plot]coordinates {
      (19.0, -14.298777702158802)
      (20.0, -14.298777702158802)
      (20.0, 0.0)
      (19.0, 0.0)
      (19.0, -14.298777702158802)
    };
    \addplot+ [line width = 1,
    solid,mark = none,
    color = black,
    mark size = 2.0,
    mark options = {
      color = black,
      fill = kth-darkgray,
      line width = 1,
      rotate = 0,
      solid
    },fill = kth-darkgray, forget plot]coordinates {
      (20.0, -14.298777702158802)
      (21.0, -14.298777702158802)
      (21.0, 0.0)
      (20.0, 0.0)
      (20.0, -14.298777702158802)
    };
    \addplot+ [line width = 1,
    solid,mark = none,
    color = black,
    mark size = 2.0,
    mark options = {
      color = black,
      fill = kth-darkgray,
      line width = 1,
      rotate = 0,
      solid
    },fill = kth-darkgray, forget plot]coordinates {
      (21.0, -14.298777702158802)
      (22.0, -14.298777702158802)
      (22.0, 0.0)
      (21.0, 0.0)
      (21.0, -14.298777702158802)
    };
    \addplot+ [line width = 1,
    solid,mark = none,
    color = black,
    mark size = 2.0,
    mark options = {
      color = black,
      fill = kth-darkgray,
      line width = 1,
      rotate = 0,
      solid
    },fill = kth-darkgray, forget plot]coordinates {
      (22.0, -14.298777702158802)
      (23.0, -14.298777702158802)
      (23.0, 0.0)
      (22.0, 0.0)
      (22.0, -14.298777702158802)
    };
    \addplot+ [line width = 1,
    solid,mark = none,
    color = black,
    mark size = 2.0,
    mark options = {
      color = black,
      fill = kth-darkgray,
      line width = 1,
      rotate = 0,
      solid
    },fill = kth-darkgray, forget plot]coordinates {
      (23.0, -14.298777702158802)
      (24.0, -14.298777702158802)
      (24.0, 0.0)
      (23.0, 0.0)
      (23.0, -14.298777702158802)
    };
    \addplot+ [line width = 1,
    solid,mark = none,
    color = black,
    mark size = 2.0,
    mark options = {
      fill = kth-darkgray,
      line width = 1,
      rotate = 0,
      solid
    },fill = kth-darkgray, forget plot]coordinates {
      (17.0, 80)
      (24.0, 80)
      (24.0, 90)
      (17.0, 90)
      (17.0, 80)
    };
    \addplot+ [line width = 1,
    solid,mark = none,
    color = black,
    mark size = 2.0,
    mark options = {
      color = black,
      fill = kth-blue,
      line width = 1,
      rotate = 0,
      solid
    },fill = kth-blue, forget plot]coordinates {
      (0.0, 0.0)
      (1.0, 0.0)
      (1.0, 26.454528621673585)
      (0.0, 26.454528621673585)
      (0.0, 0.0)
    };
    \addplot+ [line width = 1,
    solid,mark = none,
    color = black,
    mark size = 2.0,
    mark options = {
      color = black,
      fill = kth-blue,
      line width = 1,
      rotate = 0,
      solid
    },fill = kth-blue, forget plot]coordinates {
      (0.0, 26.454528621673585)
      (1.0, 26.454528621673585)
      (1.0, 46.33957073882692)
      (0.0, 46.33957073882692)
      (0.0, 26.454528621673585)
    };
    \addplot+ [line width = 1,
    solid,mark = none,
    color = black,
    mark size = 2.0,
    mark options = {
      color = black,
      fill = kth-blue,
      line width = 1,
      rotate = 0,
      solid
    },fill = kth-blue, forget plot]coordinates {
      (1.0, 0.0)
      (2.0, 0.0)
      (2.0, 25.999015907287596)
      (1.0, 25.999015907287596)
      (1.0, 0.0)
    };
    \addplot+ [line width = 1,
    solid,mark = none,
    color = black,
    mark size = 2.0,
    mark options = {
      color = black,
      fill = kth-blue,
      line width = 1,
      rotate = 0,
      solid
    },fill = kth-blue, forget plot]coordinates {
      (1.0, 25.999015907287596)
      (2.0, 25.999015907287596)
      (2.0, 45.97610227835858)
      (1.0, 45.97610227835858)
      (1.0, 25.999015907287596)
    };
    \addplot+ [line width = 1,
    solid,mark = none,
    color = black,
    mark size = 2.0,
    mark options = {
      color = black,
      fill = kth-blue,
      line width = 1,
      rotate = 0,
      solid
    },fill = kth-blue, forget plot]coordinates {
      (2.0, 0.0)
      (3.0, 0.0)
      (3.0, 25.628022993087768)
      (2.0, 25.628022993087768)
      (2.0, 0.0)
    };
    \addplot+ [line width = 1,
    solid,mark = none,
    color = black,
    mark size = 2.0,
    mark options = {
      color = black,
      fill = kth-blue,
      line width = 1,
      rotate = 0,
      solid
    },fill = kth-blue, forget plot]coordinates {
      (2.0, 25.628022993087768)
      (3.0, 25.628022993087768)
      (3.0, 45.47467153842673)
      (2.0, 45.47467153842673)
      (2.0, 25.628022993087768)
    };
    \addplot+ [line width = 1,
    solid,mark = none,
    color = black,
    mark size = 2.0,
    mark options = {
      color = black,
      fill = kth-blue,
      line width = 1,
      rotate = 0,
      solid
    },fill = kth-blue, forget plot]coordinates {
      (3.0, 0.0)
      (4.0, 0.0)
      (4.0, 28.556312987609857)
      (3.0, 28.556312987609857)
      (3.0, 0.0)
    };
    \addplot+ [line width = 1,
    solid,mark = none,
    color = black,
    mark size = 2.0,
    mark options = {
      color = black,
      fill = kth-blue,
      line width = 1,
      rotate = 0,
      solid
    },fill = kth-blue, forget plot]coordinates {
      (3.0, 28.556312987609857)
      (4.0, 28.556312987609857)
      (4.0, 45.83366751233183)
      (3.0, 45.83366751233183)
      (3.0, 28.556312987609857)
    };
    \addplot+ [line width = 1,
    solid,mark = none,
    color = black,
    mark size = 2.0,
    mark options = {
      color = black,
      fill = kth-blue,
      line width = 1,
      rotate = 0,
      solid
    },fill = kth-blue, forget plot]coordinates {
      (4.0, 0.0)
      (5.0, 0.0)
      (5.0, 29.111444969903)
      (4.0, 29.111444969903)
      (4.0, 0.0)
    };
    \addplot+ [line width = 1,
    solid,mark = none,
    color = black,
    mark size = 2.0,
    mark options = {
      color = black,
      fill = kth-blue,
      line width = 1,
      rotate = 0,
      solid
    },fill = kth-blue, forget plot]coordinates {
      (4.0, 29.111444969903)
      (5.0, 29.111444969903)
      (5.0, 46.42157333692635)
      (4.0, 46.42157333692635)
      (4.0, 29.111444969903)
    };
    \addplot+ [line width = 1,
    solid,mark = none,
    color = black,
    mark size = 2.0,
    mark options = {
      color = black,
      fill = kth-blue,
      line width = 1,
      rotate = 0,
      solid
    },fill = kth-blue, forget plot]coordinates {
      (5.0, 0.0)
      (6.0, 0.0)
      (6.0, 30.102645218077797)
      (5.0, 30.102645218077797)
      (5.0, 0.0)
    };
    \addplot+ [line width = 1,
    solid,mark = none,
    color = black,
    mark size = 2.0,
    mark options = {
      color = black,
      fill = kth-blue,
      line width = 1,
      rotate = 0,
      solid
    },fill = kth-blue, forget plot]coordinates {
      (5.0, 30.102645218077797)
      (6.0, 30.102645218077797)
      (6.0, 47.598687105370665)
      (5.0, 47.598687105370665)
      (5.0, 30.102645218077797)
    };
    \addplot+ [line width = 1,
    solid,mark = none,
    color = black,
    mark size = 2.0,
    mark options = {
      color = black,
      fill = kth-blue,
      line width = 1,
      rotate = 0,
      solid
    },fill = kth-blue, forget plot]coordinates {
      (6.0, 0.0)
      (7.0, 0.0)
      (7.0, 28.450871245384217)
      (6.0, 28.450871245384217)
      (6.0, 0.0)
    };
    \addplot+ [line width = 1,
    solid,mark = none,
    color = black,
    mark size = 2.0,
    mark options = {
      color = black,
      fill = kth-blue,
      line width = 1,
      rotate = 0,
      solid
    },fill = kth-blue, forget plot]coordinates {
      (6.0, 28.450871245384217)
      (7.0, 28.450871245384217)
      (7.0, 51.78315439040725)
      (6.0, 51.78315439040725)
      (6.0, 28.450871245384217)
    };
    \addplot+ [line width = 1,
    solid,mark = none,
    color = black,
    mark size = 2.0,
    mark options = {
      color = black,
      fill = kth-blue,
      line width = 1,
      rotate = 0,
      solid
    },fill = kth-blue, forget plot]coordinates {
      (7.0, 0.0)
      (8.0, 0.0)
      (8.0, 32.88392045783996)
      (7.0, 32.88392045783996)
      (7.0, 0.0)
    };
    \addplot+ [line width = 1,
    solid,mark = none,
    color = black,
    mark size = 2.0,
    mark options = {
      color = black,
      fill = kth-blue,
      line width = 1,
      rotate = 0,
      solid
    },fill = kth-blue, forget plot]coordinates {
      (7.0, 32.88392045783996)
      (8.0, 32.88392045783996)
      (8.0, 65.48759754558898)
      (7.0, 65.48759754558898)
      (7.0, 32.88392045783996)
    };
    \addplot+ [line width = 1,
    solid,mark = none,
    color = black,
    mark size = 2.0,
    mark options = {
      color = black,
      fill = kth-blue,
      line width = 1,
      rotate = 0,
      solid
    },fill = kth-blue, forget plot]coordinates {
      (8.0, 0.0)
      (9.0, 0.0)
      (9.0, 35.40428360366821)
      (8.0, 35.40428360366821)
      (8.0, 0.0)
    };
    \addplot+ [line width = 1,
    solid,mark = none,
    color = black,
    mark size = 2.0,
    mark options = {
      color = black,
      fill = kth-blue,
      line width = 1,
      rotate = 0,
      solid
    },fill = kth-blue, forget plot]coordinates {
      (8.0, 35.40428360366821)
      (9.0, 35.40428360366821)
      (9.0, 74.15852349325041)
      (8.0, 74.15852349325041)
      (8.0, 35.40428360366821)
    };
    \addplot+ [line width = 1,
    solid,mark = none,
    color = black,
    mark size = 2.0,
    mark options = {
      color = black,
      fill = kth-blue,
      line width = 1,
      rotate = 0,
      solid
    },fill = kth-blue, forget plot]coordinates {
      (9.0, 0.0)
      (10.0, 0.0)
      (10.0, 34.756728799819946)
      (9.0, 34.756728799819946)
      (9.0, 0.0)
    };
    \addplot+ [line width = 1,
    solid,mark = none,
    color = black,
    mark size = 2.0,
    mark options = {
      color = black,
      fill = kth-blue,
      line width = 1,
      rotate = 0,
      solid
    },fill = kth-blue, forget plot]coordinates {
      (9.0, 34.756728799819946)
      (10.0, 34.756728799819946)
      (10.0, 70.60435897172817)
      (9.0, 70.60435897172817)
      (9.0, 34.756728799819946)
    };
    \addplot+ [line width = 1,
    solid,mark = none,
    color = black,
    mark size = 2.0,
    mark options = {
      color = black,
      fill = kth-blue,
      line width = 1,
      rotate = 0,
      solid
    },fill = kth-blue, forget plot]coordinates {
      (10.0, 0.0)
      (11.0, 0.0)
      (11.0, 33.70945800590515)
      (10.0, 33.70945800590515)
      (10.0, 0.0)
    };
    \addplot+ [line width = 1,
    solid,mark = none,
    color = black,
    mark size = 2.0,
    mark options = {
      color = black,
      fill = kth-blue,
      line width = 1,
      rotate = 0,
      solid
    },fill = kth-blue, forget plot]coordinates {
      (10.0, 33.70945800590515)
      (11.0, 33.70945800590515)
      (11.0, 64.78442087187007)
      (10.0, 64.78442087187007)
      (10.0, 33.70945800590515)
    };
    \addplot+ [line width = 1,
    solid,mark = none,
    color = black,
    mark size = 2.0,
    mark options = {
      color = black,
      fill = kth-blue,
      line width = 1,
      rotate = 0,
      solid
    },fill = kth-blue, forget plot]coordinates {
      (11.0, 0.0)
      (12.0, 0.0)
      (12.0, 32.78440212631226)
      (11.0, 32.78440212631226)
      (11.0, 0.0)
    };
    \addplot+ [line width = 1,
    solid,mark = none,
    color = black,
    mark size = 2.0,
    mark options = {
      color = black,
      fill = kth-blue,
      line width = 1,
      rotate = 0,
      solid
    },fill = kth-blue, forget plot]coordinates {
      (11.0, 32.78440212631226)
      (12.0, 32.78440212631226)
      (12.0, 61.706004670491794)
      (11.0, 61.706004670491794)
      (11.0, 32.78440212631226)
    };
    \addplot+ [line width = 1,
    solid,mark = none,
    color = black,
    mark size = 2.0,
    mark options = {
      color = black,
      fill = kth-blue,
      line width = 1,
      rotate = 0,
      solid
    },fill = kth-blue, forget plot]coordinates {
      (12.0, 0.0)
      (13.0, 0.0)
      (13.0, 32.284474506378174)
      (12.0, 32.284474506378174)
      (12.0, 0.0)
    };
    \addplot+ [line width = 1,
    solid,mark = none,
    color = black,
    mark size = 2.0,
    mark options = {
      color = black,
      fill = kth-blue,
      line width = 1,
      rotate = 0,
      solid
    },fill = kth-blue, forget plot]coordinates {
      (12.0, 32.284474506378174)
      (13.0, 32.284474506378174)
      (13.0, 60.68693440226748)
      (12.0, 60.68693440226748)
      (12.0, 32.284474506378174)
    };
    \addplot+ [line width = 1,
    solid,mark = none,
    color = black,
    mark size = 2.0,
    mark options = {
      color = black,
      fill = kth-blue,
      line width = 1,
      rotate = 0,
      solid
    },fill = kth-blue, forget plot]coordinates {
      (13.0, 0.0)
      (14.0, 0.0)
      (14.0, 31.89378568458557)
      (13.0, 31.89378568458557)
      (13.0, 0.0)
    };
    \addplot+ [line width = 1,
    solid,mark = none,
    color = black,
    mark size = 2.0,
    mark options = {
      color = black,
      fill = kth-blue,
      line width = 1,
      rotate = 0,
      solid
    },fill = kth-blue, forget plot]coordinates {
      (13.0, 31.89378568458557)
      (14.0, 31.89378568458557)
      (14.0, 59.55796074885757)
      (13.0, 59.55796074885757)
      (13.0, 31.89378568458557)
    };
    \addplot+ [line width = 1,
    solid,mark = none,
    color = black,
    mark size = 2.0,
    mark options = {
      color = black,
      fill = kth-blue,
      line width = 1,
      rotate = 0,
      solid
    },fill = kth-blue, forget plot]coordinates {
      (14.0, 0.0)
      (15.0, 0.0)
      (15.0, 31.81470869445801)
      (14.0, 31.81470869445801)
      (14.0, 0.0)
    };
    \addplot+ [line width = 1,
    solid,mark = none,
    color = black,
    mark size = 2.0,
    mark options = {
      color = black,
      fill = kth-blue,
      line width = 1,
      rotate = 0,
      solid
    },fill = kth-blue, forget plot]coordinates {
      (14.0, 31.81470869445801)
      (15.0, 31.81470869445801)
      (15.0, 59.921089365066706)
      (14.0, 59.921089365066706)
      (14.0, 31.81470869445801)
    };
    \addplot+ [line width = 1,
    solid,mark = none,
    color = black,
    mark size = 2.0,
    mark options = {
      color = black,
      fill = kth-blue,
      line width = 1,
      rotate = 0,
      solid
    },fill = kth-blue, forget plot]coordinates {
      (15.0, 0.0)
      (16.0, 0.0)
      (16.0, 32.46690121746063)
      (15.0, 32.46690121746063)
      (15.0, 0.0)
    };
    \addplot+ [line width = 1,
    solid,mark = none,
    color = black,
    mark size = 2.0,
    mark options = {
      color = black,
      fill = kth-blue,
      line width = 1,
      rotate = 0,
      solid
    },fill = kth-blue, forget plot]coordinates {
      (15.0, 32.46690121746063)
      (16.0, 32.46690121746063)
      (16.0, 63.13205644745099)
      (15.0, 63.13205644745099)
      (15.0, 32.46690121746063)
    };
    \addplot+ [line width = 1,
    solid,mark = none,
    color = black,
    mark size = 2.0,
    mark options = {
      color = black,
      fill = kth-blue,
      line width = 1,
      rotate = 0,
      solid
    },fill = kth-blue, forget plot]coordinates {
      (16.0, 0.0)
      (17.0, 0.0)
      (17.0, 48.17199633004927)
      (16.0, 48.17199633004927)
      (16.0, 0.0)
    };
    \addplot+ [line width = 1,
    solid,mark = none,
    color = black,
    mark size = 2.0,
    mark options = {
      color = black,
      fill = kth-blue,
      line width = 1,
      rotate = 0,
      solid
    },fill = kth-blue, forget plot]coordinates {
      (16.0, 48.17199633004927)
      (17.0, 48.17199633004927)
      (17.0, 75.43640349805293)
      (16.0, 75.43640349805293)
      (16.0, 48.17199633004927)
    };
    \addplot+ [line width = 1,
    solid,mark = none,
    color = black,
    mark size = 2.0,
    mark options = {
      color = black,
      fill = kth-blue,
      line width = 1,
      rotate = 0,
      solid
    },fill = kth-blue, forget plot]coordinates {
      (17.0, 0.0)
      (18.0, 0.0)
      (18.0, 23.286968313389906)
      (17.0, 23.286968313389906)
      (17.0, 0.0)
    };
    \addplot+ [line width = 1,
    solid,mark = none,
    color = black,
    mark size = 2.0,
    mark options = {
      color = black,
      fill = kth-blue,
      line width = 1,
      rotate = 0,
      solid
    },fill = kth-blue, forget plot]coordinates {
      (17.0, 23.286968313389906)
      (18.0, 23.286968313389906)
      (18.0, 37.58574601554871)
      (17.0, 37.58574601554871)
      (17.0, 23.286968313389906)
    };
    \addplot+ [line width = 1,
    solid,mark = none,
    color = black,
    mark size = 2.0,
    mark options = {
      color = black,
      fill = kth-blue,
      line width = 1,
      rotate = 0,
      solid
    },fill = kth-blue, forget plot]coordinates {
      (17.0, 37.58574601554871)
      (18.0, 37.58574601554871)
      (18.0, 80.48207912202511)
      (17.0, 80.48207912202511)
      (17.0, 37.58574601554871)
    };
    \addplot+ [line width = 1,
    solid,mark = none,
    color = black,
    mark size = 2.0,
    mark options = {
      color = black,
      fill = kth-blue,
      line width = 1,
      rotate = 0,
      solid
    },fill = kth-blue, forget plot]coordinates {
      (18.0, 0.0)
      (19.0, 0.0)
      (19.0, 35.18487405967712)
      (18.0, 35.18487405967712)
      (18.0, 0.0)
    };
    \addplot+ [line width = 1,
    solid,mark = none,
    color = black,
    mark size = 2.0,
    mark options = {
      color = black,
      fill = kth-blue,
      line width = 1,
      rotate = 0,
      solid
    },fill = kth-blue, forget plot]coordinates {
      (18.0, 35.18487405967712)
      (19.0, 35.18487405967712)
      (19.0, 68.04788266778002)
      (18.0, 68.04788266778002)
      (18.0, 35.18487405967712)
    };
    \addplot+ [line width = 1,
    solid,mark = none,
    color = black,
    mark size = 2.0,
    mark options = {
      color = black,
      fill = kth-blue,
      line width = 1,
      rotate = 0,
      solid
    },fill = kth-blue, forget plot]coordinates {
      (19.0, 0.0)
      (20.0, 0.0)
      (20.0, 31.96368113899231)
      (19.0, 31.96368113899231)
      (19.0, 0.0)
    };
    \addplot+ [line width = 1,
    solid,mark = none,
    color = black,
    mark size = 2.0,
    mark options = {
      color = black,
      fill = kth-blue,
      line width = 1,
      rotate = 0,
      solid
    },fill = kth-blue, forget plot]coordinates {
      (19.0, 31.96368113899231)
      (20.0, 31.96368113899231)
      (20.0, 56.95050131722712)
      (19.0, 56.95050131722712)
      (19.0, 31.96368113899231)
    };
    \addplot+ [line width = 1,
    solid,mark = none,
    color = black,
    mark size = 2.0,
    mark options = {
      color = black,
      fill = kth-blue,
      line width = 1,
      rotate = 0,
      solid
    },fill = kth-blue, forget plot]coordinates {
      (20.0, 0.0)
      (21.0, 0.0)
      (21.0, 29.962088035583495)
      (20.0, 29.962088035583495)
      (20.0, 0.0)
    };
    \addplot+ [line width = 1,
    solid,mark = none,
    color = black,
    mark size = 2.0,
    mark options = {
      color = black,
      fill = kth-blue,
      line width = 1,
      rotate = 0,
      solid
    },fill = kth-blue, forget plot]coordinates {
      (20.0, 29.962088035583495)
      (21.0, 29.962088035583495)
      (21.0, 51.55445568559523)
      (20.0, 51.55445568559523)
      (20.0, 29.962088035583495)
    };
    \addplot+ [line width = 1,
    solid,mark = none,
    color = black,
    mark size = 2.0,
    mark options = {
      color = black,
      fill = kth-blue,
      line width = 1,
      rotate = 0,
      solid
    },fill = kth-blue, forget plot]coordinates {
      (21.0, 0.0)
      (22.0, 0.0)
      (22.0, 28.77162281513214)
      (21.0, 28.77162281513214)
      (21.0, 0.0)
    };
    \addplot+ [line width = 1,
    solid,mark = none,
    color = black,
    mark size = 2.0,
    mark options = {
      color = black,
      fill = kth-blue,
      line width = 1,
      rotate = 0,
      solid
    },fill = kth-blue, forget plot]coordinates {
      (21.0, 28.77162281513214)
      (22.0, 28.77162281513214)
      (22.0, 49.32015411317227)
      (21.0, 49.32015411317227)
      (21.0, 28.77162281513214)
    };
    \addplot+ [line width = 1,
    solid,mark = none,
    color = black,
    mark size = 2.0,
    mark options = {
      color = black,
      fill = kth-blue,
      line width = 1,
      rotate = 0,
      solid
    },fill = kth-blue, forget plot]coordinates {
      (22.0, 0.0)
      (23.0, 0.0)
      (23.0, 27.807249476909636)
      (22.0, 27.807249476909636)
      (22.0, 0.0)
    };
    \addplot+ [line width = 1,
    solid,mark = none,
    color = black,
    mark size = 2.0,
    mark options = {
      color = black,
      fill = kth-blue,
      line width = 1,
      rotate = 0,
      solid
    },fill = kth-blue, forget plot]coordinates {
      (22.0, 27.807249476909636)
      (23.0, 27.807249476909636)
      (23.0, 48.032840452621016)
      (22.0, 48.032840452621016)
      (22.0, 27.807249476909636)
    };
    \addplot+ [line width = 1,
    solid,mark = none,
    color = black,
    mark size = 2.0,
    mark options = {
      color = black,
      fill = kth-blue,
      line width = 1,
      rotate = 0,
      solid
    },fill = kth-blue, forget plot]coordinates {
      (23.0, 0.0)
      (24.0, 0.0)
      (24.0, 26.83629267024994)
      (23.0, 26.83629267024994)
      (23.0, 0.0)
    };
    \addplot+ [line width = 1,
    solid,mark = none,
    color = black,
    mark size = 2.0,
    mark options = {
      fill = kth-blue,
      line width = 1,
      rotate = 0,
      solid
    },fill = kth-blue, forget plot]coordinates {
      (17.0, 90)
      (24.0, 90)
      (24.0, 100)
      (17.0, 100)
      (17.0, 90)
    };
    \node at (axis cs:0.5, -7.149388851079401) [,
    color=white,
    rotate=0.0,
    font={\fontsize{6 pt}{15.600000000000001 pt}\bfseries\sffamily}
    ] {684.3};
    \node at (axis cs:0.5, 13.227264310836793) [,
    color=white,
    rotate=0.0,
    font={\fontsize{6 pt}{15.600000000000001 pt}\bfseries\sffamily}
    ] {0.0};
    \node at (axis cs:0.5, 36.397049680250255) [,
    color=white,
    rotate=0.0,
    font={\fontsize{6 pt}{15.600000000000001 pt}\bfseries\sffamily}
    ] {1300};
    \node at (axis cs:1.5, -7.149388851079401) [,
    color=white,
    rotate=0.0,
    font={\fontsize{6 pt}{15.600000000000001 pt}\bfseries\sffamily}
    ] {36.00};
    \node at (axis cs:1.5, 12.999507953643798) [,
    color=white,
    rotate=0.0,
    font={\fontsize{6 pt}{15.600000000000001 pt}\bfseries\sffamily}
    ] {0.0};
    \node at (axis cs:1.5, 35.98755909282309) [,
    color=white,
    rotate=0.0,
    font={\fontsize{6 pt}{15.600000000000001 pt}\bfseries\sffamily}
    ] {2000};
    \node at (axis cs:2.5, 12.814011496543884) [,
    color=white,
    rotate=0.0,
    font={\fontsize{6 pt}{15.600000000000001 pt}\bfseries\sffamily}
    ] {0.0};
    \node at (axis cs:2.5, 35.55134726575725) [,
    color=white,
    rotate=0.0,
    font={\fontsize{6 pt}{15.600000000000001 pt}\bfseries\sffamily}
    ] {2000};
    \node at (axis cs:3.5, -7.149388851079401) [,
    color=white,
    rotate=0.0,
    font={\fontsize{6 pt}{15.600000000000001 pt}\bfseries\sffamily}
    ] {35.81};
    \node at (axis cs:3.5, 14.278156493804929) [,
    color=white,
    rotate=0.0,
    font={\fontsize{6 pt}{15.600000000000001 pt}\bfseries\sffamily}
    ] {0.0};
    \node at (axis cs:3.5, 37.19499024997084) [,
    color=white,
    rotate=0.0,
    font={\fontsize{6 pt}{15.600000000000001 pt}\bfseries\sffamily}
    ] {2000};
    \node at (axis cs:4.5, -7.149388851079401) [,
    color=white,
    rotate=0.0,
    font={\fontsize{6 pt}{15.600000000000001 pt}\bfseries\sffamily}
    ] {449.5};
    \node at (axis cs:4.5, 14.5557224849515) [,
    color=white,
    rotate=0.0,
    font={\fontsize{6 pt}{15.600000000000001 pt}\bfseries\sffamily}
    ] {0.0};
    \node at (axis cs:4.5, 37.766509153414674) [,
    color=white,
    rotate=0.0,
    font={\fontsize{6 pt}{15.600000000000001 pt}\bfseries\sffamily}
    ] {1600};
    \node at (axis cs:5.5, -7.149388851079401) [,
    color=white,
    rotate=0.0,
    font={\fontsize{6 pt}{15.600000000000001 pt}\bfseries\sffamily}
    ] {592.5};
    \node at (axis cs:5.5, 15.051322609038898) [,
    color=white,
    rotate=0.0,
    font={\fontsize{6 pt}{15.600000000000001 pt}\bfseries\sffamily}
    ] {0.0};
    \node at (axis cs:5.5, 38.850666161724234) [,
    color=white,
    rotate=0.0,
    font={\fontsize{6 pt}{15.600000000000001 pt}\bfseries\sffamily}
    ] {1400};
    \node at (axis cs:6.5, -7.149388851079401) [,
    color=white,
    rotate=0.0,
    font={\fontsize{6 pt}{15.600000000000001 pt}\bfseries\sffamily}
    ] {630.9};
    \node at (axis cs:6.5, 14.225435622692109) [,
    color=white,
    rotate=0.0,
    font={\fontsize{6 pt}{15.600000000000001 pt}\bfseries\sffamily}
    ] {0.0};
    \node at (axis cs:6.5, 40.117012817895734) [,
    color=white,
    rotate=0.0,
    font={\fontsize{6 pt}{15.600000000000001 pt}\bfseries\sffamily}
    ] {1400};
    \node at (axis cs:7.5, -7.149388851079401) [,
    color=white,
    rotate=0.0,
    font={\fontsize{6 pt}{15.600000000000001 pt}\bfseries\sffamily}
    ] {665.9};
    \node at (axis cs:7.5, 16.44196022891998) [,
    color=white,
    rotate=0.0,
    font={\fontsize{6 pt}{15.600000000000001 pt}\bfseries\sffamily}
    ] {0.0};
    \node at (axis cs:7.5, 49.18575900171447) [,
    color=white,
    rotate=0.0,
    font={\fontsize{6 pt}{15.600000000000001 pt}\bfseries\sffamily}
    ] {1400};
    \node at (axis cs:8.5, -7.149388851079401) [,
    color=white,
    rotate=0.0,
    font={\fontsize{6 pt}{15.600000000000001 pt}\bfseries\sffamily}
    ] {851.1};
    \node at (axis cs:8.5, 17.702141801834106) [,
    color=white,
    rotate=0.0,
    font={\fontsize{6 pt}{15.600000000000001 pt}\bfseries\sffamily}
    ] {0.0};
    \node at (axis cs:8.5, 54.78140354845931) [,
    color=white,
    rotate=0.0,
    font={\fontsize{6 pt}{15.600000000000001 pt}\bfseries\sffamily}
    ] {1200};
    \node at (axis cs:9.5, -7.149388851079401) [,
    color=white,
    rotate=0.0,
    font={\fontsize{6 pt}{15.600000000000001 pt}\bfseries\sffamily}
    ] {885.1};
    \node at (axis cs:9.5, 17.378364399909973) [,
    color=white,
    rotate=0.0,
    font={\fontsize{6 pt}{15.600000000000001 pt}\bfseries\sffamily}
    ] {0.0};
    \node at (axis cs:9.5, 52.68054388577406) [,
    color=white,
    rotate=0.0,
    font={\fontsize{6 pt}{15.600000000000001 pt}\bfseries\sffamily}
    ] {1100};
    \node at (axis cs:10.5, -7.149388851079401) [,
    color=white,
    rotate=0.0,
    font={\fontsize{6 pt}{15.600000000000001 pt}\bfseries\sffamily}
    ] {805.8};
    \node at (axis cs:10.5, 16.854729002952574) [,
    color=white,
    rotate=0.0,
    font={\fontsize{6 pt}{15.600000000000001 pt}\bfseries\sffamily}
    ] {0.0};
    \node at (axis cs:10.5, 49.24693943888761) [,
    color=white,
    rotate=0.0,
    font={\fontsize{6 pt}{15.600000000000001 pt}\bfseries\sffamily}
    ] {1200};
    \node at (axis cs:11.5, -7.149388851079401) [,
    color=white,
    rotate=0.0,
    font={\fontsize{6 pt}{15.600000000000001 pt}\bfseries\sffamily}
    ] {925.7};
    \node at (axis cs:11.5, 16.39220106315613) [,
    color=white,
    rotate=0.0,
    font={\fontsize{6 pt}{15.600000000000001 pt}\bfseries\sffamily}
    ] {0.0};
    \node at (axis cs:11.5, 47.245203398402026) [,
    color=white,
    rotate=0.0,
    font={\fontsize{6 pt}{15.600000000000001 pt}\bfseries\sffamily}
    ] {1100};
    \node at (axis cs:12.5, -7.149388851079401) [,
    color=white,
    rotate=0.0,
    font={\fontsize{6 pt}{15.600000000000001 pt}\bfseries\sffamily}
    ] {852.3};
    \node at (axis cs:12.5, 16.142237253189087) [,
    color=white,
    rotate=0.0,
    font={\fontsize{6 pt}{15.600000000000001 pt}\bfseries\sffamily}
    ] {0.0};
    \node at (axis cs:12.5, 46.48570445432283) [,
    color=white,
    rotate=0.0,
    font={\fontsize{6 pt}{15.600000000000001 pt}\bfseries\sffamily}
    ] {1200};
    \node at (axis cs:13.5, -7.149388851079401) [,
    color=white,
    rotate=0.0,
    font={\fontsize{6 pt}{15.600000000000001 pt}\bfseries\sffamily}
    ] {857.1};
    \node at (axis cs:13.5, 15.946892842292785) [,
    color=white,
    rotate=0.0,
    font={\fontsize{6 pt}{15.600000000000001 pt}\bfseries\sffamily}
    ] {0.0};
    \node at (axis cs:13.5, 45.72587321672157) [,
    color=white,
    rotate=0.0,
    font={\fontsize{6 pt}{15.600000000000001 pt}\bfseries\sffamily}
    ] {1200};
    \node at (axis cs:14.5, -7.149388851079401) [,
    color=white,
    rotate=0.0,
    font={\fontsize{6 pt}{15.600000000000001 pt}\bfseries\sffamily}
    ] {859.2};
    \node at (axis cs:14.5, 15.907354347229004) [,
    color=white,
    rotate=0.0,
    font={\fontsize{6 pt}{15.600000000000001 pt}\bfseries\sffamily}
    ] {0.0};
    \node at (axis cs:14.5, 45.867899029762356) [,
    color=white,
    rotate=0.0,
    font={\fontsize{6 pt}{15.600000000000001 pt}\bfseries\sffamily}
    ] {1200};
    \node at (axis cs:15.5, -7.149388851079401) [,
    color=white,
    rotate=0.0,
    font={\fontsize{6 pt}{15.600000000000001 pt}\bfseries\sffamily}
    ] {894.8};
    \node at (axis cs:15.5, 16.233450608730315) [,
    color=white,
    rotate=0.0,
    font={\fontsize{6 pt}{15.600000000000001 pt}\bfseries\sffamily}
    ] {0.0};
    \node at (axis cs:15.5, 47.79947883245581) [,
    color=white,
    rotate=0.0,
    font={\fontsize{6 pt}{15.600000000000001 pt}\bfseries\sffamily}
    ] {1100};
    \node at (axis cs:16.5, -7.149388851079401) [,
    color=white,
    rotate=0.0,
    font={\fontsize{6 pt}{15.600000000000001 pt}\bfseries\sffamily}
    ] {999.3};
    \node at (axis cs:16.5, 24.085998165024634) [,
    color=white,
    rotate=0.0,
    font={\fontsize{6 pt}{15.600000000000001 pt}\bfseries\sffamily}
    ] {0.0};
    \node at (axis cs:16.5, 61.8041999140511) [,
    color=white,
    rotate=0.0,
    font={\fontsize{6 pt}{15.600000000000001 pt}\bfseries\sffamily}
    ] {1000};
    \node at (axis cs:17.5, 11.643484156694953) [,
    color=white,
    rotate=0.0,
    font={\fontsize{6 pt}{15.600000000000001 pt}\bfseries\sffamily}
    ] {0.0};
    \node at (axis cs:17.5, 30.436357164469307) [,
    color=white,
    rotate=0.0,
    font={\fontsize{6 pt}{15.600000000000001 pt}\bfseries\sffamily}
    ] {973.2};
    \node at (axis cs:17.5, 59.03391256878691) [,
    color=white,
    rotate=0.0,
    font={\fontsize{6 pt}{15.600000000000001 pt}\bfseries\sffamily}
    ] {2000};
    \node at (axis cs:18.5, -7.149388851079401) [,
    color=white,
    rotate=0.0,
    font={\fontsize{6 pt}{15.600000000000001 pt}\bfseries\sffamily}
    ] {1000};
    \node at (axis cs:18.5, 17.59243702983856) [,
    color=white,
    rotate=0.0,
    font={\fontsize{6 pt}{15.600000000000001 pt}\bfseries\sffamily}
    ] {0.0};
    \node at (axis cs:18.5, 51.61637836372857) [,
    color=white,
    rotate=0.0,
    font={\fontsize{6 pt}{15.600000000000001 pt}\bfseries\sffamily}
    ] {1000};
    \node at (axis cs:19.5, -7.149388851079401) [,
    color=white,
    rotate=0.0,
    font={\fontsize{6 pt}{15.600000000000001 pt}\bfseries\sffamily}
    ] {986.7};
    \node at (axis cs:19.5, 15.981840569496155) [,
    color=white,
    rotate=0.0,
    font={\fontsize{6 pt}{15.600000000000001 pt}\bfseries\sffamily}
    ] {0.0};
    \node at (axis cs:19.5, 44.45709122810972) [,
    color=white,
    rotate=0.0,
    font={\fontsize{6 pt}{15.600000000000001 pt}\bfseries\sffamily}
    ] {1000};
    \node at (axis cs:20.5, -7.149388851079401) [,
    color=white,
    rotate=0.0,
    font={\fontsize{6 pt}{15.600000000000001 pt}\bfseries\sffamily}
    ] {986.7};
    \node at (axis cs:20.5, 14.981044017791747) [,
    color=white,
    rotate=0.0,
    font={\fontsize{6 pt}{15.600000000000001 pt}\bfseries\sffamily}
    ] {0.0};
    \node at (axis cs:20.5, 40.75827186058936) [,
    color=white,
    rotate=0.0,
    font={\fontsize{6 pt}{15.600000000000001 pt}\bfseries\sffamily}
    ] {1000};
    \node at (axis cs:21.5, -7.149388851079401) [,
    color=white,
    rotate=0.0,
    font={\fontsize{6 pt}{15.600000000000001 pt}\bfseries\sffamily}
    ] {1000};
    \node at (axis cs:21.5, 14.38581140756607) [,
    color=white,
    rotate=0.0,
    font={\fontsize{6 pt}{15.600000000000001 pt}\bfseries\sffamily}
    ] {0.0};
    \node at (axis cs:21.5, 39.045888464152206) [,
    color=white,
    rotate=0.0,
    font={\fontsize{6 pt}{15.600000000000001 pt}\bfseries\sffamily}
    ] {1000};
    \node at (axis cs:22.5, -7.149388851079401) [,
    color=white,
    rotate=0.0,
    font={\fontsize{6 pt}{15.600000000000001 pt}\bfseries\sffamily}
    ] {1000};
    \node at (axis cs:22.5, 13.903624738454818) [,
    color=white,
    rotate=0.0,
    font={\fontsize{6 pt}{15.600000000000001 pt}\bfseries\sffamily}
    ] {0.0};
    \node at (axis cs:22.5, 37.920044964765324) [,
    color=white,
    rotate=0.0,
    font={\fontsize{6 pt}{15.600000000000001 pt}\bfseries\sffamily}
    ] {1000};
    \node at (axis cs:23.5, -7.149388851079401) [,
    color=white,
    rotate=0.0,
    font={\fontsize{6 pt}{15.600000000000001 pt}\bfseries\sffamily}
    ] {1000};
    \node at (axis cs:23.5, 13.41814633512497) [,
    color=white,
    rotate=0.0,
    font={\fontsize{6 pt}{15.600000000000001 pt}\bfseries\sffamily}
    ] {0.0};
    \node at (axis cs:23.5, 36.87570051495497) [,
    color=white,
    rotate=0.0,
    font={\fontsize{6 pt}{15.600000000000001 pt}\bfseries\sffamily}
    ] {1000};
    \node at (axis cs:21.0, 85.0) [,
    color=white,
    rotate=0.0,
    font={\fontsize{8 pt}{15.600000000000001 pt}\bfseries\sffamily}
    ] {Independent Volumes [MWh]};
    \node at (axis cs:21.0, 95.0) [,
    color=white,
    rotate=0.0,
    font={\fontsize{8 pt}{15.600000000000001 pt}\bfseries\sffamily}
    ] {Dependent Volumes [MWh]};
  \end{axis}
  \begin{axis}[
    width = {\textwidth},
    height = {0.3\textheight},
    xlabel = {Hour},
    ylabel = {},
    title = {Maintenance schedule},
    unbounded coords=jump,
    scaled x ticks = false,
    xmajorgrids = true,
    xmin = 0,
    xmax = 24,
    xtick align = inside,
    xtick = {0.0,1.0,2.0,3.0,4.0,5.0,6.0,7.0,8.0,9.0,10.0,11.0,12.0,13.0,14.0,15.0,16.0,17.0,18.0,19.0,20.0,21.0,22.0,23.0,24.0},
    xticklabels = {$0$,$1$,$2$,$3$,$4$,$5$,$6$,$7$,$8$,$9$,$10$,$11$,$12$,$13$,$14$,$15$,$16$,$17$,$18$,$19$,$20$,$21$,$22$,$23$,$24$},
    xticklabel style = {font = {\fontsize{10 pt}{10.4 pt}\selectfont}, rotate = 0.0},
    x grid style = {color = kth-lightgray,
      line width = 0.25,
      solid},
    axis x line* = left,
    x axis line style = {line width = 1,
      solid},
    scaled y ticks = false,
    ylabel style = {font = {\fontsize{10 pt}{14.3 pt}\selectfont}, rotate = 0.0},
    ymajorgrids = true,
    ytick align = inside,
    ytick = {1, 2, 3, 4, 5, 6, 7, 8, 9, 10, 11, 12, 13, 14, 15},
    yticklabels = {Rebnis,Sadva,Bergnas,Slagnas,Bastusel,Grytfors,Gallejaur,Vargfors,Rengard,Batfors,Finnfors,Granfors,Krangfors,Selsfors,Kvistforsen},
    yticklabel style = {font = {\fontsize{8 pt}{10.4 pt}\selectfont}, rotate = 0.0},
    y grid style = {color = kth-lightgray,
      line width = 0.25,
      solid},axis y line* = left,y axis line style = {line width = 1,
      solid},
    xshift = 0mm,
    title style = {font = {\fontsize{14 pt}{18.2 pt}\selectfont}, rotate = 0.0},
    legend pos = {north west},
    legend style = {line width = 2, solid,font = {\fontsize{8 pt}{10.4 pt}\selectfont}},colorbar style={title=}]
    \addplot+ [line width = 2,
    solid,mark = none,
    color = black,
    mark size = 2.0,
    mark options = {
      color = black,
      fill = kth-blue,
      line width = 2,
      rotate = 0,
      solid
    },fill = kth-blue, forget plot]coordinates {
      (2,1)
      (3,1)
      (4,1)
      (5,1)
      (6,1)
      (7,1)
      (8,1)
      (9,1)
      (10,1)
    };
    \addplot+ [line width = 2,
    solid,mark = none,
    color = black,
    mark size = 2.0,
    mark options = {
      color = black,
      fill = kth-blue,
      line width = 2,
      rotate = 0,
      solid
    },fill = kth-blue, forget plot]coordinates {
      (2,2)
      (3,2)
      (4,2)
      (5,2)
      (6,2)
      (7,2)
      (8,2)
    };
    \addplot+ [line width = 2,
    solid,mark = none,
    color = black,
    mark size = 2.0,
    mark options = {
      color = black,
      fill = kth-blue,
      line width = 2,
      rotate = 0,
      solid
    },fill = kth-blue, forget plot]coordinates {
      (2,3)
      (3,3)
      (4,3)
      (5,3)
      (6,3)
      (7,3)
      (8,3)
      (9,3)
      (10,3)
    };
    \addplot+ [line width = 2,
    solid,mark = none,
    color = black,
    mark size = 2.0,
    mark options = {
      color = black,
      fill = kth-blue,
      line width = 2,
      rotate = 0,
      solid
    },fill = kth-blue, forget plot]coordinates {
      (2,4)
      (3,4)
      (4,4)
      (5,4)
      (6,4)
      (7,4)
    };
    \addplot+ [line width = 2,
    solid,mark = none,
    color = black,
    mark size = 2.0,
    mark options = {
      color = black,
      fill = kth-blue,
      line width = 2,
      rotate = 0,
      solid
    },fill = kth-blue, forget plot]coordinates {
      (2,5)
      (3,5)
      (4,5)
      (5,5)
      (6,5)
    };
    \addplot+ [line width = 2,
    solid,mark = none,
    color = black,
    mark size = 2.0,
    mark options = {
      color = black,
      fill = kth-blue,
      line width = 2,
      rotate = 0,
      solid
    },fill = kth-blue, forget plot]coordinates {
      (2,6)
      (3,6)
      (4,6)
      (5,6)
      (6,6)
    };
    \addplot+ [line width = 2,
    solid,mark = none,
    color = black,
    mark size = 2.0,
    mark options = {
      color = black,
      fill = kth-blue,
      line width = 2,
      rotate = 0,
      solid
    },fill = kth-blue, forget plot]coordinates {
      (2,7)
      (3,7)
      (4,7)
      (5,7)
    };
    \addplot+ [line width = 2,
    solid,mark = none,
    color = black,
    mark size = 2.0,
    mark options = {
      color = black,
      fill = kth-blue,
      line width = 2,
      rotate = 0,
      solid
    },fill = kth-blue, forget plot]coordinates {
      (2,8)
      (3,8)
      (4,8)
      (5,8)
    };
    \addplot+ [line width = 2,
    solid,mark = none,
    color = black,
    mark size = 2.0,
    mark options = {
      color = black,
      fill = kth-blue,
      line width = 2,
      rotate = 0,
      solid
    },fill = kth-blue, forget plot]coordinates {
      (3,9)
      (4,9)
      (5,9)
    };
    \addplot+ [line width = 2,
    solid,mark = none,
    color = black,
    mark size = 2.0,
    mark options = {
      color = black,
      fill = kth-blue,
      line width = 2,
      rotate = 0,
      solid
    },fill = kth-blue, forget plot]coordinates {
      (2,10)
      (3,10)
      (4,10)
    };
    \addplot+ [line width = 2,
    solid,mark = none,
    color = black,
    mark size = 2.0,
    mark options = {
      color = black,
      fill = kth-blue,
      line width = 2,
      rotate = 0,
      solid
    },fill = kth-blue, forget plot]coordinates {
      (2,11)
      (3,11)
      (4,11)
    };
    \addplot+ [line width = 2,
    solid,mark = none,
    color = black,
    mark size = 2.0,
    mark options = {
      color = black,
      fill = kth-blue,
      line width = 2,
      rotate = 0,
      solid
    },fill = kth-blue, forget plot]coordinates {
      (2,12)
      (3,12)
      (4,12)
    };
    \addplot+ [line width = 2,
    solid,mark = none,
    color = black,
    mark size = 2.0,
    mark options = {
      color = black,
      fill = kth-blue,
      line width = 2,
      rotate = 0,
      solid
    },fill = kth-blue, forget plot]coordinates {
      (2,13)
      (3,13)
      (4,13)
    };
    \addplot+ [line width = 2,
    solid,mark = none,
    color = black,
    mark size = 2.0,
    mark options = {
      color = black,
      fill = kth-blue,
      line width = 2,
      rotate = 0,
      solid
    },fill = kth-blue, forget plot]coordinates {
      (1,14)
      (2,14)
    };
    \addplot+ [line width = 2,
    solid,mark = none,
    color = black,
    mark size = 2.0,
    mark options = {
      color = black,
      fill = kth-blue,
      line width = 2,
      rotate = 0,
      solid
    },fill = kth-blue, forget plot]coordinates {
      (1,15)
      (2,15)
    };
  \end{axis}
\end{tikzpicture}

%% file: capacity_one_year_confidence_intervals.tex
\begin{tikzpicture}[]
  \begin{axis}[
    width = {\textwidth},
    height = {0.4\textheight},
    xlabel = {Number of Samples $N$},
    ylabel = {Confidence interval [MEur]},
    title = {Capacity expansion confidence intervals},
    unbounded coords=jump,
    scaled x ticks = false,
    xmajorgrids = true,
    xmin = -50,
    xtick = {0.0,100,500.0,1000.0,1500.0},
    xticklabels = {$0$,$100$,$500$,$1000$,$1500$},
    xtick align = inside,
    xticklabel style = {font = {\fontsize{10 pt}{10.4 pt}\selectfont}, rotate = 0.0},
    x grid style = {color = kth-lightgray,
      line width = 0.25,
      solid},
    axis x line* = left,
    x axis line style = {line width = 1,
      solid},
    scaled y ticks = false,
    ylabel style = {font = {\fontsize{10 pt}{14.3 pt}\selectfont}, rotate = 0.0},
    ymajorgrids = true,
    ytick align = inside,
    yticklabel style = {font = {\fontsize{10 pt}{10.4 pt}\selectfont}, rotate = 0.0},
    y grid style = {color = kth-lightgray,
      line width = 0.25,
      solid},axis y line* = left,y axis line style = {line width = 1,
      solid},
    xshift = 0mm,
    title style = {font = {\fontsize{14 pt}{18.2 pt}\selectfont}, rotate = 0.0},
    legend pos = {south east},
    legend style = {line width = 2, solid,font = {\fontsize{8 pt}{10.4 pt}\selectfont}},colorbar style={title=}]
    \addplot+[draw=none, color = kth-blue,
    line width = 0,
    solid,mark = *,
    mark size = 2.0,
    mark options = {
      color = black,
      fill = kth-blue,
      line width = 1,
      rotate = 0,
      solid
    }] coordinates {
      (10.0, 114.22)
      (50.0, 113.78)
      (100.0, 113.32)
      (500.0, 114.16)
    };
    \addlegendentry{VRP confidence intervals}
    \addplot+[draw=none, color = kth-red,
    line width = 0,
    solid,mark = *,
    mark size = 2.0,
    mark options = {
      color = black,
      fill = kth-red,
      line width = 1,
      rotate = 0,
      solid
    }] coordinates {
      (500.0, 112.34)
    };
    \addlegendentry{EEV confidence interval}
    \addplot+ [color = kth-blue,
    line width = 1,
    solid,mark = -,
    mark size = 2.0,
    mark options = {
      color = black,
      line width = 1,
      rotate = 0,
      solid
    },forget plot]coordinates {
      (10.0, 107.53)
      (10.0, 120.91)
    };
    \addplot+ [color = kth-blue,
    line width = 1,
    solid,mark = -,
    mark size = 2.0,
    mark options = {
      color = black,
      line width = 1,
      rotate = 0,
      solid
    },forget plot]coordinates {
      (50.0, 110.43)
      (50.0, 117.13)
    };
    \addplot+ [color = kth-blue,
    line width = 1,
    solid,mark = -,
    mark size = 2.0,
    mark options = {
      color = black,
      line width = 1,
      rotate = 0,
      solid
    },forget plot]coordinates {
      (100.0, 111.71)
      (100.0, 114.93)
    };
    \addplot+ [color = kth-blue,
    line width = 1,
    solid,mark = -,
    mark size = 2.0,
    mark options = {
      color = black,
      line width = 1,
      rotate = 0,
      solid
    },forget plot]coordinates {
      (500.0, 113.93)
      (500.0, 114.38)
    };
    \addplot+ [color = kth-red,
    line width = 1,
    solid,mark = -,
    mark size = 2.0,
    mark options = {
      color = black,
      line width = 1,
      rotate = 0,
      solid
    },forget plot]coordinates {
      (500.0, 111.26)
      (500.0, 113.43)
    };
  \end{axis}
\end{tikzpicture}

%% file: capacity_20_year_confidence_intervals.tex
\begin{tikzpicture}[]
  \begin{axis}[
    width = {\textwidth},
    height = {0.4\textheight},
    xlabel = {Number of Samples $N$},
    ylabel = {Confidence interval [MEur]},
    title = {Capacity expansion confidence intervals},
    unbounded coords=jump,
    scaled x ticks = false,
    xmajorgrids = true,
    xmin = -5,
    xtick = {0.0,50.0,100.0},
    xticklabels = {$0$,$50$,$100$},
    xtick align = inside,
    xticklabel style = {font = {\fontsize{10 pt}{10.4 pt}\selectfont}, rotate = 0.0},
    x grid style = {color = kth-lightgray,
      line width = 0.25,
      solid},
    axis x line* = left,
    x axis line style = {line width = 1,
      solid},
    scaled y ticks = false,
    ylabel style = {font = {\fontsize{10 pt}{14.3 pt}\selectfont}, rotate = 0.0},
    ymajorgrids = true,
    ytick align = inside,
    yticklabel style = {font = {\fontsize{10 pt}{10.4 pt}\selectfont}, rotate = 0.0},
    y grid style = {color = kth-lightgray,
      line width = 0.25,
      solid},axis y line* = left,y axis line style = {line width = 1,
      solid},
    xshift = 0mm,
    title style = {font = {\fontsize{14 pt}{18.2 pt}\selectfont}, rotate = 0.0},
    legend pos = {north east},
    legend style = {line width = 2, solid,font = {\fontsize{8 pt}{10.4 pt}\selectfont}},colorbar style={title=}]
    \addplot+[draw=none, color = kth-blue,
    line width = 0,
    solid,mark = *,
    mark size = 2.0,
    mark options = {
      color = black,
      fill = kth-blue,
      line width = 1,
      rotate = 0,
      solid
    }] coordinates {
      (10.0, 7728.42)
      (50.0, 7426.77)
      (100.0, 7424.45)
    };
    \addlegendentry{VRP confidence intervals}
    \addplot+[draw=none, color = kth-red,
    line width = 0,
    solid,mark = *,
    mark size = 2.0,
    mark options = {
      color = black,
      fill = kth-red,
      line width = 1,
      rotate = 0,
      solid
    }] coordinates {
      (100.0, 6785.56)
    };
    \addlegendentry{EEV confidence interval}
    \addplot+ [color = kth-blue,
    line width = 1,
    solid,mark = -,
    mark size = 2.0,
    mark options = {
      color = black,
      line width = 1,
      rotate = 0,
      solid
    },forget plot]coordinates {
      (10.0, 7247.46)
      (10.0, 8209.39)
    };
    \addplot+ [color = kth-blue,
    line width = 1,
    solid,mark = -,
    mark size = 2.0,
    mark options = {
      color = black,
      line width = 1,
      rotate = 0,
      solid
    },forget plot]coordinates {
      (50.0, 6981.44)
      (50.0, 7872.10)
    };
    \addplot+ [color = kth-blue,
    line width = 1,
    solid,mark = -,
    mark size = 2.0,
    mark options = {
      color = black,
      line width = 1,
      rotate = 0,
      solid
    },forget plot]coordinates {
      (100.0, 7194.93)
      (100.0, 7653.97)
    };
    \addplot+ [color = kth-red,
    line width = 1,
    solid,mark = -,
    mark size = 2.0,
    mark options = {
      color = black,
      line width = 1,
      rotate = 0,
      solid
    },forget plot]coordinates {
      (100.0, 6448.34)
      (100.0, 7122.78)
    };
  \end{axis}
\end{tikzpicture}